\documentclass[11pt,parskip=half]{article}

\usepackage{geometry}                
\usepackage[parfill]{parskip}    

\usepackage[utf8]{inputenc}

\usepackage{graphicx}

\usepackage{amsmath, amssymb, amsthm, amsfonts}
\usepackage{mathtools}
\usepackage{bbm} 

\mathtoolsset{showonlyrefs=true, showmanualtags=true}
\usepackage{hyperref}

\usepackage[svgnames]{xcolor}
\usepackage[breakable,theorems]{tcolorbox}
\usepackage[font=small,labelfont=bf]{caption}

\usepackage{libertine}
\usepackage[libertine]{newtxmath}
\usepackage{inconsolata}
\usepackage{microtype}

\usepackage{enumitem}

\usepackage{booktabs}

\usepackage{csquotes}

\usepackage{pgf}

 \usepackage[backend=biber,
             giveninits=true,
             mincitenames=1,
             maxcitenames=2,
             maxbibnames=99,
             uniquename=false,
             uniquelist=false,
             style=ext-authoryear,
             url = false,
             eprint = false,
             isbn = false,
             doi = false,
             articlein = false,
             dashed = false,
            ]{biblatex}

\numberwithin{equation}{section}

\newtcbtheorem[number within=section]%
  {definition}{Definition}{theorem style=plain, 
                           coltitle=green!30!black,
                           colframe=white,
                           colback=white,
                           breakable}{def}

\newtcbtheorem[use counter from=definition, number within=section]%
  {remark}{Remark}{theorem style=plain, 
                   coltitle=red!40!black,
                   colframe=white,
                   colback=white,
                   breakable}{rem}

\newtcbtheorem[use counter from=definition, number within=section]%
  {lemma}{Lemma}{theorem style=plain, 
                 coltitle=orange!45!black,
                 colframe=white,
                 colback=white,
                 breakable}{lem}

\newtcbtheorem[use counter from=definition, number within=section]%
  {theorem}{Theorem}{theorem style=plain, 
                     coltitle=orange!45!black,
                     colframe=white,
                     colback=white,
                     breakable}{thm}

\newtcbtheorem[use counter from=definition, number within=section]%
  {corollary}{Corollary}{theorem style=plain, 
                         coltitle=orange!45!black,
                         colframe=white,
                         colback=white,
                         breakable}{cor}

\newtcbtheorem[use counter from=definition, number within=section]%
  {proposition}{Proposition}{theorem style=plain, 
                             coltitle=orange!45!black,
                             colframe=white,
                             colback=white,
                             breakable}{prop}

\newtcbtheorem[use counter from=definition, number within=section]%
  {example}{Example}{theorem style=plain, 
                     coltitle=white!20!black,
                     colframe=white,
                     colback=white,
                     breakable}{ex}

\newcommand{\PP}{\mathcal{P}}
\newcommand{\PPP}[1]{\mathcal{P}_{\!#1}}
\newcommand{\CC}{\mathcal{C}}
\newcommand{\FF}{\mathcal{F}}

\newcommand{\RR}{\mathbb{R}}
\newcommand{\NN}{\mathbb{N}}
\newcommand{\EE}{\mathbb{E}}

\newcommand{\Lip}{\mathrm{Lip}}

\newcommand{\Normal}{\mathcal{N}}
\newcommand{\Prob}{\mathbb{P}}
\newcommand{\Exp}{\mathbb{E}}

\newcommand{\Unif}{\mathrm{Unif}}

\newcommand{\mutualinfo}{M}
\newcommand{\kl}[2]{D\left(#1\,|\,#2\right)}

\newcommand{\id}{\mathrm{id}}
\newcommand{\ind}{\mathbbm{1}}

\newcommand{\cmin}{{c_{XY}}}

\newcommand{\dif}{\mathrm{d}}
\newcommand{\supp}{\mathrm{supp}}

\newcommand{\weak}{\rightharpoonup}
\newcommand{\trace}{\mathrm{trace}}

\newcommand{\dcov}{\mathrm{dcov}}
\newcommand{\dcor}{\mathrm{dcor}}

\newcommand{\ot}{T}
\newcommand{\td}{\tau}
\newcommand{\tc}{\rho}

\newcommand{\covering}{\mathcal{N}}

\bibliography{sources}
\graphicspath{{figures/}}

\title{Transport Dependency: Optimal Transport Based Dependency Measures}

\author{ 
  Thomas Giacomo Nies\thanks{Institute for Mathematical Stochastics, University of Göttingen}
  \thanks{Cluster of Excellence: Multiscale Bioimaging (MBExC), University Medical Center, Göttingen} \\[-0.25ex]
  {\small \href{mailto:thomas.nies@uni-goettingen.de}{thomas.nies@uni-goettingen.de}} \vspace{1.2ex}\\
  Thomas Staudt\footnotemark[1] \footnotemark[2] \\[-0.25ex]
  {\small \href{mailto:thomas.staudt@uni-goettingen.de}{thomas.staudt@uni-goettingen.de}} \vspace{1.2ex}\\
  Axel Munk\footnotemark[1] \footnotemark[2] \,\thanks{Max Planck Institute for Biophysical Chemistry, Göttingen} \\[-0.25ex]
  {\small \href{mailto:munk@math.uni-goettingen.de}{munk@math.uni-goettingen.de}}
}

\date{\vspace{-3ex}}

\begin{document}

\maketitle

\begin{abstract}
  \noindent
  Finding meaningful ways to measure the statistical dependency between random
  variables $\xi$ and $\zeta$ is a timeless statistical endeavor. In recent
  years, several novel concepts, like the distance covariance, have extended
  classical notions of dependency to more general settings.
  In this article, we propose and study an alternative framework 
  that is based on optimal transport. The \emph{transport
  dependency} $\td \ge 0$ applies to general Polish spaces and intrinsically
  respects metric properties.
  For suitable ground costs, independence is fully characterized by $\td = 0$.
  Via proper normalization of $\td$, three \emph{transport correlations} $\tc_\alpha$,
  $\tc_\infty$, and $\tc_*$ with values in $[0, 1]$ are defined. They attain the
  value $1$ if and only if $\zeta = \varphi(\xi)$, where $\varphi$ is an
  $\alpha$-Lipschitz function for $\tc_\alpha$, a measurable function for
  $\tc_\infty$, or a multiple of an isometry for $\tc_*$.
  The transport dependency can be estimated consistently by an empirical
  plug-in approach, but alternative estimators with the same convergence rate
  but significantly reduced computational costs are also proposed. Numerical
  results suggest that $\td$ robustly recovers dependency between data sets with
  different internal metric structures. The usage for inferential tasks, like
  transport dependency based independence testing, is illustrated on a data set
  from a cancer study.
\end{abstract}



\section{Introduction}

\begin{figure}
  \centering
  {\footnotesize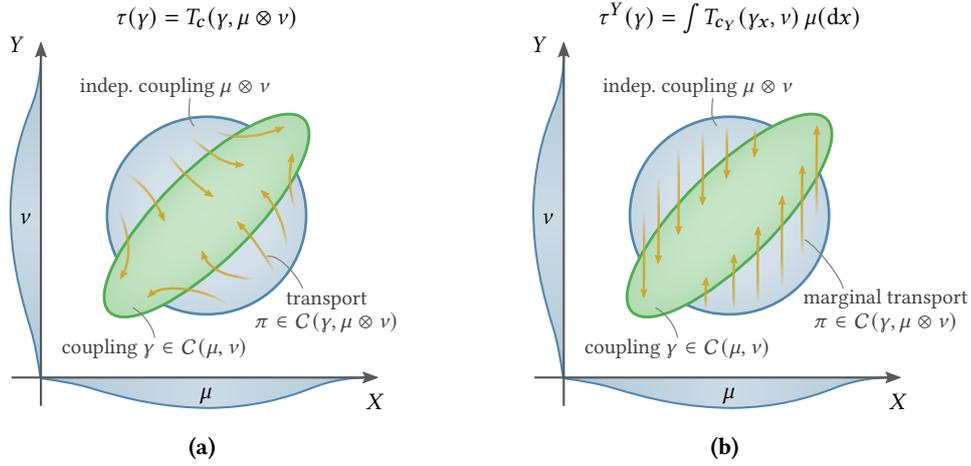}
  \caption{Illustration of the (marginal) transport dependency.
    Sketch \textbf{(a)} shows how mass could be moved in an optimal way
    when transforming $\mu\otimes\nu$ into $\gamma$ (visualized on level sets).
    The closer $\gamma$ is to $\mu\otimes \nu$, the less mass has to be
    transported and the smaller the value of $\td(\gamma)$ will be.
    Sketch \textbf{(b)} displays the same situation as (a), but this time
    transport along the space $X$ is forbidden and mass has to be transported
    vertically, which is clearly less optimal.
    }
  \label{fig:tdep-marginal-tdep}
\end{figure}

In this article, we explore a method to quantify the statistical dependence
between two random variables $\xi$ and $\zeta$ on Polish spaces $X$ and $Y$ via
optimal transport (see, e.g., \cite{rachev1998}, \cite{villani2008},
\cite{ambrosio2008}, \cite{santambrogio2015} for comprehensive analytical or
\cite{peyre2019} and \cite{panaretos2020} for computational and statistical
treatments). The core idea is to calculate the effort necessary to transform the
joint distribution $\gamma$ of $\xi$ and $\zeta$ into the product of their
marginal distributions $\mu$ and $\nu$. This motivates the definition of the
\emph{transport dependency} (Definition~\ref{def:transport-dependency})
\begin{equation}\label{eq:tdep-intro}
  \td(\xi, \zeta)
  =
  \td(\gamma)
  = \ot_c\big(\gamma, \mu\otimes\nu\big)
  = \inf_{\pi}\int\!\! c\,\dif \pi,
\end{equation}
where $\ot_c$ is the optimal transport cost with (non-negative) base costs $c$
on $X\times Y$. The infimum on the right is taken over the set of all couplings
between $\gamma$ and $\mu\otimes\nu$, i.e., probability distributions $\pi$ on
$(X\times Y)^2$ with marginals $\gamma$ and $\mu\otimes\nu$.
Figure~\ref{fig:tdep-marginal-tdep}a illustrates this concept. If the cost
function has benign properties, the transport dependency $\td$ displays many
traits that are attractive for a measure of statistical association. For
instance, if $c$ is a metric, then $\td(\gamma) = 0$ if and only if $\gamma
= \mu\otimes\nu$, which means statistical independence of $\xi$ and $\zeta$.

\paragraph{Prior work.}
The idea of evaluating an optimal transport cost between a coupling $\gamma$ and
the product $\mu\otimes\nu$ of its marginals has recently gained attention in the
statistical and machine learning literature. For example, \textcite{mori2020}
introduced the \emph{Earth mover's correlation}, a coefficient of dependency on
Polish metric spaces
that is based on a special case of \eqref{eq:tdep-intro}. As we will see later
on, several open conjectures of their work -- for example, the characterization of
couplings $\gamma$ with maximal Earth mover's correlation -- are resolved by our
theory.
In Euclidean settings, 
variants of \eqref{eq:tdep-intro} have been proposed under the names
\emph{Wasserstein dependence measure} and \emph{Wasserstein total correlation}
by \textcite{ozair2019} and \textcite{xiao2019}, who applied it to
beneficial effect in the context of representation learning. The same ansatz
also underlies recent work by \textcite{mordant2021}, who defined
\emph{Wasserstein dependency coefficients} that are powered by
\eqref{eq:tdep-intro} under squared Euclidean cost. For the purpose of
normalization, the authors divide $\tau(\gamma)$ by the supremum of
$\tau(\tilde\gamma)$ over all $\tilde\gamma$ with fixed marginals $\mu$ and
$\nu$. However, it is difficult to calculate these coefficients and
quasi-Gaussian surrogates are necessary for application.
In contrast, our approach relies on easily computable upper bounds that extend
those in \textcite{mori2020}.
Another instance of the transport dependency $\td$ has recently been explored by
\textcite{wiesel2021}, who proposed the 
association measure 
\begin{equation}\label{eq:marginal-tdep-intro}
  \td^Y(\xi, \zeta)
  =
  \td^Y(\gamma)
  =
  \int \ot_{c_Y}(\gamma_x, \nu)\,\mu(\dif x),
\end{equation}
where $(\gamma_x)_{x\in X}$ denotes the disintegration of $\gamma$ with respect
to the first coordinate (meaning that $\gamma_x$ is the law of $\zeta$ given
$\xi = x$) and $c_Y$ is the power of a metric on the space $Y$. The author
employed a similar upper bound as \textcite{mori2020} to derive a normalized
coefficient that
exhibits a number of desirable properties postulated by \textcite{dette2013} and
\textcite{chatterjee2020}.
Integrals of the form \eqref{eq:marginal-tdep-intro} have previously also
appeared in the context of generalization bounds for statistical learning
problems \parencite{zhang2018, lopez2018, wang2019}.
During our investigation of the transport dependency, it will become
clear that \eqref{eq:tdep-intro} and \eqref{eq:marginal-tdep-intro} are tightly
related.
In particular, $\td$ reduces to $\td^Y$ if transport in the space $X$ is
forbidden by costs that assume the value $\infty$ for non-vertical
movements (see Figure~\ref{fig:tdep-marginal-tdep}b). Moreover, $\td^Y$
naturally emerges both as a limit case and as an upper bound of $\td$.
We call $\td^Y$ the \emph{marginal transport dependency}
(Definition~\ref{def:marginal-transport-dependency}).

\begin{figure}
  \centering
  {\footnotesize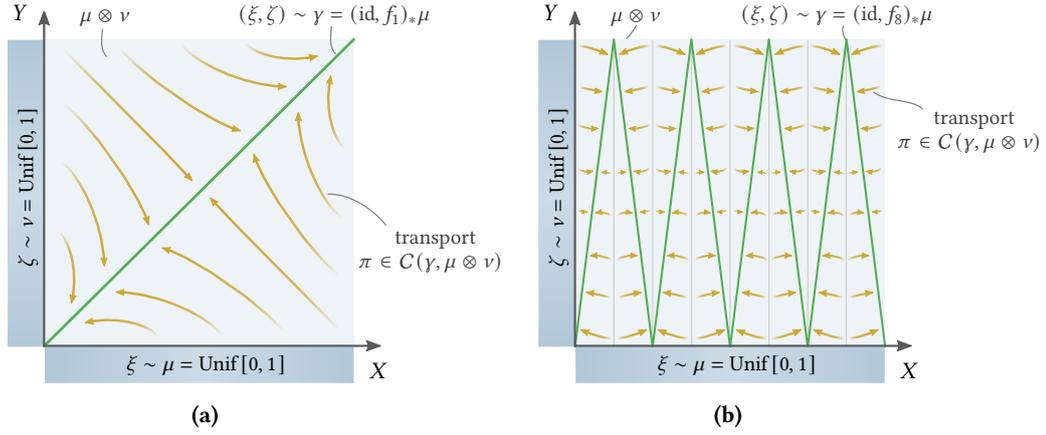}
  \caption{Marginal and joint distribution of random variables $\xi
    \sim \Unif[0, 1]$ and $\zeta = f_n(\xi) \sim \Unif[0, 1]$ for zigzag
    functions $f_n$ with $n$ linear segments for $n = 1$ in \textbf{(a)} and $n
    = 8$ in \textbf{(b)}. The drawn arrows illustrate how the optimal
    transport between $\mu\otimes\nu$ and $\gamma$ could look like. Note that
    the mutual information and related concepts that only measure the
    information content do not distinguish between the two scenarios. In
    contrast, dependency measures that are aware of metric or topological
    properties, like the transport dependency $\td$ or the distance
    covariance, assign a (much) lower degree of dependence to scenario (b). This
    discrepancy stresses an important point: should deterministic but chaotic
    relations between $\xi$ and $\zeta$ maximize a measure for statistical
    dependency? After all, one may not be able to recover the relation in
    practice and distinguish it from noise if data is limited.}
  \label{fig:zigzag}
\end{figure}

\paragraph{Mutual information.}
The general idea to compare the joint distribution of random variables to the
product of their marginals dates far back. In his landmark work,
\textcite{shannon1948} introduced the \emph{mutual information}
$\mutualinfo(\gamma) = \kl{\gamma}{\mu\otimes\nu}$, where $D$ denotes the
Kullback-Leibler divergence. The mutual information has since become an
indispensable tool for measuring the information content stored in the relation
between random variables and has found application in feature selection
\parencite{estevez2009}, image registration and alignment
\parencite{maes1997,pluim2000}, clustering \parencite{kraskov2005},
and independence testing \parencite{berrett2019}, besides others. Its immediate
use for statistical data analysis, however, is
complicated by several issues. For example, it is often inconvenient to estimate
$\mutualinfo(\gamma)$ from data, as density estimates or binning / clustering
methods are necessary and estimation suffers from the curse of dimensionality
\parencite{hall1993,paninski2008,berrett2019efficient}.
Furthermore, the mutual information does not respect topological or metric
properties of the coupling $\gamma$, as measurable rearrangements of $\xi$ and
$\zeta$ leave $M(\gamma)$ invariant. In this sense, it is not able to
distinguish \enquote{chaotic} relations between $\xi$ and $\zeta$ from
\enquote{well-behaved} ones (see Figure~\ref{fig:zigzag}). Despite these
potential drawbacks, the mutual information and its surrogates, such as the
mutual information dimension \parencite{sugiyama2013} or the maximal information
coefficient \parencite{reshef2011}, are widely used tools for detecting and
quantifying statistical dependency in data sets.

\paragraph{Distance covariance.}
A more recent approach to capture dependency 
by contrasting $\gamma$ to $\mu\otimes\nu$ is the \emph{distance covariance}
\parencite{szekely2007}.
The (Euclidean) distance covariance between
random vectors $\xi$ in $\RR^r$ and $\zeta$ in $\RR^q$ for $r,q\in\NN$ is
a weighted $L_2$ distance between the joint characteristic function $f_{\gamma}$
of $\gamma$ and the product of the marginal characteristic functions $f_{\mu}$
and $f_{\nu}$,
\begin{equation}\label{eq:euclidean-dcov}
  \dcov^2(\xi, \zeta)
  =
  \frac{1}{c_rc_q}
  \int \frac{\big|f_{\gamma}(t, s) - f_\mu(t)f_\nu(s)\big|^2}{\|t\|^{1+r}\|s\|^{1+q}}
  \,\lambda^r(\dif t)\,\lambda^q(\dif s),
\end{equation}
where $\|\cdot\|$ is the Euclidean norm, $c_d$ is a constant only depending on
the dimensions $d$, and $\lambda^d$ denotes the Lebesgue measure in $\RR^d$.
\textcite{lyons2013} 
later proposed
a generalization of \eqref{eq:euclidean-dcov} to (separable) metric spaces $(X,
d_X)$ and $(Y, d_Y)$, given by
\begin{equation}\label{eq:dcov}
\begin{aligned}
  \dcov^2(\xi, \zeta)
  &=
  \Exp\big[d_X(\xi, \xi')\,d_Y(\zeta, \zeta')\big]
  + \Exp\big[d_X(\xi, \xi')\big]\,\Exp\big[d_Y(\zeta, \zeta')\big] \\
  &\quad - 2\,\Exp\big[d_X(\xi, \xi')\,d_Y(\zeta, \zeta'')\big],
\end{aligned}
\end{equation}
where $(\xi', \zeta') \sim \gamma$ and $\zeta'' \sim \nu$ are independent copies
of $\xi$ and $\zeta$.
If $X$ and $Y$ are of strong negative
type\footnote{%
  A separable metric space $(X, d_X)$ is of \emph{negative type} iff there is an
  isometric embedding $\varphi$ of $\smash{\big(X, d_X^{1/2}\big)}$ into
  a separable Hilbert space. This condition asserts $\dcov^2 \ge 0$.
  If the mean embedding $\mu\mapsto \int\varphi\,\dif\mu$ for probability
  measures $\mu\in\PP(X)$ with finite first $d_X$-moment is additionally
  injective, the space $(X, d_X)$ is of \emph{strong negative type} and $\dcov
  = 0$ characterizes independence. Examples for spaces of negative type are
  $L_p$ spaces for $1 \le p \le 2$, ultrametric spaces, and weighted trees
  \parencite[Theorem 3.6]{meckes2013}. Known counter examples are
  $\smash{\RR^d}$ with $l_p$ norms for $p > 2$ (see the references in
  \cite{lyons2013}).
},
$\dcov^2(\xi, \zeta)$ as defined above is indeed non-negative and vanishes if
and only if $\xi$ and $\zeta$ are independent \parencite{lyons2013,
jakobsen2017}. Due to fast computability on data, a performant unbiased
estimator \parencite{gao2021asymptotic}, and a well-understood limit theory, it
is a compelling instrument for non-parametric independence
testing \parencite{yao2016,castro-prado2020,chakraborty2019} and related
problems, like independent component analysis \parencite{matteson2017}.


The distance covariance possesses the natural upper bound $\dcov^2(\xi, \zeta)
\le \dcov(\xi, \xi)\cdot\dcov(\zeta, \zeta)$, which is utilized to define the
normalized \emph{distance correlation} with values in $[0, 1]$. The distance
correlation $\dcor$, which can serve as a more general surrogate for classical
dependency coefficients like the Pearson correlation, has the following
properties \parencite{lyons2013}:
\begin{itemize}[itemindent=0pt]
  \item $\dcor(\xi, \zeta) = 0$ iff $\xi$ and $\zeta$ are
    independent,
  \item $\dcor(\xi, \zeta) = 1$ iff there is a $\beta > 0$
    and an isometry $\varphi\colon (X, \beta d_X) \to (Y, d_Y)$ with
    $\zeta = \varphi(\xi)$.
\end{itemize}
This provides a neat and tangible interpretation:
the distance
correlation measures a degree of \emph{isometric} functional dependency (up to
scalings). For the Euclidean case, this means that a value of $\dcor(\xi,
\zeta) = 1$ is assumed if and only if $\zeta = \beta\,(A\xi + a)$, where $A$ is
an orthogonal matrix, $a$ an offset vector, and $\beta > 0$ (at least if the
support of $\mu$ contains an open set). Other relations between $\xi$ and
$\zeta$, even if they are deterministic, result in smaller values $\dcor(\xi,
\zeta) < 1$. Indeed, the more chaotic the relation becomes, the further away one
is from an isometric dependency, and the lower the value of the distance
correlation will typically be (see Figure~\ref{fig:zigzag}). This draws a sharp
distinction to other (non-parametric) concepts of dependency, like the mutual
information or several recently proposed coefficients of association
\parencite{dette2013,chatterjee2020,deb2020,wiesel2021}, which assume maximal
values for \emph{any measurable} deterministic relation -- not only for
structured ones.


\paragraph{Other measures of association.}
Looking past classical concepts of correlation that are universally
applied (like the Pearson correlation, Spearman's $\rho$, or Kendall's $\tau$),
the literature on how to best measure dependency quickly becomes immensely broad
and scattered. Besides the approaches cited above, we exemplarily mention maximal
correlation coefficients \parencite{gebelein1941,koyak1987}, rank or copula
based methods \parencite{schweizer1981,marti2017}, or various measures acting on
the distribution of pairwise distances \parencite{friedman1983, heller2013}.
For a survey, see \textcite{tjostheim2018}.
More closely related to our work, optimal transport maps have recently been
utilized to define multivariate rank statistics that allow for asymptotically
distribution-free independence tests
\parencite{ghosal2019,shi2020distribution,shi2021}. Furthermore, optimal
transport induced geometries have been explored for covariance analysis in
functional data analysis \parencite{petersen2019,dubey2020}. 

A different class of data analysis and exploration techniques to be mentioned in
this context are those that quantify how multiple data sets are \emph{spatially}
associated. A prominent example is Ripley's $K$ function \parencite{ripley1976},
for which new developments have recently been advanced \parencite{amgad2015}.
A particular case of spatial association arises in colocalization problems in
cell microscopy. We mention \textcite{wang2017}, who propose a colocalization
metric based on Kendall's $\tau$, and \textcite{tameling2021}, who suggest
certain surrogates of the optimal transport plan for quantifying colocalization.

\paragraph{Transport dependency.}
A primary reason for the widely scattered literature on this topic is that the
notion of \enquote{dependency} eludes the reduction to a single real number and
heavily depends on the context.
One important demarcation line in this regard has already been stressed: do we
aim to measure dependency in a purely stochastic sense (like the mutual
information), or do we also seek to impose structural conditions, like linearity
(Pearson correlation), monotonicity (rank correlations), or metric compatibility
(distance correlation)?
Indeed, the theme of \emph{shape restrictions} is central for recent efforts
to find meaningful quantifiers of dependency (\cite{cao2020}; see also
\cite{guntuboyina2018} for related work on shape-restricted regression).

\begin{figure}
  \centering
  {\small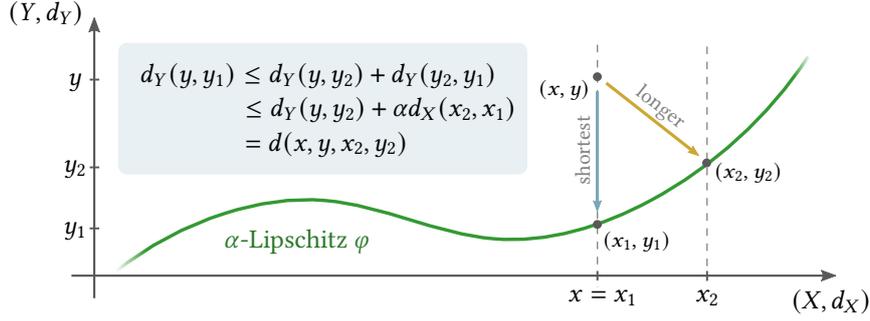}
  \caption{Projection onto the graph of an $\alpha$-Lipschitz function
    $\varphi\colon X \to Y$ under costs $c = d^p$ as in \eqref{eq:cost-intro}.
    By applying the triangle inequality of $d_Y$ and the Lipschitz property
    $d_Y\big(\varphi(x_1), \varphi(x_2)\big) \le \alpha d_X(x_1, x_2)$, one can see
    that the cost-optimal way to move any point $(x, y)\in X\times Y$ to the
    graph of $\varphi$ is to simply shift it up- or downwards.
  }
  \label{fig:contractions}
\end{figure}

In this article, we contribute to this topic by establishing the transport
dependency as a principled tool that  is flexible enough to bridge the gap
between unstructured and structured dependency quantification. To begin with,
$\td$ combines a number of attractive general properties that are desirable for
a measure of association (see Section~\ref{sec:general-properties}).
For example, the condition $\td(\gamma) = 0$
fully characterizes independence under mild assumptions on the costs
(Theorem~\ref{thm:independence}).
Furthermore, the value of $\td(\gamma)$ only relies on the intrinsic cost
structure and does not change, say, under transformations of $\xi$ and $\zeta$
that leave the cost function invariant (Proposition~\ref{prop:invariance}). The
transport dependency also behaves well under perturbations: it is
(Lipschitz) continuous (Proposition~\ref{prop:continuity} and
Theorem~\ref{thm:continuity-metric-power}) and additive
independent noise contributions (i.e., convolutions) can only decrease the value of $\td$
(Theorem~\ref{thm:convolution}). Moreover, if $t$ percent of $\gamma$ are
replaced (contaminated) by a distribution $\tilde\gamma$ with the same marginals,
then $\td\big((1-t) \gamma + t\tilde\gamma\big) \le (1-t)\td(\gamma)
+ t \td(\tilde\gamma)$ (Proposition~\ref{prop:contamination}).
A particularly interesting theory unfolds on Polish metric spaces $(X, d_X)$ and
$(Y, d_Y)$ under additive costs
\begin{equation}\label{eq:cost-intro}
  c(x_1, y_1, x_2, y_2)
  =
  \big(\alpha\cdot d_X(x_1, x_2) + d_Y(y_1, y_2)\big)^p
\end{equation}
for $\alpha, p > 0$, where $x_1, x_2 \in X$ and $y_1, y_2 \in Y$ (see
Section~\ref{sec:contractions}).
In this case, one of our core insights
(Theorem~\ref{thm:contractions-deterministic}) establishes that $\td(\gamma)$
attains the upper bound (Proposition~\ref{prop:upper-bound})
\begin{equation}\label{eq:upper-bound-intro}
  \td(\gamma)
  \le
  \int \!d_Y^p\,\dif(\nu\otimes\nu)
\end{equation}
if and only if $\gamma$ is concentrated on the graph of an $\alpha$-Lipschitz
function $\varphi\colon X \to Y$.
In other words, the equality $\td(\xi, \zeta) = \Exp[d_Y(\zeta, \zeta')^p]$, for
$\zeta, \zeta' \sim \nu$ independent, holds if and only if $\zeta
= \varphi(\xi)$, where $\varphi$ is an $\alpha$-Lipschitz map.
The underlying intuition is captured in Figure~\ref{fig:contractions}, which
shows that the cost of moving any point $(x, y)$ in $X\times Y$ to the graph of
an $\alpha$-Lipschitz function $\varphi$ is minimized by vertical movements.
Transport along vertical movements, however, leads to bound
\eqref{eq:upper-bound-intro}.

\paragraph{Transport correlation.}
Based on the cost function \eqref{eq:cost-intro} for $\alpha > 0$, we propose
a family of coefficients $\tc_\alpha \in [0, 1]$ that are tuned to detect
$\alpha$-Lipschitz associations (see Section~\ref{sec:coefficients}). The
$\alpha$\emph{-transport correlation}
(Definition~\ref{def:transport-correlation}) is given by
\begin{equation}
  \tc_\alpha(\xi, \zeta)
  =
  \tc_\alpha(\gamma)
  =
  \left(
  \frac{\td(\gamma)}{\int\! d_Y^p \,\dif(\nu\otimes\nu)}
  \right)^{1/p},
\end{equation}
where the scaling by the inverse of bound \eqref{eq:upper-bound-intro}, which we
assume to be positive, guarantees $0 \le \tc_\alpha \le 1$. Two of the
hallmark features of $\tc_\alpha$ are (Proposition~\ref{prop:transport-correlation})
\begin{itemize}
  \item $\tc_\alpha(\gamma) = 0$ iff $\xi$ and $\zeta$ are independent,
  \item $\tc_\alpha(\gamma) = 1$ iff there is an $\alpha$-Lipschitz function
    $\varphi\colon (X, d_X) \to (Y, d_Y)$ with $\zeta = \varphi(\xi)$.
\end{itemize}
This allows us to view $\tc_\alpha$ as a generalized alternative to Pearsons's
correlation coefficient that measures the degree of association in terms of best
approximation by $\alpha$-Lipschitz functions (instead of linear functions).
Later in the article, we will see that the idea behind $\tc_\alpha$ can fluently
be extended to the limit $\alpha \to \infty$
(Theorem~\ref{thm:marginal-transport-dependence}), in which the transport
dependency $\td$ becomes the marginal transport dependency $\td^Y$.
The resulting coefficient $\tc_\infty$, which we name \emph{marginal transport
correlation} (Definition~\ref{def:marginal-transport-dependency}), satisfies
(Proposition~\ref{prop:marginal-transport-correlation})
\begin{itemize}
  \item $\tc_\infty(\gamma) = 0$ iff $\xi$ and $\zeta$ are independent,
  \item $\tc_\infty(\gamma) = 1$ iff there is a measurable function
    $\varphi\colon X \to Y$ with $\zeta = \varphi(\xi)$.
\end{itemize}
Note that the marginal transport correlation is equal to the measure of
association introduced by \textcite{wiesel2021}, who already recognized the
above properties.

Both the $\alpha$-transport correlation and the marginal transport correlation
are asymmetric concepts and only measure to what extent $\zeta$ can
be understood as a function of $\xi$, not vice versa. If symmetry is
desired, the coefficients can be adjusted in various ways.
For example, we choose
$\alpha_*^p = \int d_Y^p\,\dif(\nu\otimes\nu) / \int d_X^p\,\dif(\mu\otimes\mu)$
and set $\tc_* = \tc_{\alpha_*}$ to define the \emph{isometric transport
correlation} (Definition~\ref{def:isometric-transport-correlation}), which has
the following properties
(Proposition~\ref{prop:isometric-transport-correlation}):
\begin{itemize}
  \item $\tc_*(\gamma) = 0$ iff $\xi$ and $\zeta$ are independent,
  \item $\tc_*(\gamma) = 1$ iff there is a $\beta > 0$ and an isometry
    $\varphi\colon (X, \beta d_X) \to (Y, d_Y)$ with $\zeta = \varphi(\xi)$.
\end{itemize}
This shows that $\tc_*$ assumes its extremal values for exactly the same
$\gamma$ as the distance correlation. For dependencies leading to non-extremal
values, on the other hand, $\tc_*$ can differ considerably from $\dcor$ (see our
simulation study in Section~\ref{sec:applications}). Contrary to the Earth
Movers correlation proposed by \textcite{mori2020}, which is also a symmetric
coefficient based on $\td$ (see equation
\eqref{eq:contracting-transport-correlation} in Section~\ref{sec:coefficients}),
$\tc_*$ satisfies all of the axioms put forward by \textcite{mori2019}.

\begin{figure}
  \centering
  \includegraphics[width=0.85\textwidth]{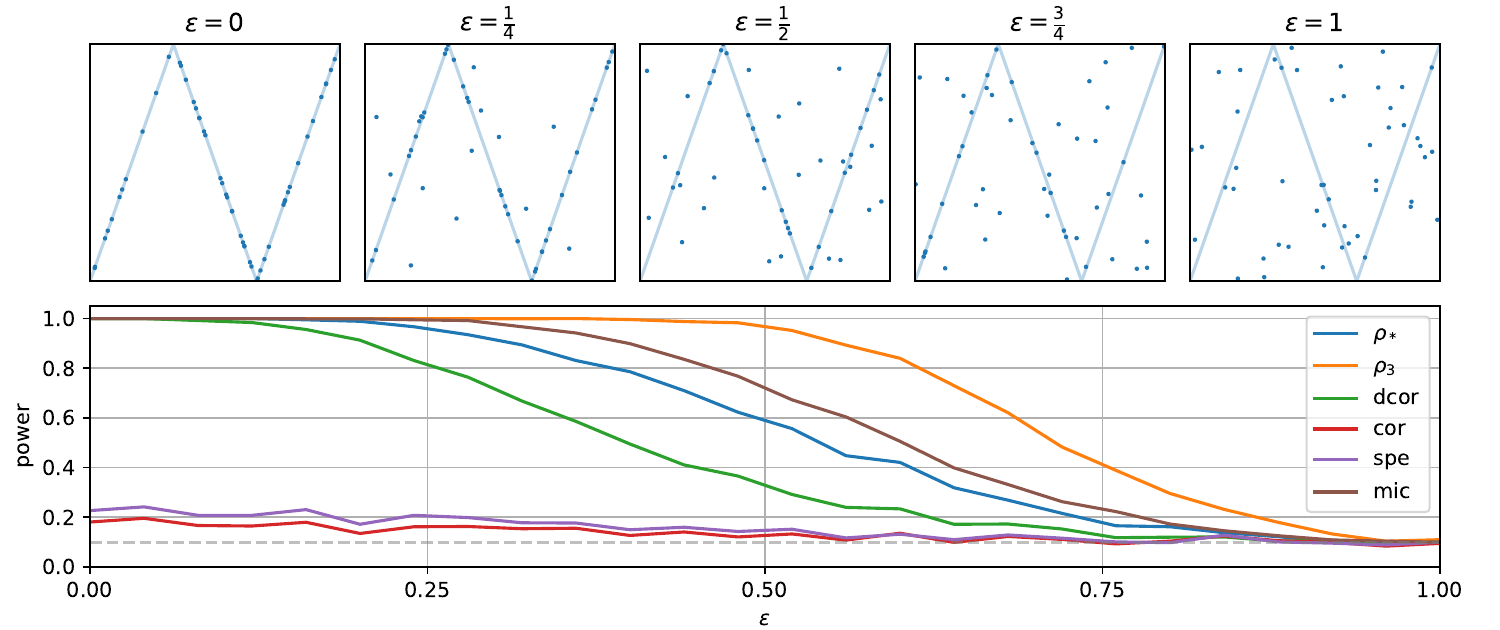}
  \caption{Power of independence tests under increasing levels of noise. The
    upper row depicts exemplary i.i.d.\ samples of size $n = 50$ at varying
    contamination noise levels $\epsilon\in[0, 1]$ with $\xi\sim\Unif[0,1]$. The
    undisturbed samples ($\epsilon = 0$) are drawn from the graph of
    a $3$-Lipschitz zigzag function. The figure on the bottom depicts the power
    of level $0.1$ permutation tests for independence that are based on $\tc_*$,
    $\tc_\alpha$ for $\alpha = 3$, the distance correlation, the Pearson
    correlation, the Spearman correlation, and the maximal information
    coefficient. See Section~\ref{sec:applications} for more details on the
    numerical setting and the applied permutation tests.}
  \label{fig:power-intro}
\end{figure}

\paragraph{Estimation, computation, and application.}
In the presence of empirical data, the transport dependency can be estimated
consistently by the plug-in estimator $\td(\hat\gamma_n)$, where $\hat\gamma_n$
is the empirical measure of $n$ independent and identically distributed
observations (see Section~\ref{sec:estimation}). Under mild conditions, the
convergence rate of this estimator is determined by the intrinsic dimension of
$\gamma$ (equation \eqref{eq:lca-metric}, Corollary~\ref{cor:lca-manifold}, and
Theorem~\ref{thm:lcaproduct}) and not by the potentially higher dimension of
$\mu\otimes\nu$ or even $X\times Y$.
This phenomenon of \emph{lower complexity
adaptation} (LCA) was recently uncovered and described by
\textcite{hundrieser2022empirical}, extending findings by
\textcite{weed2019sharp, divol2022measure}.
Consequently, the stronger the measures $\mu$ and $\nu$ are coupled, the more
LCA affects the empirical convergence rates.
Since the calculation of
$\td(\hat\gamma_n)$ requires solving an optimal transport problem of size
$n\times n^2$, which causes a high computational burden in case of large amounts
of data, we also propose re-sampling based estimators (with effective size
$n\times n$) that can be calculated by several orders of magnitude faster (for
large $n$) but achieve the same LCA rates (Theorem~\ref{thm:lcasampling}). Even
though Monte-Carlo simulations (see
Section~\ref{sec:applications}) demonstrate a sizable bias in higher dimensions,
which is inherited from the underlying optimal transport problem, the transport
dependency generally performs well in the task of recognizing dependency via
permutation tests.
In particular, the most distinctive feature of $\td$ is its flexibility to adapt
to different environments by the freedom to work with arbitrary cost functions.
In simulated settings, this is exemplified by considerable performance
advantages of $\rho_\alpha$ compared to other dependency coefficients when
detecting
(noisy) $\alpha$-Lipschitz signals (see Figure~\ref{fig:power-intro} for an
example, or Figure~\ref{fig:adapt} in Section~\ref{sec:applications}). In
real-world applications, this flexibility promises to be most beneficial
if tailor-made criteria for similarity are available.
We illustrate this on a high-dimensional gene expression dataset, where
biologically motivated correlation scores will be integrated into the cost function
(Section~\ref{sec:application-to-gene-expression-data}). Our transport
dependency based analysis is able to reproduce the core conclusions previously
obtained by much more specialized methods \parencite{behr2020}. Crucially, the
involved cost structure lacks the properties of a metric, which is required by
other dependency coefficients that are applicable in high-dimensional settings.

\section{Preliminaries}
\label{sec:definitions}

This section indroduces the notation for the rest of the manuscript and
summarizes basic definitions and results from the theory of optimal
transport on Polish (i.e., separable and completely metrizable
topological) spaces.

\paragraph{Notation.}
The set of all probability measures on a Polish space $X$ is denoted by
$\PP(X)$. The integration of a (Borel-)measurable function $f\colon X \to \RR$
with respect to $\mu\in\PP(Y)$ is flexibly written as $\int f(x) \,\mu(\dif x)$,
$\int f \,\dif\mu$, or simply $\mu f$.
Every measurable function $\varphi\colon X \to Y$ between
Polish spaces induces a pushforward map $\varphi_\#\colon \PP(X) \to
\PP(Y)$ via $\mu \mapsto \varphi_\#\mu = \mu\circ \varphi^{-1}$.
We write $\mu_n \rightharpoonup \mu$ if a sequence
$(\mu_n)_{n\in\NN}\subset\PP(X)$ converges weakly (or in distribution) to $\mu$
in $\PP(X)$. 
The product $X\times Y$ of two Polish spaces is again Polish if equipped with
its product topology. We write $p^X$ and $p^Y$ for the Cartesian projections
onto the spaces $X$ and $Y$, and we let
\begin{equation}
  (\varphi, \psi)(x)
  =
  \big(\varphi(x), \psi(x)\big)
  \qquad\text{and}\qquad
  (\varphi\otimes\psi)(x, y)
  =
  \big(\varphi(x), \psi(y)\big)
\end{equation}
for suitable functions $\varphi$ and $\psi$. The notation $\CC(\mu, \nu) \subset
\PP(X\times Y)$ is used to denote the set of couplings (or transport plans)
between measures $\mu\in\PP(X)$ and $\nu\in\PP(Y)$. Thus, $\gamma\in\CC(\mu,
\nu)$ is a joint distribution on $X\times Y$ with marginals $p^X_\#\gamma = \mu$
and $p^Y_\#\gamma = \nu$. The product measure $\mu\otimes\nu$ is always an
element of $\CC(\mu, \nu)$.
We also write $\CC(\mu, \cdot)$ or $\CC(\cdot, \nu)$ for the subsets of
$\PP(X\times Y)$ where only one marginal distribution is specified. If we want
to condition $\gamma\in\CC(\mu, \nu)$ on one of its components, we use the
shorthand notation
\begin{equation}
  \gamma(\dif x, \dif y)
  =
  \gamma(x, \dif y)\,\mu(\dif x)
\end{equation}
to indicate a disintegration of $\gamma$ along the space $X$, where $\gamma(x,
\cdot)\in\PP(Y)$ is a probability distribution for each $x\in X$ (see
\cite{chang1997,kallenberg2006}). The family
$\big(\gamma(x,\cdot)\big)_{x\in X}$ is a probability kernel (or Markov kernel
or stochastic kernel), which means that the mapping $x \mapsto \gamma(x, A)$ is
measurable for each Borel set $A\subset Y$. Analog notation will be used for
conditioning on $y\in Y$ or if products of more than two spaces are considered.

Whenever convenient, we may express our arguments in terms of random elements
$\xi$ and $\zeta$ instead of probability measures. Their joint law is usually
$(\xi, \zeta)\sim\gamma\in\CC(\mu, \nu)$. Note that $\gamma(x, \cdot)$ is the
conditional distribution of $\zeta$ given $\xi = x$ and vice versa for
$\gamma(\cdot, y)$, and that the pushforward $f_\#\gamma$ equals the law of
$f(\xi, \zeta)$ whenever $f\colon X\times Y \to Z$ is a measurable map into
a Polish space.

We call a lower semi-continuous function $c\colon X\times X \to
[0,\infty]$ a \emph{cost function} on $X$ if it is symmetric and
vanishes on the diagonal, meaning $c(x_1, x_2) = c(x_2, x_1)$ and $c(x, x) = 0$
for all $x,x_1,x_2\in X$. Sometimes we require $c(x_1, x_2) > 0$ for $x_1
\neq x_2$, in which case we call the cost function \emph{positive}.

\paragraph{Optimal transport.}
The \emph{optimal transport cost} $\ot_c$ between measures $\mu$ and $\nu$ in
$\PP(X)$ for base costs $c$ on a Polish space $X$ is defined as
\begin{equation}\label{eq:ot-definition}
  \ot_c(\mu, \nu)
  = \inf_{\pi\in\CC(\mu, \nu)} \pi c.
\end{equation}
An \emph{optimal transport plan} $\pi^*$ that attains the infimum in equation
\eqref{eq:ot-definition} always exists if $c$ is lower semi-continuous
(\cite[Theorem 4.1]{villani2008}). In general, optimal transport plans need not
be unique.
If $c = d^p$ for a metric $d$ on $X$ with $p \ge 1$, the quantity $\ot_c(\mu,
\nu)^{1/p}$ is called the \emph{$p$-Wasserstein distance} and is a metric on the
probability measures on $X$ that have finite $p$-th moments
(\cite[Theorem~6.9]{villani2008}).

When a transport plan $\pi\in\CC(\mu, \nu)$ is concentrated on the graph of
a function $\varphi\colon X\to X$, this function is called \emph{transport map}
and satisfies $\varphi_\#\mu = \nu$ and $(\id,\varphi)_\# \mu = \pi$. Under
certain conditions on $\mu$, $\nu$, and $c$, optimal transport plans $\pi^*$
that minimize \eqref{eq:ot-definition} correspond to optimal transport maps
$\varphi^*$. For example, this always holds when $\mu$ has a Lebesgue density in
$\RR^d$ for $d\in\NN$ and $c$ is given by an $l_p$ norm with $p > 1$
(\cite[Theorem~3.7]{gangbo1996geometry}).

The optimal transport problem \eqref{eq:ot-definition} can alternatively
be stated in its dual formulation,
\begin{equation}\label{eq:ot-definition-dual}
  \ot_c(\mu, \nu) = \sup_{f \oplus g \le c} \mu f + \nu g,
\end{equation}
where $(f\oplus g)(x_1, x_2) = f(x_1) + g(x_2)$ and where the supremum is taken
over bounded and continuous functions $f$ and $g$. This fact is commonly known
as Kantorovich-duality, dating back to \textcite{kantorovich1942}.
Like for the existence of transport plans, the cost function $c$ being lower
semi-continuous is sufficient for \eqref{eq:ot-definition} and
\eqref{eq:ot-definition-dual} to coincide \parencite[Theorem 5.10]{villani2008}.
If $c$ is a metric, \eqref{eq:ot-definition-dual} takes the particular form
\begin{equation}\label{eq:ot-definition-dual-metric}
  \ot_c(\mu, \nu) = \sup_{f\in\Lip_1(X)} \mu f - \nu f,
\end{equation}
where $\Lip_1(X)$ denotes all real valued 1-Lipschitz functions on $X$ with
respect to $c$. In particular, this is a special case of an integral probability
metric \parencite{muller1997, sriperumbudur2012}.

\section{Transport dependency: general properties}
\label{sec:general-properties}

The transport dependency features several desirable traits for a measure of
statistical association. In this section, we formally define it and discuss some
of its generic properties, including convexity, symmetry, continuity, and the
behavior under convolutions. We also establish three distinct upper bounds, one
of which corresponds to a special case of $\td$ when movement of mass along the
space $X$ is forbidden. This naturally leads to the definition of the marginal
transport dependency later on. Most of the proofs in this section are delegated
to Appendix~\ref{app:proofs}.

\begin{definition}{transport dependency}{transport-dependency}
Let $X$ and $Y$ be Polish spaces and $c$ a cost function on
$X\times Y$.
The ($c$-)\emph{transport dependency} $\td_c\colon\PP(X\times Y) \to
[0, \infty]$ is defined via
\begin{equation}\label{eq:transport-dependence-def}
  \td(\gamma)
  =
  \td_c(\gamma)
  =
  \ot_c\big(\gamma, p^X_\#\gamma\otimes p^Y_\#\gamma\big).
\end{equation}
\end{definition}
We usually omit the cost function $c$ in the subscript if it is apparent from
the context, and we write $\td(\xi, \zeta) = \td(\gamma)$ for
random elements $\xi$ and $\zeta$ with joint distribution $\gamma$.
One can easily see that $\td(\xi, \zeta) = 0$ in case of statistical
independence. If $c$ is a positive cost function, this criterion is even
sufficient for independence of $\xi$ and $\zeta$.

\begin{theorem}{independence}{independence}
  Let $X$ and $Y$ be Polish spaces and $c$ a positive cost function on $X\times
  Y$. Then $\td(\gamma) = 0$ if and only if $\gamma = \mu\otimes\nu$ for some
  $\mu\in\PP(X)$ and $\nu\in\PP(Y)$.
\end{theorem}
\vspace{-1em}
\begin{proof}
  If $\gamma = \mu\otimes\nu$, the transport plan $\pi
  = (\id,\id)_\#(\mu\otimes\nu)\in\CC(\gamma, \mu\otimes\nu)$ satisfies $0 \le
  \td(\gamma) \le \pi c = 0$, where we used that $c \ge 0$ and that $c$
  vanishes on the diagonal. Conversely, if $\td(\gamma) = 0$ for some
  $\gamma\in\CC(\mu, \nu)$, we find $\pi^* c = 0$ for the optimal plan $\pi^*$,
  so $c = 0$ holds $\pi^*$-almost surely. Since $c$ is positive, this implies
  $\pi^*\{(x, y, x, y)\,|\, x\in X, y\in Y\} = 1$. Consequently, it follows that
  $\pi^* = (\id,\id)_\#(\mu\otimes\nu)\in\CC(\mu\otimes\nu, \mu\otimes\nu)$,
  where $\id$ denotes the identity map on $X\times Y$. This establishes $\gamma
  = \mu\otimes\nu$.
\end{proof}

\begin{figure}
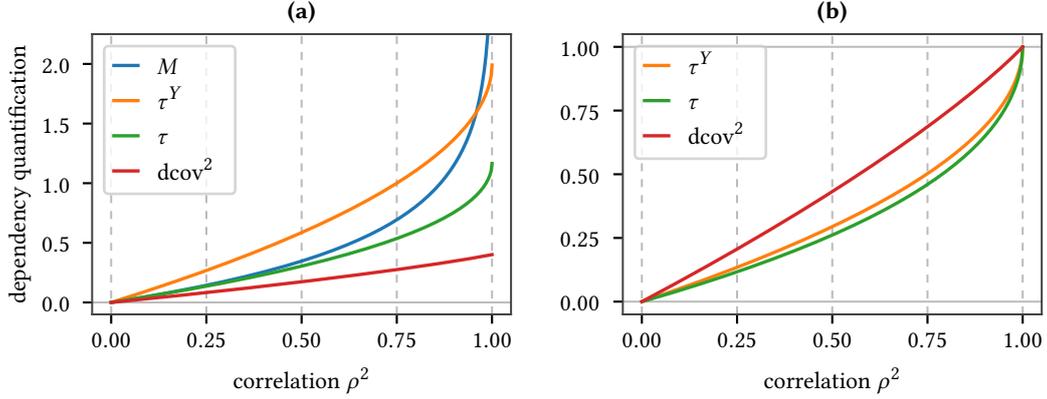

  \centering
  \input{figures/gauss-dependency-measures-absolute.pgf}
  \input{figures/gauss-dependency-measures-relative.pgf}
  \caption{Selected dependency measures in the bivariate case
    $\xi,\zeta\sim \Normal(0, 1)$ with $\mathrm{cov}(\xi, \zeta) = \rho$, see
    Example~\ref{ex:gauss}. Graph \textbf{(a)} shows the mutual
    information (which diverges as $\rho^2 \to 1$), the Euclidean distance
    covariance \eqref{eq:euclidean-dcov}, the marginal transport dependency
    \eqref{eq:marginal-tdep-intro}, and the general transport dependency as
    function of $\rho^2$. Graph \textbf{(b)} shows the latter three quantities
    normalized such that their respective maximal value equals $1$ for $\rho^2
    = 1$. See Example~\ref{ex:gauss-marginal} for a closed form of the marginal
    transport dependency $\td^Y$ in this setting.}
  \label{fig:gauss-dependency-measures}
\end{figure}

\begin{example}{multivariate Gaussian}{gauss}
  Let $(\xi, \zeta)$ be a pair of real valued and normally distributed random
  vectors on $X\times Y=\RR^{r+q}$ that follow a joint distribution
  $\gamma=\mathcal{N}(\eta, \Sigma)$ so that
  \begin{equation*}
	  \eta =
 			\begin{pmatrix}
 				\eta_1   \\
 				\eta_2  
		\end{pmatrix}, \quad \Sigma=
 			\begin{pmatrix}
				\Sigma_{11} & \Sigma_{12}  \\
				\Sigma_{21} & \Sigma_{22} 
 			\end{pmatrix}
  \end{equation*}
  with $\eta_1\in\RR^{r}, \eta_2\in \RR^{q}$, $\Sigma_{11}\in \RR^{r\times r},
  \Sigma_{22}\in \RR^{q\times q}$, $\Sigma_{12}\in \RR^{r\times q}$ and
  $\Sigma_{21} = \Sigma_{12}^T$. 
  Since independence of normal random variables is characterized by zero correlation,
  the independent coupling is given by
  \begin{equation}
    \mu\otimes\nu = \mathcal{N}(\eta, \Sigma_{\text{ind}})
    \quad\text{with}\quad
    \Sigma_{\text{ind}} =
    \begin{pmatrix}
			\Sigma_{11} & 0  \\
			0 & \Sigma_{22} 
    \end{pmatrix}.
  \end{equation}
  As costs on the space $X\times Y$, we consider the squared Euclidean distance
  \begin{equation*}
    c(x_1, y_1, x_2, y_2) = \|x_1-x_2\|^2 +\|y_1- y_2\|^2
  \end{equation*}
  for $(x_1, y_1), (x_2, y_2)\in\RR^{r}\times \RR^{q}$.
  In this setting, evaluating the transport dependence of $\gamma$ corresponds
  to computing the square of the 2-Wasserstein distance between two normal
  distributions. We obtain \parencite{dowson1982}
  \begin{equation}\label{eq:tdep-normal}
    \td(\gamma)
    =
    2\,\mathrm{trace}(\Sigma_{11})+2\,\mathrm{trace}(\Sigma_{22})
    - 2\, \mathrm{trace}\left( 
    \begin{pmatrix}
      \Sigma_{11}^2 & \Sigma_{11} \Sigma_{12}  \\
      \Sigma_{22} \Sigma_{21} & \Sigma_{22}^2 
    \end{pmatrix}^{1/2}\right).
  \end{equation}
  If $\xi$ and $\zeta$ are univariate random variables, a more explicit formula
  can easily be derived. Let the covariance matrix be given by 
  \begin{equation*}
		\Sigma=
 			\begin{pmatrix}
 				\sigma_{1}^{2} & \rho \sigma_1 \sigma_2  \\
 				\rho \sigma_{1} \sigma_{2} & \sigma_{2}^2 
 			\end{pmatrix}
  \end{equation*}
  with  $\rho\in [0, 1)$ and $\sigma_1, \sigma_2 > 0$. 
  Then, using an identity provided in \parencite{levinger1980}, we can compute the
  square root of the $2\times 2$ matrix appearing in \eqref{eq:tdep-normal}
  and obtain
  \begin{equation}\label{eq:gauss2}
    \td(\gamma)
    =
    2\,\left(\sigma_{1}^2 + \sigma_{2}^2
     - \sqrt{\sigma_1^4+\sigma_2^4+2\sigma_1^2\sigma_2^2\sqrt{1 - \rho^2}}\right).
  \end{equation}
  If $\sigma_1 = \sigma_2 = \sigma$ for some $\sigma > 0$, this expression
  simplifies to
  \begin{equation}
    \td(\gamma)
    =
    2\sigma^2\left(2 - \sqrt{2 + 2\sqrt{1-\rho^2}}\right).
  \end{equation}
  As to be expected, the transport dependency between $\xi$ and $\zeta$ is
  a strictly increasing function of the correlation $\rho^2$. Its minimal value
  is $0$ for $\rho = 0$ and its maximal value is $(4 - 2\sqrt{2})\,\sigma^2
  \approx 1.2\sigma^2$ for $\rho = \pm 1$ (if $\sigma_1 = \sigma_2 = \sigma$).
  The mutual information is given by $M(\gamma) = - \log(1-\rho^2)/2$
  \parencite{gelfand1959} and the Euclidean distance covariance by
  \parencite{szekely2007}
  \begin{equation}
    \dcov^2(\gamma)
    =
    \frac{4\sigma^2}{\pi}\left(
      \rho \arcsin \rho + \sqrt{1-\rho^2} - \rho \arcsin(\rho/2) - \sqrt{4-\rho^2} + 1
    \right).
  \end{equation}
  See Figure~\ref{fig:gauss-dependency-measures} for a comparison.
\end{example}

\paragraph{Convexity and invariance.}
The optimal transport cost $\ot_c$ in \eqref{eq:ot-definition} is a convex
functional in both of its arguments \parencite[Theorem~4.8]{villani2008}.
Similarly, the transport dependency $\td$ can be shown to be convex on subsets
of $\PP(X\times Y)$ that share (at least) one marginal.

\begin{proposition}{convexity}{convexity}
  Let $X$ and $Y$ be Polish spaces and $c$ a cost function on $X\times Y$. Fix
  marginal distributions $\mu\in\PP(X)$ and $\nu\in\PP(Y)$. Then the transport
  dependency $\td$ is convex when restricted to $\CC(\mu, \cdot)$ or $\CC(\cdot,
  \nu)$.
\end{proposition}

Note that $\td$ is not convex on its whole domain, save for trivial exceptions.
In fact, convexity on all of $\PP(X\times Y)$ would imply $0 \le
2\,\td\big(\delta_{z_1}/2 + \delta_{z_2}/2\big) \le \td(\delta_{z_1})
+ \td(\delta_{z_2}) = 0$ for point masses $\delta_{z_i}$ with $z_1, z_2\in
X\times Y$. By induction, $\td$ would have to vanish on all empirical measures.
For reasonable costs, like metrics on $X \times Y$, this only happens if $X$ and
$Y$ are singletons.

Proposition~\ref{prop:convexity} can be used to show that replacing a part of
$\gamma$ with independent contributions consistently decreases the transport
dependency.

\begin{proposition}{convex contamination}{contamination}
 Let $X$ and $Y$ be Polish spaces and $c$ a cost function on $X \times Y$. For
 $\gamma\in\CC(\mu, \nu)$ with $\mu\in\PP(X)$ and $\nu\in\PP(Y)$ define
 $\gamma_t = (1-t)\gamma + t(\mu\otimes\nu)$. Then the mapping $t\mapsto\td(\gamma_t)$
 is convex, monotonically decreasing, and satisfies 
 \begin{equation}\label{eq:contamination}
   \td(\gamma_t) \le (1-t)\,\td(\gamma).
 \end{equation}
 When $c$ is a metric, \eqref{eq:contamination} holds with equality.
\end{proposition}

We next formulate a symmetry result for $\td$ that is based on a fundamental
invariance property of optimal transport (Lemma~\ref{lem:ot-invariance} in
Appendix~\ref{app:proofs}). It can easily be generalized to cost preserving maps
between distinct Polish spaces, but we restrict to the presented setting for
simplicity.

\begin{proposition}{invariance}{invariance}
  Let $X$ and $Y$ be Polish spaces and $c$ a cost function on $X\times Y$
  of the form
  \begin{equation}\label{eq:cost-invariance-condition}
    c(x_1, y_1, x_2, y_2) = h\big(c_X(x_1, x_2), c_Y(y_1, y_2)\big)
  \end{equation}
  for marginal costs $c_X$ and $c_Y$ on $X$ and $Y$ and a measurable function
  $h\colon[0,\infty]^2\to[0,\infty]$. If $f_X\colon X\to X$ and $f_Y\colon Y\to
  Y$ are measurable maps that leave $c_X$ and $c_Y$ invariant, then any
  $\gamma\in\PP(X\times Y)$ satisfies
  \begin{equation}\label{eq:invariance}
    \td_c(f_\#\gamma) = \td_{c}(\gamma).
  \end{equation}
\end{proposition}
\vspace{-0.5em}

If $c_X$ and $c_Y$ are (based on) metrics and $\xi$ and $\zeta$ are random
elements on $X$ and $Y$, we conclude that applying isometries to either $\xi$ or
$\zeta$ does not change the value $\td(\xi, \zeta)$. In other words, $\td$
measures dependency in a way that is indifferent to isometric transformations on
the margins. This property is shared by the distance covariance \eqref{eq:dcov}.

\begin{example}{invariance in Euclidean space}{translation}
  Let $\xi$ and $\zeta$ be random vectors in $X = \RR^r$ and $Y = \RR^q$. If
  $c_X$ and $c_Y$ denote the respective Euclidean metrics and $c$ takes the
  form $h(c_X, c_Y)$, for example $c = c_X^2 + c_Y^2$ as in
  Example~\ref{ex:gauss}, Proposition~\ref{prop:invariance} asserts that
  \begin{equation}
    \td(\xi, \zeta) = \td(A\,\xi + a, B\,\zeta + b)
  \end{equation}
  for all orthogonal matrices $A\in\RR^{r\times r}$, $B\in\RR^{q\times q}$ as
  well as vectors $a\in\RR^r$, $b\in\RR^q$.
\end{example}

\paragraph{Continuity and convergence.}
The optimal transport cost $\ot_c$ with lower semi-continuous base costs $c$ is
again lower semi-continuous with respect to the weak convergence of measures.
This property carries over to the transport dependency.
\begin{proposition}{lower semi-continuity}{lower-semi-continuity}
  Let $X$ and $Y$ be Polish spaces and let $c$ be a cost function on $X\times
  Y$. Then $\td$ is lower semi-continuous with respect to weak convergence.
\end{proposition}

To show proper continuity of $\td$, we need stronger assumptions, since the
optimal transport cost $\ot_c$ is in general not continuous under weak
convergence (\cite[Proposition 7.4]{santambrogio2015}). To this end, we
equip the Polish space $X$ with a compatible metric $d_X$ that completely
metrizes its topology. For $p \ge 1$, we define the set of probability
distributions with finite $p$-th moment by
\begin{equation}
  \PPP{p}(X)
  =
  \big\{ \mu\in\PP(X)\,\big|\,
         \mu\,d_X^p(\cdot, x_0) < \infty~\text{for some}~x_0\in X \big\},
\end{equation}
and we say that a sequence $(\mu_n)_{n\in\NN}$ in $\PPP{p}(X)$ converges
$p$-weakly to $\mu\in\PPP{p}(X)$ if
\begin{equation}\label{eq:p-weak-convergence}
  \mu_n\weak\mu
  \qquad\text{and}\qquad
  \mu_n\,d_X(\cdot, x_0)^p \to \mu\, d_X(\cdot, x_0)^p
\end{equation}
as $n\to\infty$ for some $x_0\in X$. It is a well known fact
\parencite{mallows1972} that the $p$-Wasserstein distance metrizes this
particular form of convergence (see \cite[Theorem~6.9]{villani2008} for
a general proof), so \eqref{eq:p-weak-convergence} is equivalent to
$\smash{\ot_c(\mu_n, \mu) \to 0}$ as $n\to\infty$ with $c = d_X^p$. Note that
the anchor point $x_0$ in these definitions does not matter and can be replaced
by any other element of $X$.

On products $X\times Y$ of two Polish metric spaces $(X, d_X)$ and $(Y, d_Y)$, there
are many different ways to choose a compatible metric. For simplicity, we pick
\begin{equation}\label{eq:product-metric-sum}
  d(x_1, y_1, x_2, y_2)
  =
  d_X(x_1, x_2) + d_Y(y_1, y_2)
\end{equation}
for $(x_1, y_1), (x_2, y_2)\in X\times Y$ in the following statement, even
though any equivalent metric, like $d = \max(d_X, d_Y)$, works as well. 

\begin{proposition}{continuity}{continuity}
  Let $(X, d_X)$ and $(Y, d_Y)$ be Polish metric spaces and let $c$ be
  a continuous cost function on $X\times Y$ bounded by $c \le d^p$ for $p
  \ge 1$ and $d$ as in \eqref{eq:product-metric-sum}.
  If the sequence $(\gamma_n)_{n\in\NN}$ converges $p$-weakly to $\gamma$ in
  $\PPP{p}(X\times Y)$, then $\td(\gamma_n) \to \td(\gamma)$.
\end{proposition}

An important consequence of Proposition~\ref{prop:continuity} is the guarantee
of consistency when $\td(\gamma)$ is empirically estimated (see
Section~\ref{sec:estimation} for more details on estimation). We also present
a continuity statement with explicit bounds in case that $c$ is \emph{equal} to
the power of a metric. This time, since we only rely on the triangle inequality,
the metric does not have to metrize the Polish topology. In fact, it suffices if
$d$ is a pseudo-metric on $X\times Y$.

\begin{theorem}{continuity}{continuity-metric-power}
  Let $X$ and $Y$ be Polish spaces and let $c = d^p$ for a continuous
  pseudo-metric $d$ on $X\times Y$. Assume that $\gamma,\gamma'\in\PP(X\times
  Y)$ with $\td(\gamma),\td(\gamma') < \infty$. Then, for any $p \ge 1$,
  \begin{equation}
    \big| \td(\gamma)^{1/p} - \td(\gamma')^{1/p} \big|
    \le
    \ot_{c}(\gamma, \gamma')^{1/p} + \ot_{c}(\mu\otimes\nu, \mu'\otimes\nu')^{1/p}.
  \end{equation}
\end{theorem}

If the metric $d$ on $X\times Y$ in the statement above is given in terms of
marginal metrics $d_X$ and $d_Y$ on $X$ and $Y$, like $d = d_X + d_Y$, one can
apply the triangle inequality to derive
\begin{align}
  \ot_c(\mu\otimes\nu, \mu'\otimes\nu')^{1/p}
  &\le
  \ot_c(\mu\otimes\nu, \mu'\otimes\nu)^{1/p}
  + \ot_c(\mu'\otimes\nu, \mu'\otimes\nu')^{1/p} \\
  &\le
  \ot_{c_X}(\mu, \mu')^{1/p} + \ot_{c_Y}(\nu, \nu')^{1/p} \\
  &\le
  2\,\ot_{c}(\gamma, \gamma')^{1/p}, \label{eq:td-lip}
\end{align}
where $\smash{c_X = d_X^p}$ and $\smash{c_Y = d_Y^p}$.
In this case, Theorem~\ref{thm:continuity-metric-power} implies that
$\td^{1/p}$ is 3-Lipschitz continuous with respect to the $p$-Wasserstein
distance: for any $\gamma_1,\gamma_2\in\PP(X\times Y)$, it holds that
\begin{equation}\label{eq:td-lipschitz}
  \big|\td(\gamma_1)^{1/p} - \td(\gamma_2)^{1/p}\big|
  \leq
  3\,\ot_{c}(\gamma_1, \gamma_2)^{1/p}.
\end{equation}

\paragraph{Convolutions.}

The sum of independent random variables plays a distinguished role in many
applications. 
If $\xi\sim\mu\in\PP(X)$ and $\epsilon_X\sim\kappa_X\in\PP(X)$ is an independent
noise variable in a vector space $X$, the distribution of $\xi + \epsilon_X$
equals the convolution $\mu * \kappa_X$. The latter is defined by
\begin{equation}
  (\mu * \kappa_X)(A)
  =
  \int \ind_{A}(x_1 + x_2)\,(\mu\otimes\kappa_X)(\dif x_1, \dif x_2)
\end{equation}
for any Borel set $A\subset X$. If $\zeta\sim\nu\in\PP(Y)$ is also
contaminated by an independent additive noise contribution $\epsilon_Y
\sim\kappa_Y\in\PP(Y)$, the joint distribution of $(\xi + \epsilon_X, \zeta
+ \epsilon_Y)$ is given by $\gamma * \kappa$ for $\kappa = \kappa_X \otimes
\kappa_Y$.

The following theorem provides insights into how the transport dependency of
$\gamma * \kappa$ is related to the value $\td(\gamma)$ for translation
invariant costs.
We work in Polish vector spaces, by which we mean topological vector spaces that
are Polish; examples include separable Banach spaces.

\begin{theorem}{convolution}{convolution}
  Let $X$ and $Y$ be Polish vector spaces and let $c$ be a cost function on
  $X\times Y$ that satisfies $c(x_1, y_1, x_2, y_2) = h(x_1 - x_2, y_1 - y_2)$
  for $(x_1,y_1), (x_2, y_2) \in X\times Y$ and some $h\colon X\times Y\to[0,
  \infty]$. For any $\gamma\in\PP(X\times Y)$ and $\kappa
  = \kappa_X\otimes\kappa_Y$ with $\kappa_X\in\PP(X)$ and $\kappa_Y\in\PP(Y)$,
  it holds that
  \begin{equation}
    \td(\gamma * \kappa)
    \le
    \td(\gamma).
  \end{equation}
  If additionally $c = d^p$ and $\td(\gamma) < \infty$ for a continuous
  pseudo-metric $d$ on $X\times Y$ with $p \ge 1$, then
  \begin{equation}
    \td(\gamma)^{1/p} - \td(\gamma * \kappa)^{1/p} \le 2 (\kappa h)^{1/p}.
  \end{equation}
\end{theorem}

Theorem~\ref{thm:convolution} unveils a fundamental property of the transport
dependency: if a coupling $\gamma$ is blurred by the convolution with a product
kernel, then $\td$ never increases (and typically decreases). Intuitively, this
is a desirable trait. For Gaussian noise, for example, it guarantees that $\td$
monotonically decreases with increasing standard deviation of the noise.

\begin{example}{Gaussian additive noise}{convolution-gauss}
  Let $X = \RR^r$ and $Y = \RR^q$ and consider the squared Euclidean cost of the
  form $c(x_1, y_1, x_2, y_2) = h(x_1 - x_2, y_1 - y_2) = \|x_1 - x_2\|^2
  + \|y_1 - y_2\|^2$, where $x_1, x_2\in X$ and $y_1, y_2 \in Y$.
  Under Gaussian noise contributions
  $\kappa_X \sim \Normal(0, \Sigma_X)$ and $\kappa_Y \sim \Normal(0, \Sigma_Y)$,
  where $\Sigma_X$ and $\Sigma_Y$ are covariance matrices in $\RR^{r\times r}$
  and $\RR^{q\times q}$, Theorem~\ref{thm:convolution} states that
  \begin{equation}
    \td(\gamma)^{1/2}
    \ge
    \td(\gamma * \kappa)^{1/2}
    \ge
    \td(\gamma)^{1/2} - 2\big(\trace\,\Sigma_X + \trace\,\Sigma_Y\big)^{1/2}
  \end{equation}
  for any $\gamma\in\PP(X \times Y)$ with $\td(\gamma)< \infty$.
\end{example}

\begin{example}{Manhattan-type costs}{convolution-norm}
  If $\big(X, \|\cdot\|_X\big)$ and $\big(Y, \|\cdot\|_Y\big)$ are separable
  Banach spaces and the costs $c$ are given by
  $h = \|\cdot\|_X + \|\cdot\|_Y$, application of
  Theorem~\ref{thm:convolution} yields
  \begin{equation}
    \td(\gamma)
    \ge
    \td(\gamma * \kappa)
    \ge
    \td(\gamma) -
    2 \left( \int\!\| x\|_X\,\kappa_X(\dif x)
    + \int\! \| y\|_Y\,\kappa_Y(\dif y)\right)
  \end{equation}
  for any $\gamma\in\PP(X\times Y)$, $\kappa_X\in\PP(X)$, and
  $\kappa_Y\in\PP(Y)$, as long as $\td(\gamma) < \infty$ holds.
\end{example}

\begin{figure}
  \centering
  {\footnotesize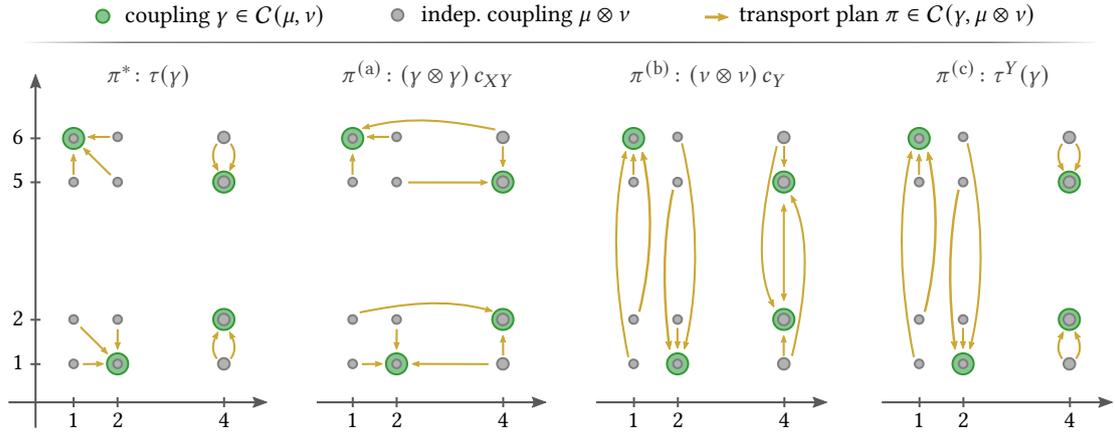}
  \caption{Different transport plans between $\gamma\in\CC(\mu\otimes\nu)$
    and $\mu\otimes\nu$ for $\gamma = \Unif\{(1, 6), (2, 1), (4, 2), (4, 5)\}$.
    Each arrow corresponds to the movement of $1/16$ parts of mass.
    The first graph shows the optimal transport plan $\pi^*$ under $l_2$ (or
    $l_1$) costs on $X \times Y = \RR^2$, while the plans $\pi^\mathrm{(a)}$ to
    $\pi^\mathrm{(c)}$ correspond to the plans behind the three
    upper bounds established in Proposition~\ref{prop:upper-bound} and
    \ref{prop:upper-bound-ot}.
    The example was chosen such that all bounds are different, and we find
    $\pi^*c < \pi^\mathrm{(a)}c < \pi^\mathrm{(c)}c < \pi^\mathrm{(b)}c$. The
    difference between $\pi^\mathrm{(b)}$ and $\pi^\mathrm{(c)}$ is that the
    vertical movement of mass restricted to the fiber $x = 4$ is optimal in case
    of $\pi^\mathrm{(c)}$, while $\pi^\mathrm{(b)}$ transports along the product
    coupling on this fiber.
    The plan $\pi^\mathrm{(a)}$ corresponds to the best assignment possible if
    transport is restricted to vertical and horizontal movements only.
  }
  \label{fig:different-bounds}
\end{figure}

\paragraph{Upper bounds.}
Our next goal is to derive upper bounds for the transport dependency. To this
end, we bound the infimum in the optimal transport problem
\eqref{eq:ot-definition} by explicitly constructing feasible transport plans
$\pi \in \CC(\mu, \nu)$ that may or may not be optimal.
The idea behind these plans is to restrict the transport in $X\times Y$ to
the fibers $X\times \{y\}$ and $\{x\}\times Y$ for $x\in X$ and $y\in Y$. 
Consequently, we often deal with costs of the form $c(x, y_1, x, y_2)$
or $c(x_1, y, x_2, y)$, and will thus assume that $c$ is
controlled by suitable marginal costs $c_X$ and $c_Y$ on $X$ and $Y$ via
\begin{equation}\label{eq:marginal-cost-sup}
  c_X(x_1, x_2)
  \ge
  \sup_{y\in Y} c(x_1, y, x_2, y)
  \quad\qquad\text{and}\qquad\quad
  c_Y(y_1, y_2)
  \ge
  \sup_{x\in X} c(x, y_1, x, y_2)
\end{equation}
for all $(x_1,y_1), (x_2,y_2) \in X\times Y$.
To gain intuition for the three upper bounds we introduce in the following (and
the transport plans $\pi^\mathrm{(a)}$ to $\pi^\mathrm{(b)}$ behind them),
Figure~\ref{fig:different-bounds} can be consulted. The first two of them,
corresponding to the plans $\pi^\mathrm{(a)}$ and $\pi^\mathrm{(b)}$, are
established in the following result.

\begin{subequations}%
\begin{proposition}{upper bounds, diameter}{upper-bound}
  Let $\gamma \in \CC(\mu, \nu)$ for $\mu\in\PP(X)$ and $\nu\in\PP(Y)$ in Polish
  spaces $X$ and $Y$. If the cost function $c$ on $X\times Y$ satisfies
  \eqref{eq:marginal-cost-sup}, then
  \label{eq:upper-bounds}
  \begin{align}
    \td(\gamma)
    &\le
    (\gamma\otimes\gamma)\,\cmin \label{eq:upper-bound}\\
    &\le
    \min\!\big( (\mu\otimes\mu)\,c_X, (\nu\otimes\nu)\,c_Y\big), \label{eq:upper-bound-margin}
  \end{align}
  where $\cmin\colon (X\times Y)^2 \to [0, \infty]$ is the cost
  function defined by
  \begin{equation}\label{eq:cmin}
    \cmin(x_1, y_1, x_2, y_2) = \min\!\big(c_X(x_1, x_2), c_Y(y_1, y_2)\big).
  \end{equation}
\end{proposition}

For convenience, we refer to integrals of the form $(\nu\otimes\nu)\,c_Y$ as the
($c_Y$)-\emph{diameter} of $\nu$. Hence, inequality
\eqref{eq:upper-bound-margin} shows that the transport dependency $\td(\gamma)$
can be bounded in terms of the diameters of the marginals $\mu$ and $\nu$.

The third upper bound, which is based on the improvement $\pi^\mathrm{(c)}$ of
the plan $\pi^\mathrm{(b)}$ in Figure~\ref{fig:different-bounds}, is
a worthwhile measure of association in its own right and receives its own
definition.

\begin{definition}{marginal transport dependency}{marginal-transport-dependency}
  Let $X$ and $Y$ be Polish spaces and $c_Y$ a continuous cost function on $Y$.
  The ($c_Y$)-\emph{marginal transport dependency} is defined via
  \begin{equation}
    \td^Y(\gamma)
    =
    \td^Y_{c_Y}(\gamma)
    =
    \int \ot_{c_Y}\big(\gamma(x, \cdot), \nu\big)\,\mu(\dif x).
  \end{equation}
\end{definition}

Again, we typically suppress the dependency on $c_Y$ in our notation if the
costs are apparent from the context.
Note also that we have to assume that the mapping $x\mapsto
\ot_{c_Y}\big(\gamma(x, \cdot), \nu\big)$ is measurable for this definition to
make sense. When $c_Y$ is continuous, measurability is guaranteed by
\textcite[Corollary~5.22]{villani2008}.
Since $\td^Y$ disregards the cost landscape on $X$, it is especially suited to
quantify dependency in asymmetric settings where points in $X$ cannot be
transported from one to another in a meaningful way, for instance if $X$ is
categorical.

\begin{proposition}{upper bounds, marginal transport}{upper-bound-ot}
  Let $\gamma \in \CC(\mu, \nu)$ for $\mu\in\PP(X)$ and $\nu\in\PP(Y)$ in Polish
  spaces $X$ and $Y$.
  If the cost function $c$ on $X\times Y$ satisfies the marginal bounds
  \eqref{eq:marginal-cost-sup} for a continuous $c_Y$, then
  \begin{equation}\label{eq:upper-bound-ot}
    \td(\gamma)
    \le
    \td^Y(\gamma)
    \le
    (\nu\otimes\nu)\,c_Y.
  \end{equation}
\end{proposition}
\end{subequations}%

It is worth pointing out that the marginal transport dependency is in fact
a special case of the (general) transport dependency under costs of the form
\begin{equation*}
  c_\infty(x_1, y_1, x_2, y_2)
  =
  \begin{cases}
    c_Y(y_1, y_2) &\text{if}~x_1 = x_2, \\
    \infty & \text{else}.
  \end{cases}
\end{equation*}
This indicates that $\td^Y$ arises as limit case if movements within the space
$X$ become prohibitively expensive, such that all transport eventually withdraws
to the fibers $\{x\}\times Y$ for $x\in X$. The following result confirms this
intuition to be accurate.
\begin{theorem}{marginal transport dependency as limit}{marginal-transport-dependence}
  Let $X$ and $Y$ be Polish spaces. For $\alpha > 0$, let $c_\alpha$ be a cost
  function on $X\times Y$ that satisfies
  \begin{equation}
    c_\alpha(x_1, y_1, x_2, y_2)
    \begin{cases}
      \,=\, c_Y(y_1, y_2)          &\text{if}~x_1 = x_2,\\
      \,\ge\, \max\big(\alpha\,c_X(x_1, x_2), c_Y(y_1, y_2)\big) &\text{else},
    \end{cases}
  \end{equation}
  where $c_Y$ is a continuous cost function on $Y$ and $c_X$ a positive cost
  function on $X$. Then, for any $\gamma\in\PP(X\times Y)$,
  \begin{equation}
    \lim_{\alpha\to\infty} \td_{c_\alpha}(\gamma)
    =
    \td_{c_\infty}(\gamma)
    =
    \td^Y_{c_Y}(\gamma).
  \end{equation}
\end{theorem}

\begin{example}{multivariate Gaussian, part 2}{gauss-marginal}
  We revisit the Gaussian setting of Example~\ref{ex:gauss}. Recall that $(\xi,
  \zeta)\sim\gamma = \Normal(\eta, \Sigma)$, where $\eta$ is a mean vector and
  $\Sigma$ a covariance matrix with blocks $\Sigma_{11}$, $\Sigma_{12}$,
  $\Sigma_{21}= \Sigma_{12}^T$ and $\Sigma_{22}$. This time, we work with costs
  of the form
  \begin{equation}\label{eq:gauss-marginal-costs}
    c_\alpha(x_1, y_1, x_2, y_2)
    =
    \alpha\,\|x_1 - x_2\|^2 + \|y_1 - y_2\|^2
  \end{equation}
  for $\alpha > 0$. Since these costs can be interpreted as scaling
  $\xi$ by the factor $\sqrt{\alpha}$ under the usual squared Euclidean distance
  $c_1$, we can employ the same arguments as in Example~\ref{ex:gauss} and find
  expressions for $\td_{c_\alpha}$ via replacing $\Sigma_{11}$ by
  $\alpha\,\Sigma_{11}$ and $\Sigma_{12}$ by $\sqrt{\alpha}\,\Sigma_{12}$.
  Adapting equation \eqref{eq:tdep-normal} in this way yields
  \begin{equation}
    \td_{c_\alpha}(\gamma)
    =
    2\alpha\,\mathrm{trace}(\Sigma_{11})+2\,\mathrm{trace}(\Sigma_{22})
    - 2\, \mathrm{trace}\left( 
    \begin{pmatrix}
      \alpha^2\Sigma_{11}^2 & \alpha^{3/2}\,\Sigma_{11} \Sigma_{12}  \\
      \alpha^{1/2}\,\Sigma_{22} \Sigma_{21} & \Sigma_{22}^2 
    \end{pmatrix}^{1/2}\right).
  \end{equation}
  According to Theorem~\ref{thm:marginal-transport-dependence}, the right-hand
  side converges to the marginal transport dependency $\td_{c_Y}^Y(\gamma)$ as
  $\alpha \to \infty$, where $c_Y$ is the squared Euclidean distance. In the
  bivariate case \eqref{eq:gauss2}, this limit can readily be calculated and
  reads
  \begin{align}
    \td_{c_Y}^Y(\gamma)
    &=
    \lim_{\alpha\to\infty} \td_{c_\alpha}(\gamma) \\
    &= 
    \lim_{\alpha\to\infty} 2\,\left(\alpha\,\sigma_{1}^2 + \sigma_{2}^2
     - \sqrt{\alpha^2\sigma_1^4+\sigma_2^4+2\alpha\sigma_1^2\sigma_2^2\sqrt{1 - \rho^2}}\right) \\
    &=
    2 \sigma_2^2 \left(1 - \sqrt{1 - \rho^2}\right). \label{eq:gauss2-marginal}
  \end{align}
  From this, it becomes apparent that $\td_{c_Y}^Y$ \enquote{loses} the metric
  information on the space $X$ carried via $\sigma_1$.
\end{example}
The three upper bounds \eqref{eq:upper-bound} to \eqref{eq:upper-bound-ot}
established above play a crucial role for the remainder of the manuscript. In
Section~\ref{sec:contractions}, we investigate under which conditions the
inequalities in
\eqref{eq:upper-bounds} are actually equalities. This uncovers an intimate
relation between maximal values of transport dependency on the one hand and
contracting couplings on the other. Afterwards, in
Section~\ref{sec:coefficients}, we apply this understanding to derive normalized
coefficients of dependency.

\section{Transport dependency: contractions and maximal values}
\label{sec:contractions}
Up to this point, we were concerned with general properties of the transport dependency
for generic cost functions. We now focus on a particular additive cost structure
that enables us to characterize under which conditions the upper bounds of the
previous section are attained. All proofs in this section are delegated to
Appendix~\ref{app:proofs}.

We equip the Polish space $X$ with a (general) cost function $k_X$ and the
Polish space $Y$ with a lower semi-continuous (pseudo-)metric $d_Y$, and we
consider costs of the form
\begin{subequations}\label{eq:additive-costs-and-margins}
\noeqref{eq:additive-costs}
\noeqref{eq:additive-costs-margins}
\begin{equation}\label{eq:additive-costs}
  c(x_1, y_1, x_2, y_2)
  =
  h\big(k_X(x_1, x_2) + d_Y(y_1, y_2)\big)
\end{equation}
for $(x_1, y_1), (x_2, y_2) \in X\times Y$, where $h\colon[0, \infty) \to
[0,\infty)$ is a strictly increasing function that satisfies $h(0) = 0$. We also
fix the marginal costs
\begin{equation}\label{eq:additive-costs-margins}
  c_X
  =
  h\circ k_X
  \qquad\text{and}\qquad
  c_Y
  =
  h\circ d_Y,
\end{equation}
\end{subequations}
which fulfill condition \eqref{eq:marginal-cost-sup}. Therefore, $c_X$ and $c_Y$
are suited for the upper bounds in Proposition~\ref{prop:upper-bound}.
To state our findings, we need to define couplings that describe contractions
between $X$ and $Y$. We say that $\gamma\in\PP(X\times Y)$ is \emph{contracting}
(on its support) if 
\begin{equation}\label{eq:def-contracting-general}
  d_Y(y_1, y_2) \le k_X(x_1, x_2)
  \qquad\text{for~all}\qquad(x_1,y_1), (x_2, y_2)\in\supp\,\gamma.
\end{equation}
We also need a slightly weaker condition and say that $\gamma$ is
\emph{almost surely contracting} if \eqref{eq:def-contracting-general} holds
$(\gamma\otimes\gamma)$-almost surely (but not necessarily on each single pair
of points of the support).

\begin{theorem}{contracting couplings}{contractions}
  Let $X$ and $Y$ be Polish spaces, $c$ be a cost function on $X\times Y$
  of the form~\eqref{eq:additive-costs-and-margins}, and $\gamma\in\CC(\cdot,
  \nu)\subset\PP(X\times Y)$ for $\nu\in\PP(Y)$. If $\gamma$ is contracting,
  \begin{equation}\label{eq:contraction-condition-general}
    \td(\gamma) = (\nu\otimes\nu)\, c_Y.
  \end{equation}
  Conversely, if $(\nu\otimes\nu)\,c_Y < \infty$ and
  \eqref{eq:contraction-condition-general} holds, then $\gamma$ is almost surely
  contracting.
\end{theorem}
This theorem reveals that vertical movements of mass, which the upper bound
$(\nu\otimes\nu)\,c_Y$ is based on, are optimal when mass is transported to
a contracting coupling $\gamma$ (recall Figure~\ref{fig:contractions} in this
context). Conversely, vertical movements are not optimal whenever $\gamma$ is
not (almost surely) contracting.
Note that the distinction between \enquote{contracting} and \enquote{almost
surely contracting} in Theorem~\ref{thm:contractions} is necessary, as
one can construct almost surely contracting couplings $\gamma$ such
that $\td(\gamma) = 0$ and $(\nu\otimes\nu)\,c_Y > 0$. If $d_Y$ and $k_X$ are
continuous,
the two concepts coincide (see
Lemma~\ref{lem:continuous-contraction} in Appendix~\ref{app:proofs}). In this
case, Theorem~\ref{thm:contractions} implies the equivalence
\begin{equation}
  \gamma\in\CC(\cdot, \nu)\text{ is contracting}
  \qquad\Longleftrightarrow\qquad
  \td(\gamma) = (\nu\otimes\nu)\,c_Y.
\end{equation}
Under mild assumptions, contracting couplings are always deterministic and
correspond to contracting functions.
We say that a measurable function $\varphi\colon X\to Y$ is \emph{$\mu$-almost
surely contracting} if there is a Borel set $A\subset X$ with $\mu(A) = 1$ and
\begin{equation}\label{eq:contraction-condition}
  d_Y\big(\varphi(x_1), \varphi(x_2)\big)
  \le
  k_X(x_1, x_2)
  \qquad
  \text{for all}
  \qquad
  x_1, x_2\in A.
\end{equation}
Combining the previous observations with Theorem~\ref{thm:contractions} yields
the following statement.

\begin{theorem}{contracting deterministic couplings}{contractions-deterministic}
  Let $X$ and Y be Polish spaces and $c$ be a cost function on $X\times Y$ of
  the form~\eqref{eq:additive-costs-and-margins}, where $h$ and $k_X$ are
  continuous and $d_Y$ is a continuous metric. Let $\gamma\in\CC(\mu, \nu)$ with
  $\mu\in\PP(X)$, $\nu\in\PP(Y)$, and $(\nu\otimes\nu)\,c_Y < \infty$. Then
  \begin{equation}
    \td(\gamma) = (\nu\otimes\nu)\,c_Y
  \end{equation}
  holds if and only if $\gamma = (\id,\varphi)_\#\mu$ for a $\mu$-almost surely
  contracting function $\varphi\colon X\to Y$.
\end{theorem}

A point deserving emphasis is that
Theorem~\ref{thm:contractions-deterministic} actually provides
a characterization of Lipschitz and Hölder functions, as well as of isometries.
If $k_X$ is set to $\alpha\cdot d_X$ for a metric $d_X$ on $X$, then
condition~\eqref{eq:contraction-condition} is equivalent to $\varphi$ being
$\alpha$-Lipschitz $\mu$-almost surely. Similarly, under the choice $k_X
= \alpha \cdot \smash{d_X^\beta}$ for $\beta\in(0, 1)$, condition
\eqref{eq:contraction-condition} coincides with the $\beta$-Hölder criterion for
$\varphi$ (with constant $\alpha$). Finally, the simultaneous equality
\begin{equation}
  \td(\gamma) = (\mu\otimes\mu)\, c_X = (\nu\otimes\nu)\,c_Y
\end{equation}
for $k_X = d_X$ holds if and only if $\varphi$ is $\mu$-almost surely an
isometry from $(X, d_X)$ to $(Y, d_Y)$.
In this context, we highlight that a $\mu$-almost sure contraction
$\varphi$ in the setting of Theorem~\ref{thm:contractions-deterministic} can
always be uniformly extended to the full support of $\mu$ if $(Y, d_Y)$ is
a Polish metric space (see Lemma~\ref{lem:uniform-extension} in
Appendix~\ref{app:proofs}).

It is interesting to contrast the condition $\tau(\gamma)
= (\nu\otimes\nu)\,c_Y$ at the heart of
Theorem~\ref{thm:contractions-deterministic} to the weaker condition
$\tau^Y(\gamma) = (\nu\otimes\nu)\,c_Y$, where $\tau^Y$ denotes the marginal
transport dependency established in
Definition~\ref{def:marginal-transport-dependency}. In the following statement,
$c_Y$ can be a general cost function and does not have to satisfy
\eqref{eq:additive-costs-margins}.

\begin{theorem}{measurable deterministic couplings}{marginal-transport-dependence-measurable}
  Let $X$ and $Y$ be Polish spaces and $c_Y$ a positive continuous cost function
  on $Y$. Let $\gamma\in\CC(\mu, \nu)$ with $\mu\in\PP(X)$, $\nu\in\PP(Y)$,
  and $(\nu\otimes\nu)\,c_Y < \infty$. Then
  \begin{equation}
    \td^Y(\gamma) = (\nu\otimes\nu)\,c_Y
  \end{equation}
  holds if and only if $\gamma = (\id,\varphi)_\#\mu$ for a measurable function
  $\varphi\colon X \to Y$.
\end{theorem}

In this context, we note that Conjecture~2 of \textcite{mori2020} states that
the equality $\td(\gamma) = \min\big(\td^Y(\gamma), \td^X(\gamma)\big)$ holds
under ``general conditions'' for costs $c = d_X + d_Y$ in Polish metric spaces.
If, however, $\gamma$ is concentrated on the graph of a bijective function
$\varphi\colon X\to Y$ such that $\varphi$ and $\varphi^{-1}$ are both
\emph{not} 1-Lipschitz, combining Theorem~\ref{thm:contractions-deterministic}
and Theorem~\ref{thm:marginal-transport-dependence} shows $\td(\gamma)
< \min\big(\td^Y(\gamma), \td^X(\gamma)\big)$. In fact, we believe (and
numerical tests suggest) that the claimed equality only holds in somewhat
special situations.

To conclude this section, we want to highlight that the preceding results equip
the transport dependency with a meaningful interpretation as quantifier of
dependence: the larger the transport dependency $\td(\xi, \zeta)$ between two
random variables $\xi$ and $\zeta$ is, the more they have to be associated in
a contracting manner. Intuitively, this means that the conditional law of $\zeta
\,|\, \xi = x$ must behave well as a function of $x\in X$, judged in terms of
$k_X$ and $d_Y$.
In fact, the highest possible degree of transport dependency for fixed marginals
is (under continuity assumptions) only assumed if $\xi$ and $\zeta$ are
deterministically related by $\zeta = \varphi(\xi)$ for a contraction $\varphi$.
Other deterministic relations between $\zeta$ and $\xi$, which exhibit rapid
changes that break condition \eqref{eq:contraction-condition}, are assigned
a lower degree of dependency. This way of dependency quantification can often be
desirable, especially in situations where quickly oscillating or chaotic
relations between $\xi$ and $\zeta$ practically cannot (or should not) be
distinguished from actual noise.

\section{Transport correlation}
\label{sec:coefficients}

\begin{figure}
  \centering
  {\footnotesize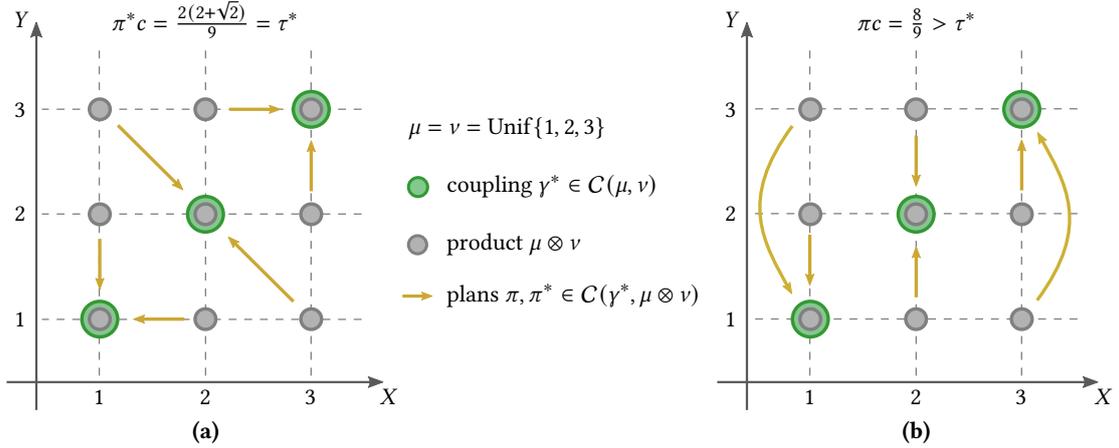}
  \caption{Example for which the upper bounds \eqref{eq:upper-bounds}
    are not sharp. The marginal distributions are $\mu
    = \nu = \Unif\{1, 2, 3\}$.
    \textbf{(a)} One can show that $\gamma^* = \Unif\{(1,1), (2,2), (3,3)\}$
    maximizes $\td(\gamma)$ over $\gamma\in\CC(\mu, \nu)$, assuming the value
    $\tau^* = 2(2 + \sqrt{2})/9$ if $c$ is the Euclidean distance on $\RR^2$.
    The visualized plan $\pi^*$ is optimal.
    \textbf{(b)} All of the established upper bounds move mass along a plan
    $\pi$ that is restricted to vertical or horizontal transports. The total
    transportation cost is $8/9 > \tau^*$. Note that the Euclidean distance on
    $X\times Y$ can not be expressed in the additive form
    \eqref{eq:additive-costs}.
  }
  \label{fig:sharpness-counter-example}
\end{figure}

In this section, we introduce several coefficients of association that are based
on the transport dependency. The central ingredient is upper
bound~\eqref{eq:upper-bound-margin} in Proposition~\ref{prop:upper-bound}, which
can be used to scale $\td$ to the interval $[0, 1]$ in a way that only depends
on the marginal distributions.
Without further assumptions, however, bound \eqref{eq:upper-bound-margin} is not
necessarily sharp, and values close to $1$ may be impossible (see
Figure~\ref{fig:sharpness-counter-example}).
We therefore focus on costs of the form \eqref{eq:additive-costs}, for which the
sharpness of the upper bounds is well understood
(Theorem~\ref{thm:contractions-deterministic} in
Section~\ref{sec:contractions}). Detailed proofs of the statements in this
section can be found in Appendix~\ref{app:proofs}.

For the sake of clarity, we work in a slightly less general setting than in
Section~\ref{sec:contractions} and right away assume $(X, d_X)$ and $(Y, d_Y)$
to be Polish metric spaces, restricting to costs
\begin{equation}\label{eq:additive-costs-coefficients}
  c(x_1, y_1, x_2, y_2) = \big(\alpha\cdot d_X(x_1, x_2) + d_Y(y_1, y_2)\big)^p
\end{equation}
for $x_1, x_2\in X$, $y_1, y_2\in Y$, and $\alpha, p > 0$. In what follows,
$d_X$, $d_Y$, and $p$ are usually considered to be fixed, and we mainly explore
the influence of $\alpha$.
We call a map $\varphi\colon X \to Y$
a \emph{dilatation} if there exists some $\beta > 0$ such that
\begin{equation}
  d_Y\big(\varphi(x_1), \varphi(x_2)\big) = \beta\,d_X(x_1, x_2)
\end{equation}
for all $x_1, x_2\in X$. A dilatation can be thought of as an isometry from $(X,
\beta d_X)$ to $(Y, d_Y)$, i.e., an isometry up to the correct scaling. We also
recall the notion of $p$-weak convergence on $\PPP{p}(X\times Y)$, which was
introduced and discussed in the context of continuity
(Section~\ref{sec:general-properties}).

\begin{definition}{$\alpha$-transport correlation}{transport-correlation}
  Let $(X, d_X)$ and $(Y, d_Y)$ be Polish metric spaces and consider costs $c$
  of the form \eqref{eq:additive-costs-coefficients} for $\alpha, p > 0$. For
  $\gamma\in\CC(\mu, \nu)$ with $\mu\in\PP(X)$ and $\nu\in\PP(Y)$ such that $0
  < (\nu\otimes\nu)\,d_Y^p < \infty$, the \emph{$\alpha$-transport correlation}
  is defined via
  \begin{equation}\label{eq:transport-correlation}
    \tc_\alpha(\gamma)
    =
    \left(\frac{\td(\gamma)}{(\nu\otimes\nu)\,d_Y^p}\right)^{1/p}.
  \end{equation}
\end{definition}
\begin{proposition}{}{transport-correlation}
  The $\alpha$-transport correlation $\tc_\alpha$ in
  Definition~\ref{def:transport-correlation} satisfies
  \begin{enumerate}[itemsep=.5ex,topsep=1.5ex]
    \item $\tc_\alpha(\gamma) = 0$ iff $\gamma = \mu\otimes\nu$,
    \item $\tc_\alpha(\gamma) = 1$ iff $\gamma = (\id,\varphi)_\#\mu$ for an
      $\alpha$-Lipschitz function $\varphi: X \to Y$,
    \item $\tc_\alpha(\gamma) = \tc_\alpha(\gamma')$ if $\gamma' = (f_X,
      f_Y)_\#\gamma$ for isometries $f_X\colon X\to X$ and $f_Y\colon Y \to Y$,
    \item $\gamma \mapsto \tc_\alpha(\gamma)^p$ is convex when restricted to
      $\CC(\cdot, \nu) \subset \PP(X\times Y)$ for fixed $\nu$,
    \item $\tc_{\alpha_n}(\gamma_n) \to\tc_\alpha(\gamma)$ as $n\to\infty$ if
      $(\gamma_n)_{n\in\NN}$ converges $p$-weakly to $\gamma$ and
      $\alpha_n \to \alpha$,
    \item $\alpha \mapsto \tc_\alpha(\gamma)^p$ is monotonically increasing for
      all $p > 0$ and concave if $p \le 1$,
  \end{enumerate}
  where the functions $\varphi$, $f_X$, and $f_Y$ only have to be defined $\mu$-
  or $\nu$-almost surely.
\end{proposition}
Properties 1 and 2 lend the transport correlation a distinctive interpretation
as dependency coefficient that identifies Lipschitz relations
between random variables $\xi$ and $\zeta$. 
Based on our earlier findings (Theorems~\ref{thm:marginal-transport-dependence}
and \ref{thm:marginal-transport-dependence-measurable}) on the marginal
transport dependency $\tau^Y$ (see Definition~\ref{def:marginal-transport-dependency}),
it is possible to extend $\tc_\alpha(\gamma)$ to $\alpha = \infty$.

\begin{definition}{marginal transport correlation}{marginal-transport-correlation}
  Let $X$ and $(Y, d_Y)$ be Polish (metric) spaces and $c_Y = d_Y^p$ for $p
  > 0$. For $\gamma\in\CC(\mu, \nu)$ with $\mu\in\PP(X)$ and $\nu\in\PP(Y)$
  such that $0 < (\nu\otimes\nu)\,d_Y^p < \infty$, the \emph{marginal transport
  correlation} is defined via
  \begin{equation}\label{eq:marginal-transport-correlation}
    \tc_\infty(\gamma)
    =
    \left(\frac{\td^Y(\gamma)}{(\nu\otimes\nu)\, d_Y^p}\right)^{1/p}.
  \end{equation}
\end{definition}

\begin{proposition}{}{marginal-transport-correlation}
  The marginal transport correlation $\tc_\infty$ in
  Definition~\ref{def:marginal-transport-correlation} satisfies
  \begin{enumerate}[itemsep=.5ex, topsep=1.5ex]
    \item $\tc_\infty(\gamma) = 0$ iff $\gamma = \mu\otimes\nu$,
    \item $\tc_\infty(\gamma) = 1$ iff $\gamma = (\id,\varphi)_\#\mu$ for a
      measurable function $\varphi: X \to Y$,
    \item $\tc_\infty(\gamma) = \tc_\infty(\gamma')$ if $\gamma' = (f_X,
      f_Y)_\#\gamma$ for a measurable injection $f_X\colon X\to X$ and
      a dilatation $f_Y\colon Y \to Y$,
    \item $\gamma \mapsto \tc_\infty(\gamma)^p$ is convex when restricted to
      $\CC(\cdot, \nu)\subset\PP(X\times Y)$ for fixed $\nu$,
  \end{enumerate}
  where the functions $\varphi$, $f_X$, and $f_Y$ only have to be defined $\mu$-
  or $\nu$-almost surely.
\end{proposition}
The differences in Proposition~\ref{prop:transport-correlation} and
\ref{prop:marginal-transport-correlation} reflect that $\tc_\alpha$
progressively loses the sense for metrical structure in $X$ when $\alpha$ is
increased. In the limit $\alpha = \infty$, the Lipschitz restrictions appearing
in Proposition~\ref{prop:transport-correlation} are dissolved and only
conditions of measurability remain.

Since we normalize with the diameter of $\nu$ instead of the one of $\mu$ in
Definition~\ref{def:transport-correlation}, the transport correlation
$\tc_\alpha$ singles out functional relations in the direction $X \to Y$. One
possibility to compensate for this asymmetry is to adapt the value of $\alpha$
to the diameters.

\begin{definition}{isometric transport correlation}{isometric-transport-correlation}
  Let $(X, d_X)$ and $(Y, d_Y)$ be Polish metric spaces.
  For $\gamma\in\CC(\mu, \nu)$ with $\mu\in\PP(X)$ and $\nu\in\PP(Y)$, consider
  costs of the form \eqref{eq:additive-costs-coefficients} with $p > 0$ and
  \begin{equation}\label{eq:isometric-alpha}
    \alpha
    =
    \alpha_*
    =
    \left(\frac{(\nu\otimes\nu)\,d_Y^p}{(\mu\otimes\mu)\,d_X^p}\right)^{1/p}.
  \end{equation}
  For $0 < \alpha_* < \infty$, the \emph{isometric transport correlation} is
  defined via $\tc_*(\gamma) = \tc_{\alpha_*}(\gamma)$.
\end{definition}

\begin{proposition}{}{isometric-transport-correlation}
  The isometric transport correlation $\tc_*$ in
  Definition~\ref{def:isometric-transport-correlation} satisfies
  \begin{enumerate}[itemsep=.5ex, topsep=1.5ex]
    \item $\tc_*(\gamma) = 0$ iff $\gamma = \mu\otimes\nu$,
    \item $\tc_*(\gamma) = 1$ iff $\gamma = (\id,\varphi)_\#\mu$ or $\gamma
      = (\psi, \id)_\#\nu$ for dilatations $\varphi: X \to Y$ or $\psi: Y\to X$,
    \item $\tc_*(\gamma) = \tc_*(\gamma')$ if $\gamma' = (f_X,
      f_Y)_\#\gamma$ for dilatations $f_X\colon X\to X$ and $f_Y\colon Y \to Y$,
    \item $\gamma \mapsto \tc_*(\gamma)^p$ is convex when restricted to $\CC(\mu,
      \nu)\subset\PP(X\times Y)$ for fixed $\mu$ and $\nu$,
    \item $\tc_*(\gamma_n) \to\tc_*(\gamma)$ as $n\to\infty$ if
      $(\gamma_n)_{n\in\NN}$ converges $p$-weakly to $\gamma$,
    \item $\tc_*(\gamma) = \tc_*(\gamma')$ if $\gamma' = f_\#\gamma$ for the
      symmetry map $f(x, y) = (y, x)$,
  \end{enumerate}
  where the functions $\varphi$, $\psi$, $f_X$, and $f_Y$ only have to be
  defined $\mu$- or $\nu$-almost surely.
\end{proposition}

We can interpret this choice of $\alpha_*$ as first normalizing
the metric measure spaces $(X, d_X, \mu)$ and $(Y, d_Y, \nu)$ by their
$p$-diameters before calculating the transport dependency.

To conclude, we mention another symmetric coefficient of
association which can be derived from the transport dependency.
Instead of dividing by the
diameter of $\nu$ in \eqref{eq:transport-correlation}, one can divide by the
minimum of the diameters of $\nu$ and $\mu$. Setting $\alpha = 1$ in
\eqref{eq:additive-costs-coefficients}, this results in the coefficient
\begin{equation}\label{eq:contracting-transport-correlation}
  \left(
    \frac{\td_c(\gamma)}{\min\big((\mu\otimes\mu)\,d_X^p, (\nu\otimes\nu)\,d_Y^p\big)}
  \right)^{1/p},
\end{equation}
which has, for $p = 1$, been introduced as \emph{Earth mover's correlation} by 
\textcite{mori2020}. It enjoys similar properties to the isometric transport
correlation. In Conjecture~3 of \textcite{mori2020} it is speculated that the
expression in \eqref{eq:contracting-transport-correlation} is $1$ if and only if
$\gamma = (\id,\varphi)_\#\mu$ for a dilatation $\varphi\colon X \rightarrow Y$.
However, our results in Section~\ref{sec:contractions} clarify that this is the
case if and only if $\gamma$ is concentrated on the graph of a 1-Lipschitz
function from $X$ to $Y$ \emph{or} a 1-Lipschitz function from $Y$ to $X$.

\section{Estimation and computation}
\label{sec:estimation}

We now discuss various strategies to estimate the transport dependency
$\td(\gamma)$ in the presence of $n\in\NN$ empirical i.i.d.\ observations $(\xi_i,
\zeta_i)_{i=1}^n \sim \gamma^{\otimes n}$ for $\gamma\in\CC(\mu, \nu)$, where
$\mu\in\PP(X)$ and $\nu\in\PP(Y)$. The general approach will be to construct an
estimator $\gamma_n$ of $\gamma$ as well as an estimator $(\mu\otimes\nu)_n$ of
$\mu\otimes\nu$ and plug them into the optimal transport functional, yielding
estimates of the general form
\begin{equation}\label{eq:plugin-estimator}
  \hat\td_n
  =
  \ot_c\big(\gamma_n, (\mu\otimes\nu)_n\big).
\end{equation}
A striking observation, which we substantiate below, is that such estimators
naturally seem to exhibit \emph{lower complexity adaptation} (LCA,
\cite{hundrieser2022empirical}), meaning that the intrinsic dimension of
$\gamma$ determines the statistical convergence rate to the true value -- and
not the dimension of the product $\mu\otimes\nu$ that is potentially higher.

\paragraph{Product estimator.}
The most immediate estimates of the probability measures $\gamma$ and
$\mu$, and $\nu$ are given by their empirical counterparts $\hat\gamma_n
= \frac{1}{n} \sum_{i=1}^n \delta_{(\xi_i, \zeta_i)}$, $\hat\mu_n = \frac{1}{n}
\sum_{i=1}^n \delta_{\xi_i}$, and $\hat\nu_n = \frac{1}{n} \sum_{i=1}^n
\delta_{\zeta_i}$. This results in the estimator
\begin{equation}\label{eq:product-estimator}
  \td(\hat\gamma_n)
  =
  \ot_c(\hat\gamma_n, \hat\mu_n\otimes\hat\nu_n).
\end{equation}
We call this the \emph{product estimator} as it relies on the product of the
empirical measures $\hat\mu_n$ and $\hat\nu_n$ to estimate $\mu\otimes\nu$.
Since this estimator is simply the transport dependency of $\hat\gamma_n$,
several properties follow immediately from established results. For example,
strong consistency, meaning that $\lim_{n\to\infty}\td(\hat\gamma_n)
= \td(\gamma)$ almost surely, follows from Proposition~\ref{prop:continuity}
under modest moment requirements.
Further continuity-related properties can be derived from
Theorem~\ref{thm:continuity-metric-power}, or, for additive costs,
equation~\eqref{eq:td-lipschitz}. If the latter can be applied, the product
estimator satisfies 
\begin{equation}\label{eq:lca-metric}
  |\td(\hat\gamma_n)^{1/p} - \td(\gamma)^{1/p}|
  \lesssim
  \ot_c(\gamma, \hat\gamma_n)^{1/p},
\end{equation}
which implies a convergence rate that only depends on $\gamma$ (and not on
$\mu\otimes\nu$). This demonstrates lower complexity adaptation. Indeed, as we
explore in Appendix~\ref{app:lca}, LCA of $\td(\hat\gamma_n)$ holds for
non-metric and non-additive costs as well. For instance, the following
statement is a corollary of Theorem~\ref{thm:lcaproduct} together with metric
entropy bounds derived in \textcite{hundrieser2022empirical}.

\begin{corollary}{lower complexity adaptation}{lca-manifold}
  Let $X$ and $Y$ be smooth manifolds and $\gamma\in\PP(X\times Y)$ such that
  $\supp\,\gamma$ is contained in a compact smooth manifold of dimension
  $s\in\NN$. If the cost function $c$ on $X\times Y$ is twice continuously
  differentiable, then
  \begin{equation}
    \Exp\, |\td(\hat\gamma_n) - \td(\gamma)| \lesssim n^{-2/s}.
  \end{equation}
\end{corollary}

The major practical drawback of the product estimator consists of the associated
computational burden. Its calculation relies on solving an optimal assignment
problem between $n^2$ and $n$ points. Even for fast algorithms like the network
simplex or cost scaling, this implies a (worst-case) runtime of $O(n^5)$ (up to
logarithmic factors, see \cite{peyre2019}).
We next discuss alternative estimators that use only $O(n)$ instead of $n^2$
points to estimate $\mu\otimes\nu$, reducing the worst-case time complexity to
$O(n^3)$ (again, up to logarithmic factors).

\paragraph{Splitting estimators.}
Using a suitable sample splitting procedure, i.i.d.\ observations can be
constructed to estimate $\gamma$ and $\mu\otimes\nu$.
In a setting with $2n$ observations $(\xi_i, \zeta_i)_{i=1}^{2n}
\sim\gamma^{\otimes 2n}$, for example, this can be realized by letting
$\gamma_n$ be the empirical measure of the first $n$ data points and
$(\mu\otimes\nu)_n$ be the empirical measure of $(\xi_i, \zeta_{n+i})_{i=1}^n
\sim(\mu\otimes\nu)^{\otimes n}$ in the estimator \eqref{eq:plugin-estimator}.
The resulting estimator is of the form of an empirical optimal transport cost,
so the LCA framework developed in \textcite{hundrieser2022empirical} is 
applicable.
However, splitting estimators are inefficient from a practical point of view and
primarily serve as a proof-of-concept that the LCA-rates of the product
estimator can also be expected when $\mu\otimes\nu$ is estimated by $O(n)$
points only.

\paragraph{Sampling estimators.}
A more efficient approach to estimate $\mu\otimes\nu$ is to randomly
sample points from the empirical product set $A = (\xi_i,
\zeta_j)_{i,j=1}^n\subset X\times Y$.
This can be done in various ways, for example by drawing $N = O(n)$ times from
$A$ either with or without replacement. Alternatively,
samples of the form $\big(\xi_i, \zeta_{\sigma(i)}\big)_{i=1}^n$, where $\sigma$
is uniformly distributed over the set of permutations of $\{1, \ldots, n\}$, can
be employed.
To improve the estimate, we can also repeat this scheme for $k$
different random permutations $\sigma = (\sigma_1, \ldots, \sigma_k)$, ending up
with a total of $N = kn$ sampled points. This procedure yields estimators of the
form
\begin{equation}\label{eq:permutation-estimator}
  \hat\td^\sigma_n
  =
  \ot_c\big(\hat\gamma_n, (\mu\otimes\nu)_{n}^\sigma\big)
  \qquad\text{where}\qquad
  (\mu\otimes\nu)_{n}^{\sigma}
  =
  \frac{1}{kn} \sum_{r=1}^k \sum_{i=1}^n \delta_{(\xi_i, \zeta_{\sigma_r(i)})},
\end{equation}
which we call \emph{permutation estimators}. In Appendix~\ref{app:lca}, we prove
that these estimators satisfy the LCA property as well
(Theorem~\ref{thm:lcasampling} and Corollary~\ref{cor:lcarepeated}) and thus
exhibit no worse convergence rates than the product estimator. In particular,
Corollary~\ref{cor:lca-manifold} still holds when the estimate
$\td(\hat\gamma_n)$ is replaced by $\hat\td^\sigma_n$.

\paragraph{Other estimators.}
All of the previous proposals for estimates of $\td(\gamma)$ are based on
(randomly) picking points in the set $A = (\xi_i, \zeta_j)_{i,j=1}^n$ in order
to approximate $\mu\otimes\nu$. We focus on these estimators in our simulation
study (Section~\ref{sec:applications}) and data application
(Section~\ref{sec:application-to-gene-expression-data}). However, they are not
the only reasonable options and various alternatives could be explored. For
example, instead of randomly sampling from $A$, the estimator
$(\mu\otimes\nu)_n$ may be based on a systematic approximation of $A$ with
$O(n)$ weighted support points, e.g., realized via clustering. Based on our
previous observations, it is to be expected that estimates of this form should
display LCA properties as well.

Still, the comparably large bias of optimal transport costs in high dimensions
will cause the expected statistical convergence rate to be slow in complex
settings. As a remedy, further steps to exploit the potential smoothness of the
measure $\gamma$ could be pursued, like kernel density or wavelet estimators for
$\gamma$ and $\mu\otimes\nu$
\parencite{weed2019estimation,deb2021rates,manole2021sharp,divol2022measure}.
Presently, these types of estimators are mostly confined to the realm of theory,
however, as their computation poses several difficulties.

\paragraph{Estimating the marginal transport dependency.}
We finally want to discuss estimation strategies for the marginal transport
dependency $\td^Y(\gamma)$, see Definition~\ref{def:marginal-transport-dependency},
for which the plug-in approach is often not feasible: in settings where
$\mu$ is diffuse, the empirical data $\hat\gamma_n$ likely obeys functional
relations $\hat\gamma_n = (\id, \hat\varphi_n)_\#\hat\mu_n$ for
$\hat\varphi_n\colon X \to Y$. Thus, $\td^Y(\hat\gamma_n)$ always equals
$(\hat\nu_n\otimes\hat\nu_n)\,c_Y$ and is not informative (recall
Theorem~\ref{thm:marginal-transport-dependence-measurable}).

Similar problems also affect other quantifiers of (unstructured) dependency,
like the mutual information, and the usual remedy is to preprocess the data. In
Euclidean settings, for example, one can first estimate a smooth kernel-density
or to bin the data with a kernel size or bin width $h_n$ that decreases to $0$
as $n$ grows in order to obtain a consistent estimator (see
\cite{tsybakov2008}). In a certain sense, the limiting procedure of $\alpha
\to\infty$ in Theorem~\ref{thm:marginal-transport-dependence} has a similar
effect: a finite $\alpha < \infty$ allows for some leeway along the space $X$
when matching $\hat\gamma_n$ with $\hat\mu_n\otimes\hat\nu_n$, and this leeway
becomes progressively smaller as $\alpha$ increases. Indeed, it is possible to
find suitable sequences $\alpha_n\to\infty$ so that
$\td_{c_{\alpha_n}}(\hat\gamma_n)$ is a consistent estimator of $\td^Y(\gamma)$
under mild assumptions.

\begin{proposition}{consistent estimation of $\td^Y$}{marginal-consistency}
  Let $(X, d_X)$ and $(Y, d_Y)$ be Polish metric spaces, $c_Y=d_Y^p$, and
  $c_\alpha = (\alpha d_X+d_Y)^p$ for $p \ge 1$ and $\alpha > 0$. Moreover,
  assume $\gamma\in \PP(X\times Y)$ to have a finite $p$-moment with respect to
  $d = d_X+d_Y$. If $(\alpha_n)_{n\in\NN}$ is a positive diverging sequence that
  satisfies
  $\alpha_n \cdot \EE\,\ot_{c_1}(\hat\gamma_n, \gamma)^{1/p} \to 0$ as $n\to\infty$,
  then
  \begin{equation}
    \Exp\,\big|\td_{c_{\alpha_n}}(\hat\gamma_n)^{1/p}
    - \td^Y_{c_Y}(\gamma)^{1/p}\big| \to 0.
  \end{equation}
\end{proposition}

Another route to consistently estimate the marginal transport dependency from
data is described in \textcite{wiesel2021}, who derives convergence rates if
$\td^Y(\gamma)$ is estimated by $\td^Y(\tilde\gamma_n)$ for a so-called adapted
empirical measure $\tilde\gamma_n$. An example for such an adaption is the
projection of the individual observations $\xi_i$ and $\zeta_i$ to a grid in the
unit cube (see \cite{backhoff2020}).

\section{Simulations}
\label{sec:applications}

We now investigate the performance of the transport dependency on simulated
data. After briefly outlining our settings and numerical methodologies, we
proceed to compare the product and permutation estimators introduced in the
previous section. In particular, we confirm the LCA property for both of them
and show that their rates of convergence are identical.
Afterwards, we conduct a series of benchmarks that illuminate properties of the
coefficients $\tc_\alpha$ and $\tc_*$ in Euclidean spaces. We also investigate
commonalities as well as differences to other coefficients of association, for
example by comparing their performance in permutation tests for independence.
Additional simulations that cover alternative parameter choices and adopt
different dependency models can be found in Appendix~\ref{app:simulations}.

\paragraph{Setting.}
In the following, we restrict to Euclidean spaces $X = \RR^r$ and $Y = \RR^q$
for $r,q\in\NN$, and equip them with their respective Euclidean metrics $d_X$
and $d_Y$. We focus on joint distributions $\gamma\in\PP(\RR^r \times \RR^q)$
that are either given deterministically via $\gamma = (\id, \varphi)_\#\mu$ for
$\varphi\colon[0, 1]^r \to [0, 1]^q$ and $\mu = \Unif[0, 1]^r$ (i.e.,
concentrated on the graph of the function $\varphi$), or by placing (uniform)
mass on more general shapes in $[0, 1]^{r + q}$.
The number of samples independently drawn from $\gamma$ for the purpose of
empirical estimation is denoted by $n \in \NN$, and we write $\hat\gamma_n$ to
refer to the corresponding empirical measure, see Section~\ref{sec:estimation}.
To study the influence of statistical noise, we later in the section also
consider convex contamination models of the form
\begin{equation}\label{eq:convex-contamination}
  \gamma^\epsilon = (1-\epsilon)\,\gamma + \epsilon\,(\mu\otimes\nu),
\end{equation}
where $\epsilon \in[0, 1]$ denotes the noise level distorting $\gamma\in\CC(\mu,
\nu)$. Observations $(x,y)\in\RR^r\times\RR^q$ sampled from $\gamma^\epsilon$
are with probability $(1-\epsilon)$ randomly drawn from $\gamma$ and with
probability $\epsilon$ randomly drawn from $\mu\otimes\nu$.
Further simulations that use additive Gaussian noise instead are provided in
Appendix~\ref{app:simulations}.

\paragraph{Computation.}
Depending on the size of the optimal transport problem, we alternate between two
different computational approaches. In sufficiently small settings, like $n \le
100$ for the product estimator \eqref{eq:product-estimator} and all simulated
$n$ with the permutation estimator \eqref{eq:permutation-estimator}, we employ
a network simplex based solver provided by the python optimal transport (POT)
package \parencite{flamary2017}.
To still be able to capture the behavior of the product estimator for larger
values of $n$ (up to $n = 1000$), we additionally employ our own
implementation\footnote{%
  The code is available in the julia programming language under
  \href{https://gitlab.gwdg.de/thomas.staudt/otter.jl}{gitlab.gwdg.de/thomas.staudt/otter.jl}.
}
of an approximation scheme proposed by \textcite{schmitzer2019}. It is based on
the Sinkhorn algorithm for entropically regularized optimal transport
\parencite{cuturi2013} and operates by successively decreasing the
regularization constant $\eta > 0$ until a suitable approximation of the
non-regularized problem is obtained. While $\eta$ is scaled down, in our case
from $\eta = 10^{-1}$ to $\eta = 10^{-3}$, the increasing sparsity of the
transport plan (after negligible entries are removed) is exploited by using
sparsity-optimized data structures. On systems equipped with a decent GPU, the
product estimator can thereby be calculated accurately for $n = 1000$ (as is
done in Figure~\ref{fig:bias} below) within a couple of seconds to minutes for
generic costs.

\paragraph{Lower complexity adaptation.}

\begin{figure}
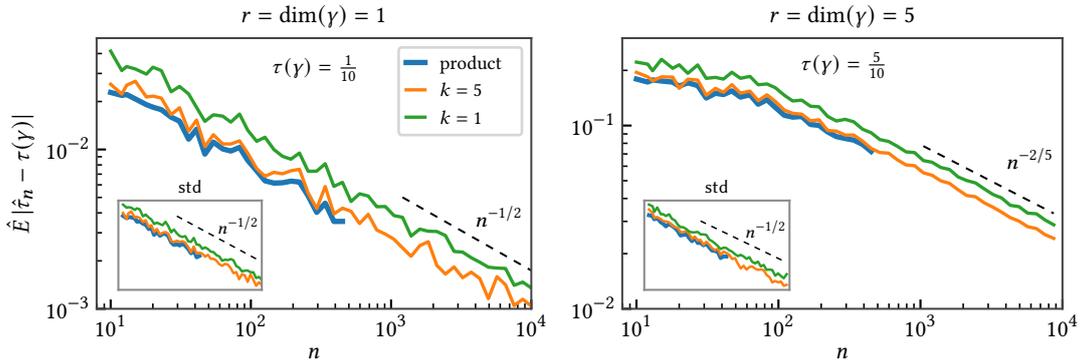

  \centering
  \input{figures/lca-permutation-1.pgf}\hspace{-2.0em}
  \input{figures/lca-permutation-5.pgf}

  \vspace{-1.5em}

  \caption{Lower complexity adaptation of the product estimator $\hat{\td}_n
    = \td(\hat\gamma_n)$ as well as the permutation estimator $\hat\td_n
    = \hat\td_n^\sigma$, defined in \eqref{eq:permutation-estimator}, based on $k \in\{1,
    5\}$ permutations.
    Depicted is a Monte-Carlo estimate (50 independent repetitions) of the mean
    absolute error $\EE|\hat\td_n - \td(\gamma)|$ as a function of the sample
    size $n$. The coupling $\gamma = (\id, \id)_\#\,\Unif[0, 1]^r$ corresponds
    to the choice $\xi = \zeta \sim \Unif[0, 1]^r$ and has intrinsic dimension
    $r$, while the product of its marginals $\mu\otimes\nu$ has dimension $2r$.
    The embedded plots show that the empirical standard deviation of $\hat\td_n$
    decays with the rate $n^{-1/2}$ independent of the dimension.
  }
  \label{fig:lca-permutation}
\end{figure}

In order to investigate the convergence behavior of the estimators proposed in
Section~\ref{sec:estimation}, we have to operate in a setting where
$\td(\gamma)$ can be calculated explicitly (or at least be approximated
very well). This is trivially the case for independent couplings $\gamma
= \mu\otimes\nu$ with $\td(\gamma) = 0$, but since the dimensionality of
$\gamma$ and $\mu\otimes\nu$ coincides for such $\gamma$, we would not be able
to discern lower complexity adaptation of the estimators. Instead, we consider
the case $\gamma = (\id, \id)_\# \mu$ and $\mu\sim\Unif[0,1]^r$, where analytical
solutions are feasible for squared Euclidean costs (see
Appendix~\ref{app:explicit}). The intrinsic dimensionality of $\gamma$ is equal
to $r$. Therefore, according to Corollary~\ref{cor:lca-manifold} and the
associated results in Appendix~\ref{app:lca}, we expect that
\begin{equation}\label{eq:lca-rate-simulation}
  \Exp\,|\hat\td_n - \td(\gamma)| \lesssim n^{-2/r}
\end{equation}
if $r > 4$ for both the product estimator $\hat\td_n = \td(\hat\gamma_n)$ and
the permutation estimator $\hat\td_n = \hat\td_n^\sigma$.
Figure~\ref{fig:lca-permutation} depicts the results of
Monte-Carlo simulations of $\Exp\,|\hat\td_n - \td(\gamma)|$ for $r\in\{1, 5\}$
under different choices of $\hat\td_n$. As expected, the product estimator
performs best. However, it is only marginally better than the permutation
estimator with $k = 5$ random permutations.
The simulations confirm that the upper bound \eqref{eq:lca-rate-simulation}
correctly characterizes the decay of the mean absolute error, and show that the
LCA property of the estimators affects the finite sample regime considered in
Figure~\ref{fig:lca-permutation}.

\paragraph{Dependency coefficients.}

We continue our numerical study by comparing different dependency coefficients
that assume values in $[0, 1]$.
Besides the isometric and $\alpha$-Lipschitz transport correlations $\tc_*$ and
$\tc_\alpha$ for $\alpha > 0$ (see Section~\ref{sec:coefficients}), which we
estimate via the product estimator \eqref{eq:product-estimator} unless specified
otherwise, we consider the following commonly applied coefficients:
\begin{description}[labelindent=1em, labelwidth=2em, parsep=0.5ex, leftmargin=!]
  \item[\textbf{cor}] the Pearson correlation. It is only
    applicable if $r = q = 1$. Since it assumes values in $[-1, 1]$, we always
    report its absolute value.
  \item[\textbf{spe}] the Spearman rank correlation coefficient. It is only
    applicable if $r = q = 1$ and is comparable to Kendall's $\tau$,
    another popular rank based correlation coefficient.
  \item[\textbf{dcor}] the Euclidean distance correlation
    \parencite{szekely2007}, based on the distance covariance defined in
    equation~\eqref{eq:euclidean-dcov}. It is applicable for all $r, q\in\NN$.
    In its generalized form \eqref{eq:dcov}, it is also applicable in generic
    (separable) metric spaces of (strong) negative type. As discussed previously,
    some of its properties make it comparable to $\tc_*$.
    Note that
    we use the vanilla empirical distance correlation in our simulations.
    Related estimators, like an unbiased estimator of $\dcor^2$ proposed by
    \textcite{szekely2013distance}, showed a comparable performance when testing
    for independence and are not included in this study.
  \item[\textbf{mic}]  the maximal information coefficient
    \parencite{reshef2011}. It is only applicable if $r = q = 1$.
    For its estimation, we use the estimator $\mathrm{mic_e}$ \parencite{reshef2016measuring},
    which we compute via the tools provided by
    \textcite{albanese2012}. The two algorithmic parameters $c$ and $\alpha$ were set
    to 5 and 0.75, respectively (as recommended in
    \cite{albanese2018}).
\end{description}
For each of these coefficients, generically called $\rho$ for the moment, we are
interested in several features.
Apart from the actual value of $\rho(\gamma)$, which signifies the amount of
dependency attributed to $\gamma$, we look at the variance and bias of
$\rho(\hat\gamma_n)$ as an estimator of $\rho(\gamma)$ when data is limited.

To check how well the coefficients are able to distinguish structure from noise,
we also include the results of permutation tests for independence (see
\cite[Section~15.2]{lehmann2006}, or \cite{janssen2003} for background on
permutation based tests).
The $\rho$-based permutation test we employ works as follows:
for given data $z = (x_i, y_i)_{i=1}^n \in (X\times Y)^n$, assumed to be sampled
from $\gamma^{\otimes n}$, we write $\rho(z)$ for the empirical estimate
of $\rho(\gamma)$ based on $z$. Furthermore, we denote $z_\sigma = (x_i,
y_{\sigma(i)})_{i=1}^n$, where $\sigma$ is a permutation of $n$ elements. For
the test, we randomly select $m$ permutations $\sigma_1, \ldots, \sigma_m$ and
reject the null hypothesis that $z$ is sampled from an independent coupling
$\gamma = \mu\otimes\nu$ if
\begin{equation}\label{eq:permutation-test}
  \big|\big\{i\,\colon\rho(z_{\sigma_i}) > \rho(z)\big\}\big| \le k
\end{equation}
for some $k \in\{0, \ldots, m\}$. Since a permutation of the second components
does not affect the distribution of $\gamma^{\otimes n}$ if $\gamma = \mu\otimes\nu$ is
a product coupling, this leads to a level $(k+1) / (m+1)$ test.
In all of our applications, we choose $k$ and $m$ such that this level is $\le 0.1$.
Note that the power of the permutation test will usually increase if $m$ is
increased while the level is held constant.

\begin{figure}[t!]
  \centering
  \includegraphics[width=.9\textwidth]{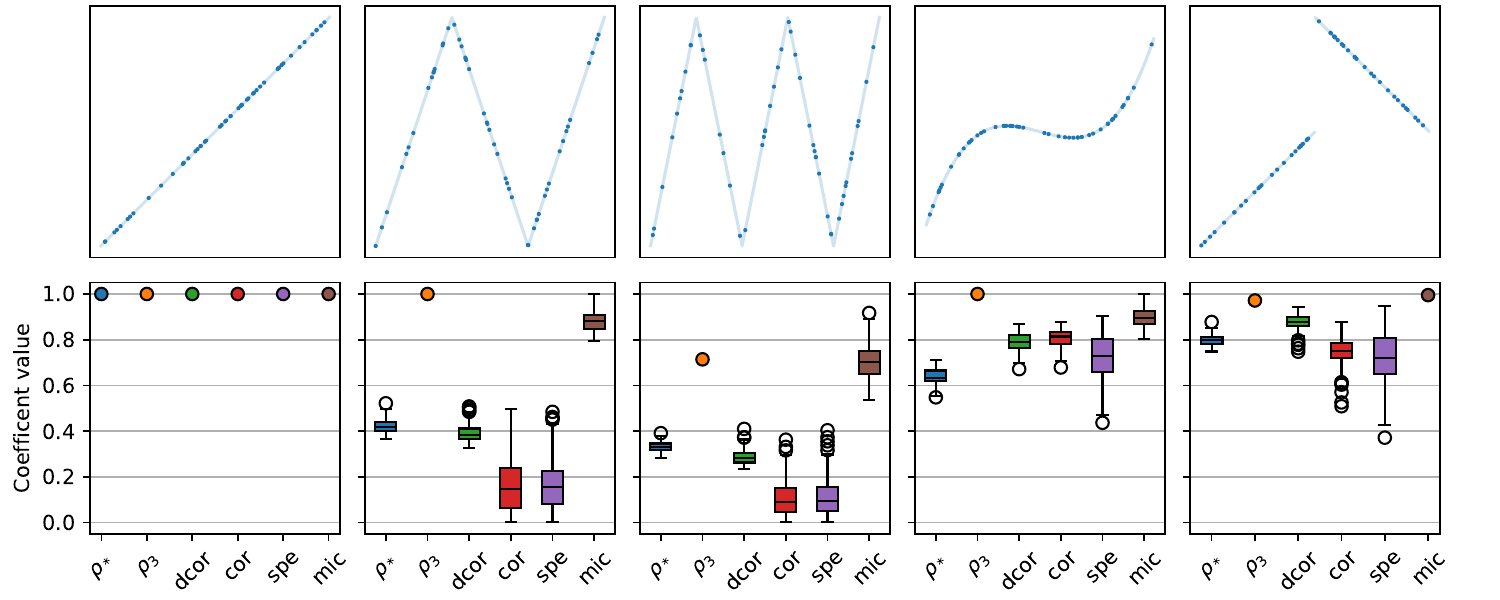}
  \includegraphics[width=.9\textwidth]{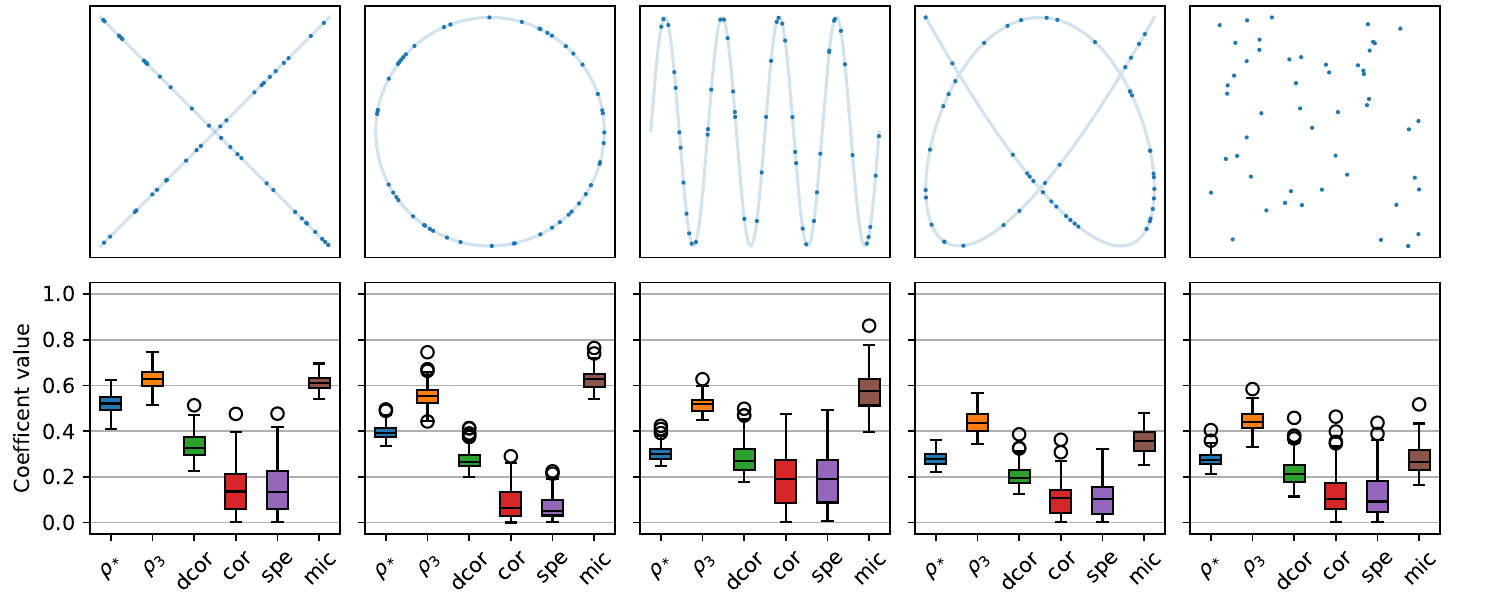}
  \caption{Comparison of several dependency coefficients on two dimensional
    geometries. As reference, independent $\Unif[0, 1]^2$ noise has also been
    simulated (bottom right). For each geometry, the scatter plot on top displays an exemplary
    sample of size $n = 50$. The box plots below summarize the values of the
    empirically estimated coefficients based on $100$ such samples. For
    the sake of visualization, boxes in the box plot are replaced by single dots if
    their extent would be smaller than 0.1.}
  \label{fig:geometric}
\end{figure}

\paragraph{Recognizing shapes.}
We begin with an investigation of the behavior of the aforementioned
coefficients in two dimensional settings, meaning $r = q = 1$.
Figure~\ref{fig:geometric} contains box plots depicting the coefficients'
performance on simple geometries, like lines or circles, for $n = 50$ samples.
As a point of reference, we also include uniform noise on $[0, 1]^2$.

To showcase the Lipschitz-selectivity of $\tc_\alpha$, we chose to include the coefficient
$\tc_3$ for $\alpha = 3$. Recall that this implies $\tc_3(\gamma) = 1$ whenever
$\gamma$ is concentrated on the graph of a function whose slope is at most 3.
Figure~\ref{fig:geometric} confirms this property of $\tc_3$, which is the
only coefficient to assign maximal dependency to the zigzag function with
slope 3 and the polynomial. On the zigzag function with 5 segments (and thus
slope 5), it has already decreased to about 0.7. In this context, the
$\mathrm{mic}$, which generally achieves high values on all geometries, performs
notably well. It is also the coefficient that most clearly distinguishes the
pretzel example from noise. Another noteworthy observation is that $\tc_*$ and
$\mathrm{dcor}$ indeed behave comparably, especially on functional relations.
On non-functional patterns, like the circle or the cross, $\tc_*$ assumes
somewhat higher values than $\mathrm{dcor}$. At the same time, the bias of
$\tc_*$ on independent noise (for which values of $0$ are expected in the limit
of large $n$) is slightly higher than the one of $\mathrm{dcor}$, albeit with
a smaller variance.
Finally, as to be expected, the Pearson and Spearman correlation coefficients do
a poor job at discerning non-monotonic structures. In case of the circle,
for example, the values of these coefficients are systematically lower
than when confronted with random noise.

\begin{figure}[t]
  \centering
  \includegraphics[width=.9\textwidth]{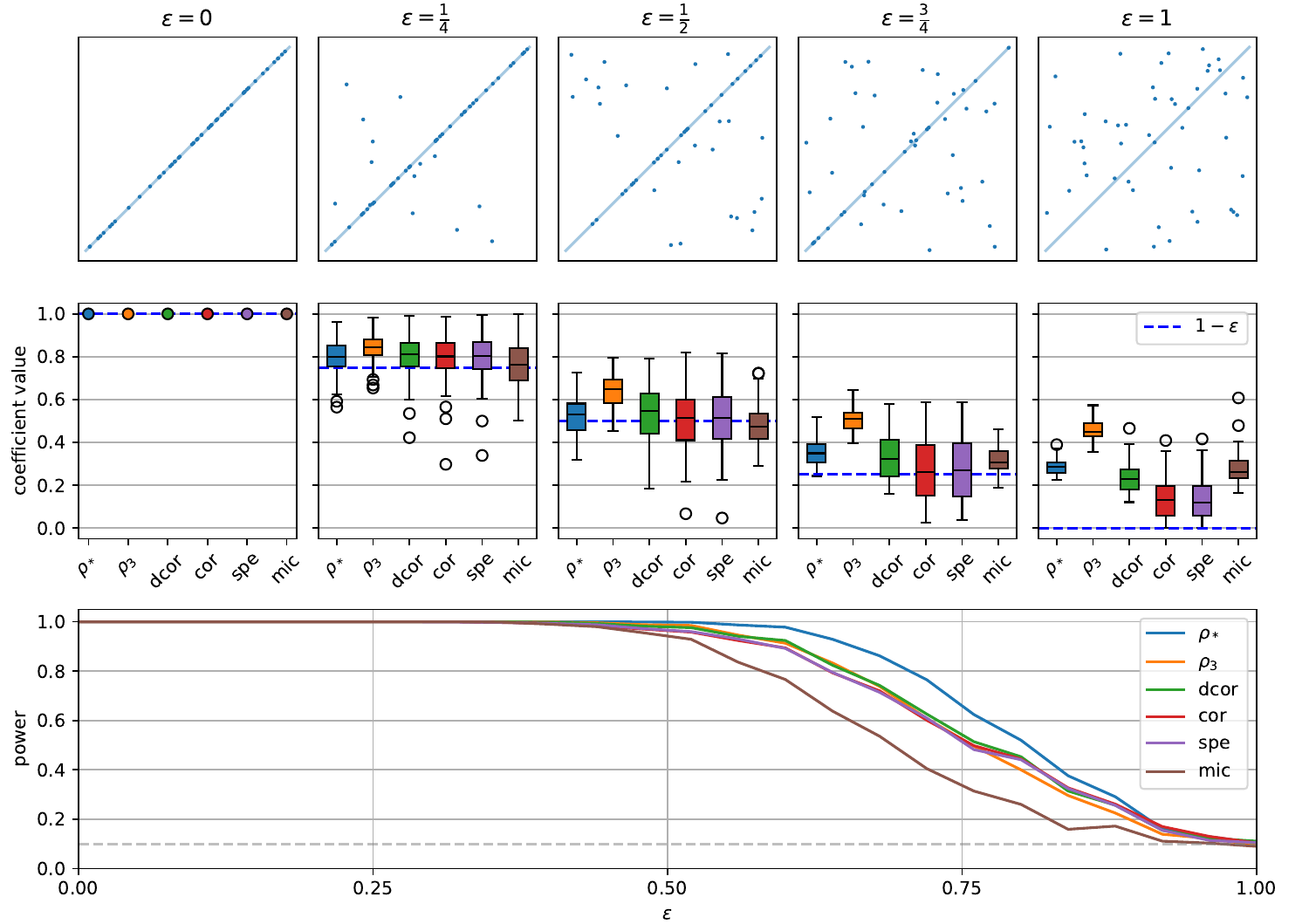}
  \vspace{-2ex}
  \caption{Dependency coefficients applied to increasingly noisy datasets
    according to the convex noise model in \eqref{eq:convex-contamination},
    where $\gamma$ is given in terms of the identity on $[0, 1]$.
    The scatter plots on the top show exemplary samples drawn from $\gamma^\epsilon$ of size $n = 50$.
    The box plots are based on 100 such samples. The power curves on the
    bottom display the results of the permutation tests described in
    \eqref{eq:permutation-test}. To estimate the power, 1000 tests were conducted
    per value of $\epsilon$ for each coefficient. The significance level (dashed line) of
    these tests is $10\%$.}
  \label{fig:noise1}
\end{figure}

\begin{figure}[t!]
  \centering
  \includegraphics[width=.90\textwidth]{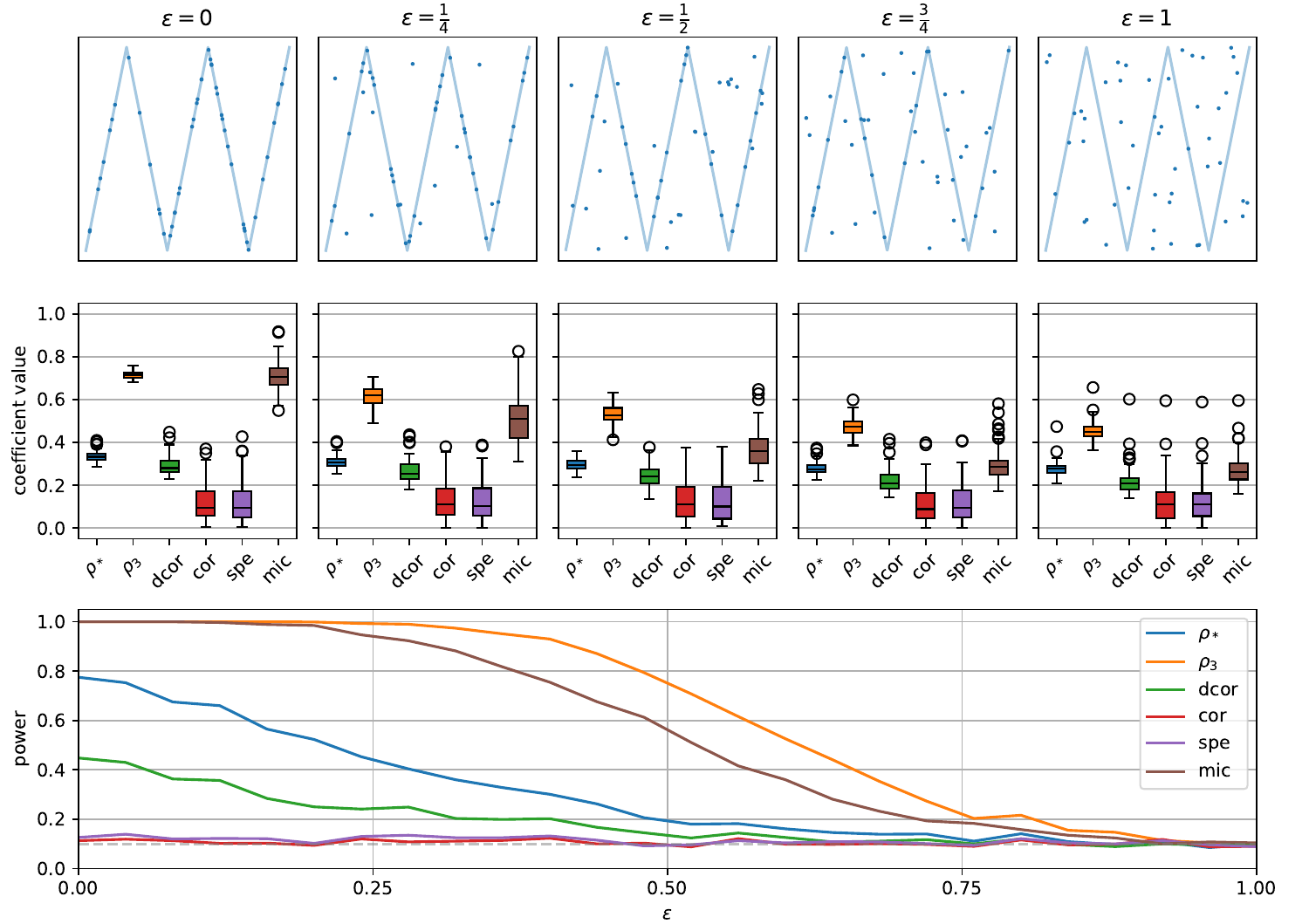}
  \vspace{-2ex}
  \caption{Dependency coefficients applied to increasingly noisy datasets
    according to the convex noise model in \eqref{eq:convex-contamination},
    where $\gamma$ is given in terms of a zigzag function with maximal slope $5$.
    See Figure \ref{fig:noise1} for a description of the individual graphs.}
  \label{fig:noise2}
\end{figure}

\paragraph{Behavior under noise.}
Our next simulations concern the performance under convex noise models
$\gamma^\epsilon$ as defined in equation~\eqref{eq:convex-contamination}.
Figure~\ref{fig:noise1} and \ref{fig:noise2} illustrate the coefficients'
empirical estimates and their power when used for independence testing in case
that $\gamma$ is deterministically given in terms of the identity
(Figure~\ref{fig:noise1}) or the zigzag function with maximal slope $5$
(Figure~\ref{fig:noise2}).
In Appendix~\ref{app:simulations}, this type of comparison can be found for all
other distributions considered in Figure~\ref{fig:geometric} as well, and we
also present results under additive Gaussian noise.

In case of the identity, all coefficients seem to behave roughly similar,
especially for small noise levels. Under pure noise ($\epsilon
= 1$), for which the coefficients should attain the value $0$ in the limit
of large $n$, the transport correlation exhibits a comparably large bias at
a relatively small variance. This trend of high biases becomes even more serious
in higher dimensions and is further investigated below (Figure \ref{fig:bias}).
Regarding the test performance, the power curves in Figure~\ref{fig:noise1}
reveal that all coefficients except $\tc_*$ and $\mathrm{mic}$ perform
comparably. The power of $\tc_*$ is consistently higher than its competitors',
while $\mathrm{mic}$ performs notably worse.

The picture changes substantially for the zigzag example in
Figure~\ref{fig:noise2}. Due to the absence of monotonicity, $\mathrm{cor}$ and
$\mathrm{spe}$ are not able to distinguish data points originated from the
zigzag function from the ones coming from the independence coupling of its
marginals.
The coefficient $\mathrm{dcor}$ and, to a lesser extent, $\mathrm{\tc_*}$ also
assume lower values and can only partially discern the dependency structure
under noise.
Meanwhile, the coefficients $\mathrm{mic}$ and specifically $\tc_3$ lie
systematically higher and are still able to recognize dependency under high
noise levels.

\begin{figure}
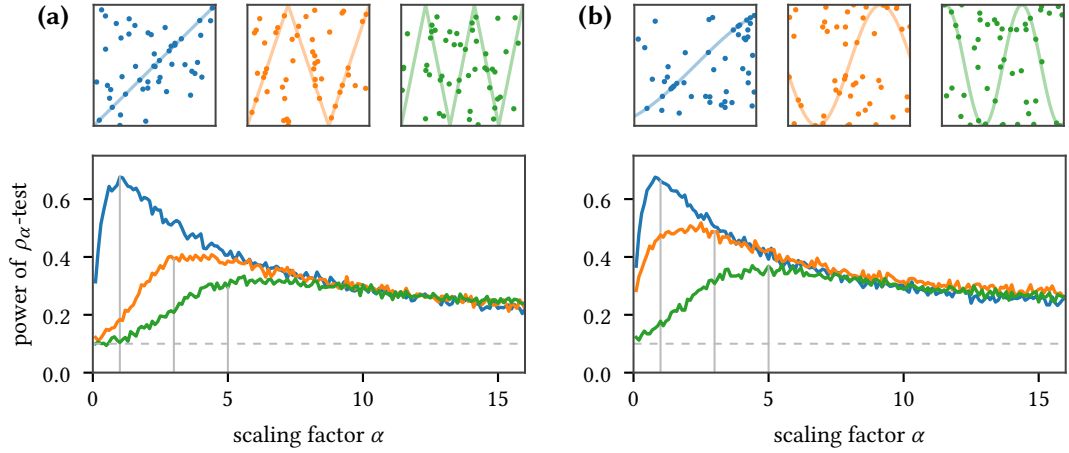

  \centering
  \input{figures/adapt-zigzag-50.pgf}\hspace{-0.75em}
  \input{figures/adapt-sinus-50.pgf}
  \caption{Influence of $\alpha$ on the test performance of
    $\tc_\alpha$ for linear \textbf{(a)} and sine \textbf{(b)} dependencies at
    a noise level of $\epsilon = 0.75$, see \eqref{eq:convex-contamination}. The
    graphs on the top display exemplary samples together with the actual functional
    relations encoded in $\gamma$ for different slopes. From left to right, the
    respective curves have maximal slopes of $1$, $3$, and $5$ in both (a) and
    (b). The respective values of $\alpha\in\{1, 3, 5\}$ are marked as vertical
    lines in the power plot below. For each value of $\alpha$, the power is
    estimated based on $500$ samples of size $n = 50$.}
  \label{fig:adapt}
\end{figure}

\paragraph{The role of $\boldsymbol{\alpha}$.}
In Figure~\ref{fig:power-intro} of the introduction, we already noted that
$\tc_3$ performs very well on a $3$-Lipschitz functional relation. Additionally,
Figure~\ref{fig:noise2} testifies that $\tc_3$ also outperforms all other
considered coefficients when applied to a $5$-Lipschitz
relation. It stands to reason, however, that $\tc_5$ would do even better than
$\tc_3$ in this setting. To further examine the claim that adapting $\alpha$ to
the slopes inherent to $\gamma$ improves the performance, we conducted a series
of
numerical experiments with functions of different maximal slope and varying
$\alpha$. The findings are displayed in Figure~\ref{fig:adapt}, where the power
of $\tc_\alpha$-based permutation tests is plotted as a function of $\alpha$ (for
fixed noise levels $\epsilon = 0.75$).  The results clearly suggest that
choosing a value of $\alpha$ close to the Lipschitz constant of the actual
functional relation in $\gamma$ can significantly improve the test performance.
If $\alpha$ is chosen too small or too large, the power systematically declines.

\paragraph{Bias and high dimensions.}
Another prevalent trend in the numerical results so far is the notable bias of
$\tc_*$ and $\tc_\alpha$ on independent noise (i.e., for $\epsilon = 1$), where
we expect a value of $0$ for $n\to\infty$. In fact, empirical
estimators of optimal transport distances are known to be susceptible to
a certain degree of bias, especially in high dimensions.  In
Figure~\ref{fig:bias}, we therefore compare the empirical estimation of
$\tc_*$ and $\mathrm{dcor}$ in settings of different
dimensions for sample sizes $n$ running from $10$ to $1000$.
We observe that the bias of $\tc_*$ and $\mathrm{dcor}$ seem comparable for $r
= q = 1$. If the dimensions are chosen higher, the bias increases much quicker
for $\tc_*$ than for $\mathrm{dcor}$. At the same time, the variance of the
$\tc_*$ estimates is (in part much) smaller for all choices of dimensions.
Indeed, in case of $r = q = 5$ and $n = 1000$, the estimated standard
deviation of $\tc_*$ is multiple times smaller than the one of
$\mathrm{dcor}$, even though the bias is substantial ($\tc_* \approx 0.72$,
while $\mathrm{dcor}\approx 0.13$).

\begin{figure}
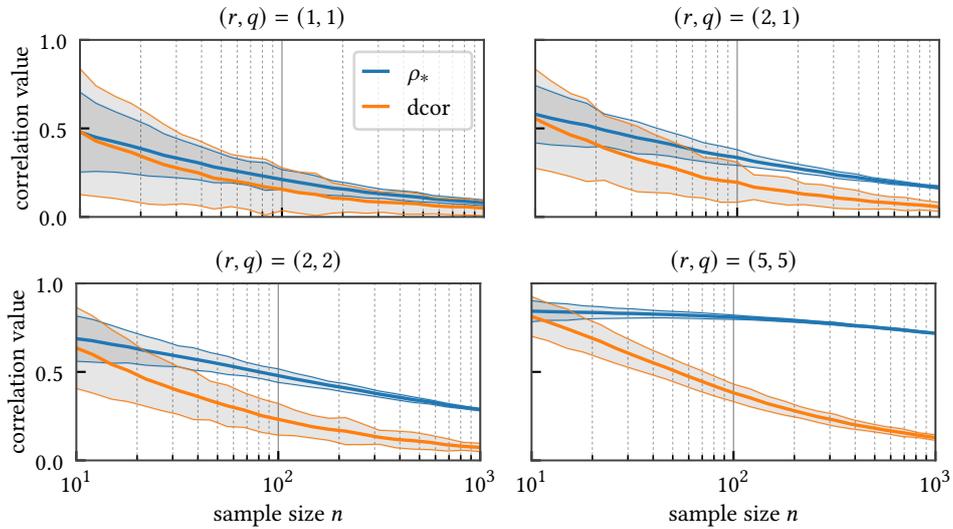

  \centering
  \input{figures/bias-iunif-1-1.pgf}\hspace{-3em}
  \input{figures/bias-iunif-2-1.pgf}\vspace{-6ex}
  \input{figures/bias-iunif-2-2.pgf}\hspace{-3em}
  \input{figures/bias-iunif-5-5.pgf}
  \vspace{-1.5ex}
  \caption{Empirical estimates of the isometric transport correlation and the
    distance correlation as a function of the sample size $n$ in different
    dimensions $r$ and $q$. The true distribution in the respective graphs is
    given by the product measure $\gamma = \Unif[0, 1]^r\otimes\Unif[0, 1]^q$.
    Therefore, $\tc_*(\gamma) = \dcor(\gamma) = 0$ is the (unbiased)
    value to be expected for $n\to\infty$. The shaded error regions correspond
    to $\pm$ three times the estimated standard deviation. The number of samples
    employed for the estimation of the mean and variance was adaptively decreased
    from $500$ for $n = 10$ to $20$ for $n = 1000$.}
  \label{fig:bias}
\end{figure}

\begin{figure}[t]
  \centering
  \input{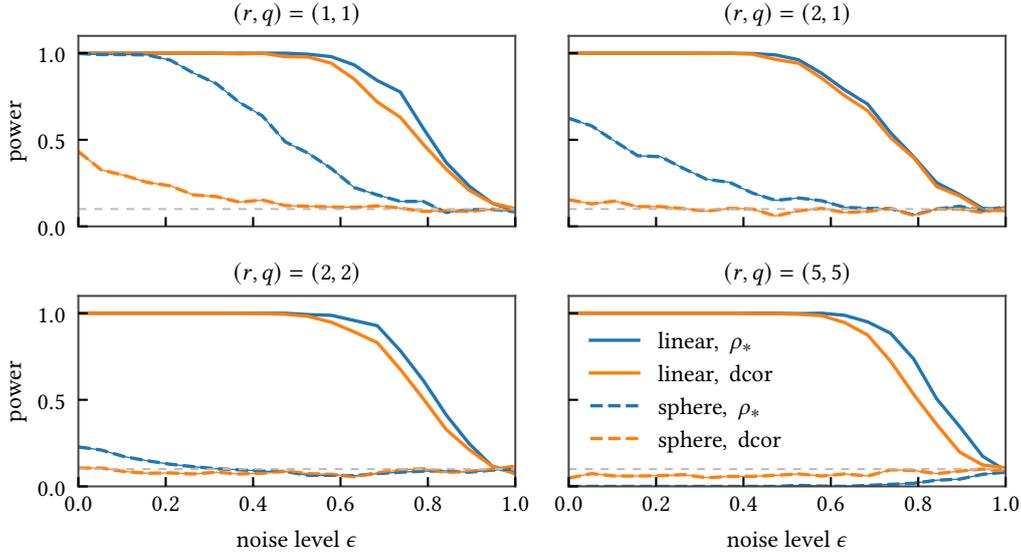}
  \vspace{-1.5ex}
  \caption{Test performance of the isometric transport correlation and the
    distance correlation in different dimensions $r$ and $q$. Shown are results
    for both linear and spherical relations under the noise
    model~\eqref{eq:convex-contamination}. In the former case, $\gamma$ is given
    by $(\id, \id)_\#\Unif[0, 1]^r$ if $r = q$ and by $(\id,
    \mathrm{pr}_1)_\#\Unif[0, 1]^r$ in case of $r = 2$ and $q = 1$, where
    $\mathrm{pr}_1$ denotes the projection on the first coordinate. In case of
    the sphere, $\gamma$ is uniformly distributed on the surface
    $\mathbb{S}^{r+q-1}\subset\RR^{r+q}$. For all power curves, $1000$ samples of
    size $n = 50$ were used.}
  \label{fig:higher-dims}
\end{figure}

A high bias in itself does not necessarily mean that
$\tc_*$ is blind to dependencies in high dimensions, however.
Figure~\ref{fig:higher-dims} reveals that simple linear structures with
noise of the form \eqref{eq:convex-contamination} are still recognized somewhat
better via $\tc_*$ than via $\dcor$ for $n = 50$, particularly in the setting
$r = q = 5$. To examine a more involved example, we also look at spherical
dependencies, where $\gamma$ is the uniform distribution on the sphere
$\mathbb{S}^{r+q-1}\subset\RR^{r}\times\RR^q$. In this case, the results are
more unintuitive and prompt several questions.
First of all, spherical dependencies are separated from noise (much) better by
the transport correlation than by $\dcor$ in settings with $r, q\in\{1, 2\}$.
At the same time, for $r = q = 5$, both the distance correlation and the
transport correlation consistently exhibit test powers that are smaller than
$0.1$, meaning that a random sample from the sphere regularly results in
\emph{lower} estimates of $\tc_*$ and $\dcor$ than drawing from the
corresponding marginals $\mu\otimes\nu$ would. This effect is particularly
severe in case of the transport dependency.
Even though this observation could in part be caused by the small sample size $n
= 50$ and the rather weak dependency in the spherical setting, it demonstrates
that detecting dependency in more complex situations remains an open issue that
merits further investigation. In fact, the data application in the next
section showcases that the transport dependency can yield meaningful results
even for $r=5000$ and $n < 100$.

\section{Application to gene expression data}
\label{sec:application-to-gene-expression-data}

In this section, we re-analyze a breast cancer gene expression study
\parencite{van2002gene} using the transport dependency, confirming the findings
by \textcite{behr2020}, who developed a specialized method for this data set.
Our results demonstrate that the transport dependency is able to sensibly detect
structural relations even when relying on non-metric similarity criteria in high
dimensions.

\begin{figure}[h]
  \centering
	\includegraphics[scale=0.7]{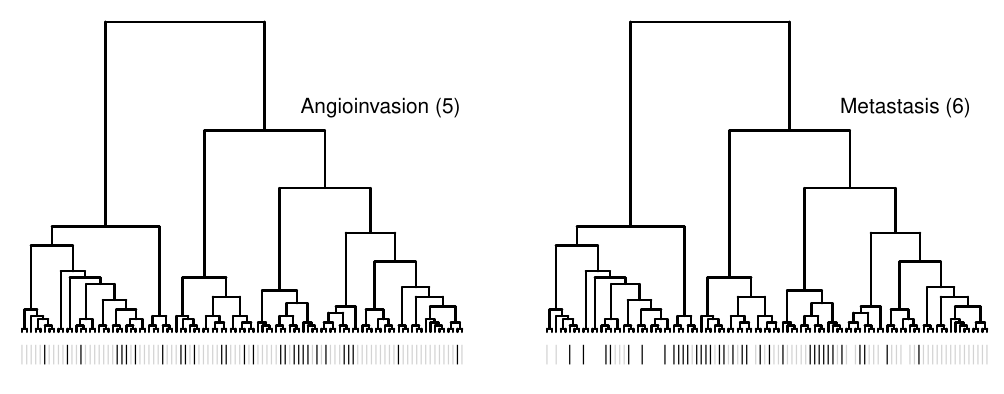}
  \caption{%
    Plot of the gene expression level derived hierarchical tree (with
    patients as leaves) and the corresponding binarized responses angioinvasion
    ($k=5$) and development of distant metastasis ($k=6$). Visually, the
    presence of dependency between the respective tree structure and the binary
    responses is not clearly evident, but the transport dependency based
    permutation test confidently rejects the null-hypothesis of independence for
    $k = 6$ while not rejecting it for $k = 5$, confirming the results obtained
    by \textcite{behr2020}.
  }
	\label{fig:tree}
\end{figure}

In the original study \parencite{van2002gene}, gene expression levels of
breast cancer samples were collected from 98 patients along with six clinical
responses: BRCA mutation ($k = 1$), estrogen receptor expression ($k = 2$),
histological grade ($k = 3$), lymphocytic infiltration ($k = 4$), angioinvasion
($k = 5$), and development of distant metastasis within 5 years ($k = 6$). 
We model the genetic data as i.i.d.\ random variables $\xi = (\xi_1, \dots,
\xi_{98})$, each taking values in the high dimensional space $X = \RR^{5000}$.
The components in this space correspond to the expression levels of $5000$
different genes. We furthermore let $\zeta^k = (\zeta_1^{k}, \dots,
\zeta_{98}^k)$ for $k\in\{1, \dots, 6\}$ denote the $k$-th responses
with values in $Y = \RR$, again assumed to be i.i.d.
To quantify the similarity between the gene expressions of different patients,
\textcite{van2002gene} employ a biologically motivated gene expression
correlation coefficient $\kappa\colon X^2\rightarrow [-1, 1]$ for a hierarchical
clustering scheme, see Figure~\ref{fig:tree}.
Based on the same approach, \textcite{behr2020} propose a highly specialized
test for independence between the resulting tree structure and the
responses $\zeta$.

The hierarchical clustering of the patient data provides benefits for
visualization and is a crucial preprocessing step for the techniques in
\textcite{behr2020} and other dependency measurement tools that require metric
properties, like the distance correlation. However, it introduces an additional
source of error, as the subsequent data analysis has to be performed
conditionally on this clustering. In contrast, clustering is not
necessary when working with the transport dependency: we can directly resort to
the gene expression correlation $\kappa$ by, e.g., considering the (non-metric)
cost function
\begin{equation}
  c_X(x_1, x_2)
  =
  1 - |\mathrm{\kappa}(x_1, x_2)|
\end{equation}
on $X$. On the space $Y$, we choose the absolute value $d_Y(y_1, y_2) = |y_1
- y_2|$ as cost. For each of the response variables, we then evaluate the
transport dependency $\td(\gamma_k)$, where
\begin{equation}
  \gamma_k = \frac{1}{n} \sum_{i=1}^n \delta_{(\xi_i, \zeta_i^{k})},
\end{equation}
under the standardized cost function
\begin{equation}
  c(x_1, y_1, x_2, y_2)
  =
  \frac{c_X(x_1, x_2)}{(\hat\mu_n\otimes\hat\mu_n)\, c_X}
  +
  \frac{d_Y(y_1, y_2)}{(\hat\nu_n\otimes\hat\nu_n)\, d_Y}
\end{equation}
on $X \times Y$, where $\hat\mu_n = \frac{1}{n}\sum_{i=1}^n \delta_{\xi_i}$
and $\hat\nu_n = \frac{1}{n}\sum_{i=1}^n \delta_{\zeta_i^k}$.
The results of a permutation test for independence with $m = 1000$ permutations
are recorded in Table~\ref{tab:gene-expression-cor}.

Looking at the values of $\td$, we can interpret that the estrogen
receptor expression ($k = 2$) is particularly strongly related to the gene
expressions.
Moreover, from the $p$-values, we conclude that we can decidedly reject the
null-hypothesis of independence between the gene expression and all the response
variables, with the exception of the angioinvasion response ($k = 5$), which has
also a low transport dependency value. This conclusion coincides with the
findings by \textcite{behr2020},
demonstrating that the transport dependency can be used to detect and quantify
dependencies even when the data is high dimensional, provided that meaningful
criteria of similarity are available on the marginal spaces.

\begin{table}
  \begin{center}
  \caption{}
  \label{tab:gene-expression-cor}
  \begin{tabular}{ccccccc}
    response $k$	& 1 & 2 & 3 & 4 & 5 & 6 \\
    \toprule
    $\td(\gamma_k) $& 0.275 & 0.630 & 0.412 & 0.314 & 0.284 & 0.378 \\
	  $p$-value & 0.001 & 0.001 & 0.001 &0.001 & 0.220& 0.006
  \end{tabular}
  \end{center}
\end{table}

\section*{Acknowledgements}
Thomas Giacomo Nies acknowledges support of DFG (Deutsche
Forsch\-ungsgemeinschaft, German Research Foundation) CRC 1456.
Thomas Staudt acknowledges support of DFG RTG 2088.
Axel Munk is supported by, and Thomas Giacomo Nies and Thomas Staudt were in
part funded by, the DFG under Germany's Excellence Strategy EXC 2067/1-390729940.
We also want to thank Facundo Memoli and Yoav Zemel for insightful discussions
in the initial phase of this project and Mathias Drton for helpful comments.

\printbibliography
\clearpage

\begin{appendix}

\section{Omitted statements and proofs}
\label{app:proofs}

This appendix provides detailed proofs for all statements that are not proven in
the main text. Note that we rely on disintegration arguments to construct or
destruct transport plans on several occasions. This is often not the only viable
method, however, and several of our results could alternatively be tackled by
other techniques, like exploiting the cyclical monotonicity of optimal transport
plans (see \cite{villani2008}).

\paragraph{Convexity and invariance.}

In the following, we provide proofs for Proposition~\ref{prop:convexity} and
Proposition~\ref{prop:contamination} (related to the convexity of the transport
dependency). We furthermore establish a fundamental invariance property of
optimal transport (Lemma~\ref{lem:ot-invariance}), which then powers the proof
of Proposition~\ref{prop:invariance}.

\begin{proof}[Proof of Proposition~\ref{prop:convexity}]
  Let $\gamma_0, \gamma_1\in\CC(\mu, \cdot)$ with second marginals $\nu_i \in \PP(Y)$ for
  $i \in \{0, 1\}$, and let $\pi^*_i$ be an optimal transport plan between
  $\gamma_i$ and $\mu\otimes\nu_i$ with respect to $c$. Then, defining $\gamma_t
  = (1-t) \gamma_0 + t \gamma_1$ for $t\in[0, 1]$ (and similarly $\nu_t$),
  we find $\gamma_t\in\CC(\mu, \nu_t)$ and $\pi_t = (1-t)\pi^*_0
  + t \pi^*_1 \in\CC(\gamma_t, \mu\otimes\nu_t)$. Hence,
  \begin{equation}
    \td(\gamma_t)
    =
    \ot_c(\gamma_t, \mu\otimes\nu_t)
    \le
    \pi_t c
    =
    (1-t)\,\td(\gamma_0) + t\,\td(\gamma_1),
  \end{equation}
  which establishes convexity on $\CC(\mu, \cdot)$. The result on $\CC(\cdot,
  \nu)$ follows analogously.
\end{proof}

\begin{proof}[Proof of Proposition~\ref{prop:contamination}]
  Since $\td(\mu\otimes\nu) = 0$, the first part of the statement follows from
  Proposition~\ref{prop:convexity}. For equality in \eqref{eq:contamination} if
  $c$ is a metric, we consider the dual formulation
  \eqref{eq:ot-definition-dual-metric} and deduce
  \begin{equation}
    \td(\gamma_t)
    =
    \sup \big(\gamma_t f - (\mu\otimes\nu) f\big)
    =
    (1-t)\,\sup \big(\gamma f - (\mu\otimes\nu) f\big)
    =
    (1-t)\,\td(\gamma),
  \end{equation}
  where both suprema are taken over $f\in\Lip_1(X\times Y)$ with respect to $c$.
\end{proof}

\begin{lemma}{}{ot-invariance}
  Let $X$ and $Y$ be Polish spaces and let $c$ be a cost function on $Y$.
  Consider a measurable function $f\colon X\to Y$ and let $c_f \colon X^2 \to
  [0, \infty]$ be defined by $c_f(x, x') = c\big(f(x), f(x')\big)$. Then for any
  $\mu,\nu\in\PP(X)$,
  \begin{equation}
    \ot_c(f_\#\mu, f_\#\nu) = \ot_{c_f}(\mu, \nu).
  \end{equation}
\end{lemma}

\begin{proof}[Proof of Lemma~\ref{lem:ot-invariance}]
  By a change of variables, we can
  immediately establish $\pi c_f = (f\otimes f)_\#\pi c$ for any $\pi\in\CC(\mu,
  \nu)$, where $(f\otimes f)_\#\pi\in\CC(f_\#\mu, f_\#\nu)$. This shows
  \begin{equation}
    \ot_{c_f}(\mu, \nu)
    =
    \inf_{\pi\in\CC(\mu, \nu)} \pi c_f
    =
    \inf_{\pi\in\CC(\mu, \nu)} (f\otimes f)_\#\pi c
    \ge
    \inf_{\tilde{\pi}\in\CC(f_\#\mu, f_\#\nu)} \tilde\pi c
    =
    \ot_{c}(f_\#\mu, f_\#\nu).
  \end{equation}
  To prove the reverse inequality, we fix some $\tilde\pi\in\CC(f_\#\mu,
  f_\#\nu)$ and explicitly construct a measure $\pi\in\CC(\mu, \nu)$ that
  satisfies $\pi c_f = \tilde\pi c$.

  We construct $\pi$ in two steps. First, we define an intermediate
  measure $\pi'$ by the relation $\pi'(\dif x_1, \dif y_2)
  = \tilde{\pi}\big(f(x_1), \dif y_2\big)\,\mu(\dif x_1)$. In the second step,
  we set $\pi(\dif x_1, \dif x_2) = \pi'\big(\dif x_1, f(x_2)\big)\,\nu(\dif
  x_2)$. It is straightforward to check $\pi'\in\CC(\mu, f_\#\nu)$ and
  $\pi\in\CC(\mu, \nu)$ by applying substitution and utilizing the properties of
  conditioning. In a similar vein, consecutive steps of substitution also show
  \begin{align}
    \pi c_f
    &= \int c_f(x_1, x_2)\,\pi(\dif x_1, \dif x_2) \\
    &= \int c\big(f(x_1), f(x_2)\big)\,\pi'\big(\dif x_1, f(x_2)\big)\,\nu(\dif x_2) \\
    &= \int c\big(f(x_1), y_2\big)\,\pi'(\dif x_1, y_2)\,(f_\#\nu)(\dif y_2) \\
    &= \int c\big(f(x_1), y_2\big)\,\pi'(\dif x_1, \dif y_2) \\
    &= \int c\big(f(x_1), y_2\big)\,\tilde{\pi}\big(f(x_1), \dif y_2\big)\,\mu(\dif x_1) \\
    &= \int c(y_1, y_2)\,\tilde{\pi}(y_1, \dif y_2)\,(f_\#\mu)(\dif y_1)
    = \tilde\pi c.
  \end{align}
  This establishes $\ot_{c_f}(\mu, \nu) \le \ot_c(f_\#\mu, f_\#\nu)$ and thus the
  equality of the two transport costs. Note that $\ot_c(f_\#\mu, f_\#\nu)$ and
  $\ot_{c_f}(\mu, \nu)$ do not have to be finite for this result to hold, as all
  integrands are non-negative.
\end{proof}

\begin{proof}[Proof of Proposition~\ref{prop:invariance}]
  For $\gamma\in\CC(\mu, \nu)$ with $\mu\in\PP(X)$ and $\nu\in\PP(Y)$, we
  observe that $c_f = c$ as well as $f_\#(\mu\otimes\nu) = f_{X\#}\mu \otimes
  f_{Y\#}\nu$ and $f_\#\gamma\in\CC(f_{X\#}\mu, f_{Y\#}\nu)$. The conclusion
  then follows from Lemma~\ref{lem:ot-invariance}.
\end{proof}

\paragraph{Continuity.}

This section covers the proofs for Proposition~\ref{prop:lower-semi-continuity},
Proposition~\ref{prop:continuity}, and Theorem~\ref{thm:continuity-metric-power}
(continuity of the transport dependency). We furthermore provide two auxiliary
results (Proposition~\ref{prop:varying-costs} and Lemma~\ref{lem:uniform-costs})
that control optimal transport costs under uniform changes of the base costs.
These are needed in Section~\ref{sec:coefficients}.

\begin{proof}[Proof of Proposition~\ref{prop:lower-semi-continuity}]
  Let $(\gamma_n)_{n\in\NN} \subset\PP(X\times Y)$ be a sequence of probability
  measures that weakly converges to $\gamma\in\CC(\mu, \nu)$ for $\mu\in\PP(X)$
  and $\nu\in\PP(Y)$. By the continuous mapping theorem, the respective
  marginals $(\mu_n)_{n\in\NN}$ and $(\nu_n)_{n\in\NN}$ weakly converge to $\mu$
  and $\nu$. According to \textcite[Theorem 2.8]{billingsley2013}, this implies
  weak convergence of the product measures $\mu_n\otimes\nu_n$ to $\mu\otimes\nu$.
  The result of the proposition now follows along the lines of the proof of the
  lower semi-continuity of the optimal transport cost
  (see for example \cite[Proposition 7.4]{santambrogio2015}, which requires only
  slight adaptations to make it work for sequences in both arguments of
  $\ot_c$).
\end{proof}

\begin{proof}[Proof of Proposition~\ref{prop:continuity}]
In order to prove this claim, we first establish an auxiliary result.
\begin{lemma}{}{dominated-uniform-integrability}
  Let $X$ be a Polish space and $f,g\colon X \to [0, \infty)$ continuous
  functions with $g \le f$. For a sequence $(\mu_n)_{n\in\NN}\subset\PP(X)$
  converging weakly to $\mu\in\PP(X)$, it holds that
  \begin{equation}
    \lim_{n\to\infty} \mu_n f = \mu f < \infty
    \qquad\text{implies}\qquad
    \lim_{n\to\infty} \mu_n g = \mu g < \infty.
  \end{equation}
\end{lemma}
\begin{proof}
  Assume $\mu_n f \to \mu f < \infty$ as $n\to\infty$. According to
  \textcite[Lemma 5.1.7]{ambrosio2008}, this implies that $f$ is uniformly
  integrable with respect to $(\mu_n)_{n\in\NN}$. Fix some $k > 0$. Since $g(x)
  > k$ implies $f(x) > k$ for all $x\in X$, we find
  $\mu_n\big(\ind_{g>k}\,g\big) \le \mu_n\big(\ind_{f>k}\,f\big)$
  for all $n$ and conclude that $g$ is also uniformly integrable with respect to
  $(\mu_n)_{n\in\NN}$.
  This in turn asserts that $\mu_n g \to \mu g$ as $n\to\infty$.
\end{proof}

We now return to the proof of Proposition~\ref{prop:continuity}. We know
that $\liminf_{n\to\infty} \td(\gamma_n) \ge \td(\gamma)$ due to the 
semi-continuity of $\td$ (Proposition~\ref{prop:lower-semi-continuity}).  Thus,
it is enough to show $\limsup_{n\to\infty} \td(\gamma_n) =: \bar\tau \le
\td(\gamma)$. By passing to a (not explicitly named) subsequence, we can assume
that $\td(\gamma_n)$ converges to $\bar\tau$.

Let $\gamma_n \in \CC(\mu_n, \nu_n)$ and $\gamma \in\CC(\mu, \nu)$ for
suitable marginal distributions $\mu, \mu_n\in\PP(X)$ and $\nu,
\nu_n\in\PP(Y)$. Since $\gamma_n\weak\gamma$ by assumption, we note that
$\mu_n\weak \mu$ and $\nu_n\weak\nu$ as $n\to\infty$ by the continuous mapping
theorem.
Like in the proof of Proposition~\ref{prop:lower-semi-continuity}, we conclude
that $\mu_n \otimes \nu_n$ converges weakly to $\mu \otimes \nu$.
Next, we fix a point $(x_0, y_0)\in X\times Y$. According to Lemma~\ref{lem:dominated-uniform-integrability},
$\mu_n\,d_X(\cdot, x_0)^p = \gamma_n\,d_X(\cdot, x_0)^p \to \gamma\,d_X(\cdot,
x_0)^p = \mu\,d_X(\cdot, x_0)^p$ as $n\to\infty$, where we set $g(x,y)
= d_X(x, x_0)^p$, $f(x,y) = d\big((x,y), (x_0, y_0)\big)^p$, and
used that $\gamma_n f \to \gamma f < \infty$ (by assumption).
Therefore, $\mu_n$ converges $p$-weakly to $\mu$, and the same holds for $\nu_n$
by a similar argument.

In the next step, we bound the metric $d$ on $X\times Y$ by applying the
triangle inequality and using that $\big(\sum_{i=1}^4 a_i\big)^p \le 4^p
\sum_{i=1}^4 a_i^p$ for numbers $a_i \ge 0$ with $1 \le i \le 4$. We find
\begin{align}
  d(x_1, y_1, x_2, y_2)^p
  &\le
  \big( d_X(x_1, x_0) + d_X(x_2, x_0) + d_Y(y_1, y_0) + d_Y(y_2, y_0) \big)^p \\
  &\le
  4^p \big( d_X(x_1, x_0)^p + d_X(x_2, x_0)^p + d_Y(y_1, y_0)^p + d_Y(y_2, y_0)^p\big) \\
  &=: \delta(x_1, y_1, x_2, y_2).
\end{align}
Let $\pi^*_n$ denote an (arbitrary) optimal transport plan between $\gamma_n$
and $\mu_n\otimes\nu_n$, meaning $\td(\gamma_n) = \pi_n^*c$.
Applying the previous inequality yields
\begin{equation}
  \pi^*_n\,c
  \le
  \pi^*_n\,d^p
  \le
  \pi^*_n\,\delta
  =
  2\cdot 4^p \big(\mu_n\,d_X(\cdot, x_0)^p + \nu_n\,d_Y(\cdot, y_0)^p\big),
\end{equation}
where the right-hand side converges as $n\to\infty$, since we have established
that $\mu_n$ and $\nu_n$ converge $p$-weakly. As the cost $c$ is continuous,
we can use the stability result in \parencite[Theorem~5.20]{villani2008} and
conclude that there is a subsequence $\pi_{n_k}^*$ of optimal transport plans
weakly converging to some optimal $\pi^*\in\CC(\gamma, \mu\otimes\nu)$, such
that $\td(\gamma) = \pi^* c$. Exploiting that
\begin{equation}
  \pi_{n_k}^* \delta
  \to
  2\cdot 4^p \big(\mu\,d_X(\cdot, x_0)^p + \nu\,d_Y(\cdot, y_0)^p\big)
  =
  \pi^* \delta
  < \infty
\end{equation}
as $k\to\infty$, we can apply Lemma~\ref{lem:dominated-uniform-integrability}
with $f = \delta$ and $g = c$ to finally conclude
\begin{equation}
  \bar\tau
  =
  \lim_{n\to\infty} \pi_n^*\,c
  =
  \lim_{k\to\infty} \pi_{n_k}^*c
  =
  \pi^*c
  =
  \td(\gamma),
\end{equation}
which completes the proof.
\end{proof}

\begin{proof}[Proof of Theorem~\ref{thm:continuity-metric-power}]
  Since $d$ satisfies the triangle inequality, so does $\rho = \ot_c^{1/p}$
  (see \cite[Definition~6.1]{villani2008}). Furthermore, the stated inequality
  is trivial if the right-hand side is $\infty$. Thus, we can assume $\rho(\gamma, \gamma') < \infty$
  and $\rho(\mu\otimes\nu, \mu'\otimes\nu') < \infty$. Then also $\rho(\gamma',
  \mu\otimes\nu) \le \rho(\gamma', \mu'\otimes\nu') + \rho(\mu\otimes\nu, \mu'\otimes\nu')
  < \infty$. Now, the result follows from applying the reverse triangle
  inequality:
  \begin{align}
    \big| \td(\gamma)^{1/p} - \td(\gamma')^{1/p} \big|
    &=
    \big| \rho(\gamma, \mu\otimes\nu) - \rho(\gamma', \mu'\otimes\nu') \big| \\
    &\le
    \big| \rho(\gamma, \mu\otimes\nu) - \rho(\gamma', \mu\otimes\nu) \big|
    +
    \big| \rho(\gamma', \mu\otimes\nu) - \rho(\gamma', \mu'\otimes\nu') \big| \\
    &\le
    \rho(\gamma, \gamma') + \rho(\mu\otimes\nu, \mu'\otimes\nu'). \qedhere
  \end{align}
\end{proof}

We write $\|\cdot\|_\infty$ to denote the $\sup$-norm of a real valued function
and use the convention $0/0 = 1$ in the following statement.

\begin{proposition}{varying costs}{varying-costs}
  Let $X$ and $Y$ be Polish spaces and let $c$ and $c_n$ be cost functions
  on $X\times Y$ that satisfy $\| c/c_n - 1 \|_\infty \to 0$ as $n\to\infty$. 
  Let $(\gamma_n)_{n\in\NN}$ be a sequence in $\PP(X\times Y)$ and $\gamma \in
  \PP(X\times Y)$ such that $\td_c(\gamma) < \infty$. Then
  \begin{alignat}{2}
    \lim_{n\to\infty} \td_c(\gamma_n) &= \td_c(\gamma)
    \qquad&&\text{implies}\qquad
    \lim_{n\to\infty} \td_{c_n}(\gamma_n)= \td_c(\gamma).
  \end{alignat}
\end{proposition}
\begin{proof}
  We set $a_n = \max\big(\|c/c_n - 1\|_\infty, \|c_n/c - 1\|_\infty\big)$ and
  observe $a_n \to 0$ as $n\to\infty$ by the assumption of uniform convergence.
  Applying Lemma~\ref{lem:uniform-costs} below, we can control the deviation of
  $\td_{c_n}$ from $\td_c$ and find, as $n\to\infty$,
  \begin{align}
    |\td_{c_n}(\gamma_n) - \td_c(\gamma)|
    &\le
    |\td_{c_n}(\gamma_n) - \td_c(\gamma_n)| + |\td_{c}(\gamma_n) - \td_c(\gamma)| \\
    &\le
    a_n(1 + a_n)\,\td_c(\gamma_n) + |\td_{c}(\gamma_n) - \td_c(\gamma)|
    \to 0.
    \qedhere
  \end{align}
\end{proof}

\begin{lemma}{}{uniform-costs}
  Let $X$ be a Polish space and let $c_1$ and $c_2$ be cost functions
  on $X$ that satisfy
  \begin{equation}
    \max\big(\|c_1/c_2 - 1\|_\infty, \|c_2/c_1 - 1\|_\infty\big) \le a
  \end{equation}
  under the convention $0/0 = 1$, where $\|\cdot\|_\infty$ is the sup norm and
  $a > 0$. Then for all $\mu,\nu\in\PP(X)$ with $\ot_{c_2}(\mu, \nu) < \infty$,
  we find $\ot_{c_1}(\mu, \nu) \le (1+a)\ot_{c_2}(\mu, \nu) < \infty$ and
  \begin{equation}
    |\ot_{c_1}(\mu, \nu) - \ot_{c_2}(\mu, \nu)|
    \le
    a\, \max\big(\ot_{c_1}(\mu, \nu), \ot_{c_2}(\mu, \nu)\big)
    \le
    a(1+a) \, \ot_{c_2}(\mu, \nu).
  \end{equation}
\end{lemma}
\begin{proof}[Proof of Lemma~\ref{lem:uniform-costs}]
  We first note that $c_1 \le (1+a) c_2$, which implies $\ot_{c_1}(\mu, \nu) \le
  (1+a)\ot_{c_2}(\mu, \nu)$. 
  Furthermore, if $\pi^*_1$ and $\pi^*_2$ denote optimal transport plans
  between $\mu$ and $\nu$ under the costs $c_1$ and $c_2$, then
  \begin{equation}
    \ot_{c_1}(\mu, \nu) - \ot_{c_2}(\mu, \nu)
    =
    \pi^*_1\,c_1 - \pi^*_2\,c_2
    \le
    \pi^*_2 |c_1 - c_2| \le a\,\pi^*_2\,c_2 = a\,\ot_{c_2}(\mu, \nu)
  \end{equation}
  with an analogous result for $\ot_{c_2}(\mu, \nu) - \ot_{c_1}(\mu, \nu)$, which
  establishes the claims of the lemma.
\end{proof}

\paragraph{Convolutions.}

In the following, a general statement about the behavior of the optimal
transport cost under convolution is formulated and proved
(Lemma~\ref{lem:ot-convolution}). This result is then applied to prove
Theorem~\ref{thm:convolution}.

\begin{lemma}{}{ot-convolution}
  Let $X$ be a Polish vector space and $c(x_1, x_2) = h(x_1 - x_2)$ for
  $x_1,x_2\in X$ a translation invariant cost function on $X$. For any
  probability measures $\mu$, $\nu$, and $\kappa$ in $\PP(X)$,
  \begin{equation}
    \ot_c(\mu * \kappa, \nu * \kappa) \le \ot_c(\mu, \nu)
  \qquad\text{and}\qquad
    \ot_c(\mu, \mu * \kappa) \le \kappa h.
  \end{equation}
\end{lemma}

\begin{proof}[Proof of Lemma~\ref{lem:ot-convolution}]
  To prove the first inequality, we reach for the dual formulation
  \eqref{eq:ot-definition-dual} of optimal transport. For any continuous and
  bounded potential $f\colon X \to \RR$, we define $f_\kappa(x) = \kappa\, f(x
  + \cdot)$ for $x\in X$. It is easy to check that $f_\kappa$ is again
  continuous and bounded. If $g$ is another potential such that $f\oplus g \le
  c$, we find
  \begin{equation}
    f_\kappa(x_1) + g_\kappa(x_2)
    =
    \kappa \big(f(x_1 + \cdot) + g(x_2 + \cdot)\big)
    \le
    \int\! c(x_1 + y, x_2 + y)\,\kappa(\dif y)
    =
    c(x_1, x_2)
  \end{equation}
  for any $x_1, x_2\in X$ due to the translation invariance of $c$. Therefore,
  $f_\kappa \oplus g_\kappa \le c$. This implies
  \begin{equation}
    (\mu * \kappa)\,f + (\nu * \kappa)\,g
    = 
    \mu f_\kappa + \nu g_\kappa
    \le
    \sup \mu f' + \nu g' = \ot_c(\mu, \nu),
  \end{equation}
  where the supremum is taken over continuous and bounded potentials $f'$ and $g'$
  with $f' \otimes g' \le c$. Since $f$ and $g$ were arbitrary, $\ot_c(\mu * \kappa,
  \nu * \kappa) \le \ot_c(\mu, \nu)$ follows. Note that this arguments holds even
  if $\ot_c(\mu, \nu) = \infty$, so there are no restrictions on $\mu,
  \nu\in\PP(X)$.

  For the upper bound in the second result, we construct an explicit transport
  plan $\pi$ between $\mu$ and $\mu * \kappa$. It is defined by $\pi f = \int
  f(x_1, x_1+x_2)\,\mu(\dif x_1) \kappa(\dif x_2)$ for any measurable map
  $f\colon X\times X \to [0, \infty)$. It is straightforward to check that
  $\pi\in\CC(\mu, \mu * \kappa)$, and we can conclude
  \begin{equation}
    \ot_c(\mu, \mu * \kappa) \le \pi c
    =
    \int c(x_1, x_1 + x_2)\,\mu(\dif x_1)\,\kappa(\dif x_2)
    = \kappa h,
  \end{equation}
  where we made use of the translation invariance $c(x_1, x_2) = h(x_1 - x_2)$
  and the symmetry $c(x_1, x_2) = c(x_2, x_1)$ of $c$ for any $x_1, x_2\in X$.
  Again, the argument stays valid even if $T_c(\mu, \mu * \kappa) = \infty$,
  so the results holds for any $\mu\in\PP(X)$.
\end{proof}

\begin{proof}[Proof of Theorem~\ref{thm:convolution}]
  One can easily check that $\gamma * \kappa \in \CC(\mu * \kappa_X, \nu
  * \kappa_Y)$ and $(\mu\otimes\nu)*\kappa = (\mu * \kappa_X) \otimes (\nu
  * \kappa_Y)$. Therefore, $\td(\gamma * \kappa) = \ot_c\big(\gamma * \kappa,
  (\mu\otimes\nu) * \kappa\big)$ and Lemma~\ref{lem:ot-convolution} can be
  applied, which yields the first result. For the second result, we
  recall Theorem~\ref{thm:continuity-metric-power} and use the second part
  of Lemma~\ref{lem:ot-convolution} to find
  \begin{equation}
    \td(\gamma)^{1/p} - \td(\gamma * \kappa)^{1/p}
    \le
    \ot_c(\gamma, \gamma * \kappa)^{1/p}
    + \ot_c\big(\mu\otimes\nu, (\mu\otimes\nu)*\kappa\big)^{1/p}
    \le 
    2 (\kappa h)^{1/p},
  \end{equation}
  which finishes the proof.
\end{proof}

\paragraph{Upper bounds and marginal transport dependency.}
This segment establishes the upper bounds in
Proposition~\ref{prop:upper-bound} and \ref{prop:upper-bound-ot} and
contains proofs for Theorem~\ref{thm:marginal-transport-dependence},
\ref{thm:marginal-transport-dependence-measurable}, and
Proposition~\ref{prop:marginal-consistency} concerning the marginal transport
dependency.

\begin{proof}[Proof of Proposition~\ref{prop:upper-bound}]
  It is evident that $\cmin$ is a cost function (non-negative, symmetric, lower
  semi-continuous). Furthermore, the second inequality in \eqref{eq:upper-bound}
  follows trivially. To prove the first inequality, we construct a coupling
  $\pi_2\in\CC(\gamma, \mu\otimes\nu)$ that aims to prevent either horizontal or
  vertical movements, depending on which of the associated marginal costs is
  larger. We therefore define the set
  \begin{equation*}
    S = \big\{(x_1,y_1,x_2,y_2)\,\big|\,c_X(x_1,x_2)
      \ge
      c_Y(y_1,y_2)\big\} \subset (X\times Y)^2,
  \end{equation*}
  and note that $S$ is symmetric under exchanging $(x_1,y_1)$ and $(x_2,y_2)$
  due to the symmetry of the costs $c_X$ and $c_Y$. We write $R$ to denote
  the complement of $S$, which is also symmetric. Next, we introduce the
  function $r\colon (X\times Y)^2 \to (X\times Y)^2$ given by
  \begin{equation*}
    r(x_1,y_1,x_2,y_2) = \begin{cases}
      (x_1,y_1, x_1, y_2) &\text{if}~(x_1,y_1,x_2,y_2)\in S, \\
      (x_2,y_2, x_1, y_2) &\text{else}.
    \end{cases}
  \end{equation*}
  The proof is completed once we show that the coupling defined by $\pi_2
  = r_\#(\gamma\otimes\gamma)$ has the correct marginals, meaning
  $\pi_2\in\CC(\gamma, \mu\otimes\nu)$, and that $\pi_2 c$, which is an upper bound
  for $\td(\gamma)$, is in turn upper bounded by $(\gamma\otimes\gamma)\,\cmin$.
  The second marginal $\mu\otimes\nu$ is an immediate consequence of the
  definition of $r$ and $\pi_2$. To check the first marginal, we consider an
  arbitrary positive and measurable function $f\colon
  X \times Y\to\RR$. Then, if $q$ denotes the two first components of $r$,
  \begin{align}
    \int f(x_1,y_1)\,\dif\pi_2(x_1,y_1,x_2,y_2) 
    &=
    \int f\big(q(x_1,y_1,x_2,y_2)\big)\,\dif\gamma(x_1,y_1)\,\dif\gamma(x_2,y_2) \\
    &=
    \int \mathbbm{1}_{S}(x_1,y_1,x_2,y_2)\,f(x_1,y_1)\,\dif\gamma(x_1,y_2)\,\dif\gamma(x_2,y_2) \\
    &\qquad\quad +
      \int \mathbbm{1}_{R}(x_1,y_1,x_2,y_2)\,f(x_2,y_2)\,\dif\gamma(x_1,y_1)\,\dif\gamma(x_2,y_2) \\
    &= \int f(x_1,y_1)\,\dif\gamma(x_1,y_1) = \gamma f,
  \end{align}
  where we swapped the roles of $(x_1,y_1)$ and $(x_2,y_2)$ to establish the
  third equality. This is permissible due to the symmetry of $S$ (and $R$).
  Similarly, we observe
  \begin{align}
    \td(\gamma)
    &\le
    \pi_2 c \\
    &=
    \int_S c(x_1,y_1,x_1,y_2)\,\dif\gamma(x_1,y_1)\,\dif\gamma(x_2,y_2)
    +
    \int_{R} c(x_2,y_2,x_1,y_2)\,\dif\gamma(x_1,y_1)\,\dif\gamma(x_2,y_2) \\
    &\le 
    \int_S c_Y(y_1, y_2)\,\dif\gamma(x_1,y_1)\,\dif\gamma(x_2,y_2)
    +
    \int_{R} c_X(x_1,x_2)\,\dif\gamma(x_1,y_1)\,\dif\gamma(x_2,y_2) \\
    &=
    \int \min\!\big(c_X(x_1,x_2), c_Y(y_1,y_2)\big)\,\dif\gamma(x_1,y_1)\,\dif\gamma(x_2,y_2)
    = (\gamma\otimes\gamma)\,\cmin,
  \end{align}
  where we bounded $c$ by $c_Y$ and $c_X$ via condition
  \eqref{eq:marginal-cost-sup}.
\end{proof}

\begin{proof}[Proof of Proposition~\ref{prop:upper-bound-ot}]
  In order to prove this result, we first introduce an alternative
  characterization of the measurability of probability kernels that is employed
  in \textcite{villani2008}.

  \begin{lemma}{measurability of kernels}{kernel-measurability}
    Let $X$ be Polish and let $(\mu_\omega)_{\omega\in \Omega}
    \subset\PP(X)$ be a family of probability measures indexed in a measurable
    space $(\Omega, \mathcal{F})$. Then 
    \begin{equation}
      \omega \mapsto \mu_\omega(A)~\text{is measurable for all Borel sets}~A\subset X
      \quad\Longleftrightarrow\quad
      \omega \mapsto \mu_\omega~\text{is measurable},
    \end{equation}
    where $\PP(X)$ is equipped with the Borel $\sigma$-algebra with respect to the
    topology of weak convergence of measures.
  \end{lemma}
  \begin{proof}
    Theorem~17.24 in \textcite{kechris2012} asserts that the Borel
    $\sigma$-algebra of $\PP(X)$ is generated by functions $r_A\colon
    \PP(X)\to[0,1]$ of the form $\nu \mapsto \nu(A)$ for Borel sets $A\subset X$.
    This means that
    $\mathcal{G} = \big\{r_A^{-1}(B)\,\big|\,B\subset[0, 1]~
    \text{Borel and}~A\subset X~\text{Borel}\big\}$
    is a generator of the Borel $\sigma$-algebra of $\PP(X)$. In particular, each
    $r_A$ is measurable.

    Thus, if $\mu\colon\omega \mapsto \mu_\omega$ is measurable,
    then $\omega \mapsto (r_A\circ\mu)(\omega) = \mu_\omega(A)$ is also measurable as
    composition of measurable functions. Conversely, if $\omega \mapsto \mu_\omega(A)$
    is measurable for each Borel set $A\subset X$, then $\mu^{-1}(G) \in
    \mathcal{F}$ for each $G\in\mathcal{G}$. Since $\mathcal{G}$ is a generator, this
    suffices to show that $\omega \mapsto \mu_\omega$ is measurable.
  \end{proof}

  We now return to the proof of Proposition~\ref{prop:upper-bound-ot}.
  Let $\pi^*_x\in\CC\big(\gamma(x, \cdot), \mu\big)$ be an optimal transport
  plan with respect to the base costs $c_Y$ for each $x\in X$. Corollary~5.22 in
  \textcite{villani2008} together with the continuity of $c_X$ and
  Lemma~\ref{lem:kernel-measurability} above guarantee that $\pi^*_x$ can be
  selected such that $(\pi^*_x)_{x\in X}$ is a probability kernel. We can thus
  define $\pi_3\in\PP(X\times Y)$ via
  \begin{equation}
    \dif\pi_3(x_1, y_1, x_2, y_2)
    =
    \pi^*_{x_1}(\dif y_1, \dif y_2)\,\delta_{x_1}(\dif x_2) \,\mu(\dif x_1).
  \end{equation}
  It can easily be checked that the marginals of $\pi_3$ match $\gamma$ and
  $\mu\otimes\nu$. Using condition~\eqref{eq:marginal-cost-sup}, we find
  \begin{equation}
    \td(\gamma)
    \le
    \pi_3 c
    \le
    \int \left(\int c_Y(y_1, y_2)\,\pi^*_{x}(\dif y_1,\dif y_2)\right)\,\mu(\dif x)
    =
    \int (\pi_{x}^*\,c_Y) \,\mu(\dif x),
  \end{equation}
  which shows the first inequality of the proposition. To assert the second
  inequality, one just has to note that $\pi_x^*\,c_Y \le \big(\gamma(x,
  \cdot)\otimes\nu\big)\,c_Y$ for each $x\in X$ by construction of $\pi^*_x$ as
  optimal plan.
\end{proof}

\begin{proof}[Proof of Theorem~\ref{thm:marginal-transport-dependence}]
  We begin by defining a set that contains all vertical movements along the
  fibers $\{x\}\times Y$, given by $S = \{(x, y_1, x, y_2)\,|\, x\in X, \,y_1,
  y_2 \in Y\} \subset (X\times Y)^2$. Due to the definition of $c_\infty$, it is
  evident that 
  \begin{equation}\label{eq:td-c-infinity}
    \td_{c_\infty}(\gamma)
    =
    \inf_{\substack{\pi\in\CC(\gamma, \mu\otimes\nu)\\ \pi(S) = 1}} \pi c_Y.
  \end{equation}
  We use this characterization to show that the equality $\td_{c_\infty}(\gamma)
  = \td^Y_{c_Y}(\gamma)$ holds. First, we define $f\colon S\to X\times Y^2$
  via $f(x, y_1, x, y_2) = (x, y_1, y_2)$. For each $\pi$ that is feasible
  in the infimum in \eqref{eq:td-c-infinity}, we furthermore define $\pi_x
  = (f_\#\pi)(x, \cdot, \cdot)\in\PP(Y\times Y)$ for $x\in X$. One can
  check that $\pi_x \in \CC\big(\gamma(x, \cdot), \nu\big)$ holds
  $\mu$-almost surely due to the (almost sure) uniqueness property of
  disintegrations. If $\pi_x^*\in\CC\big(\gamma(x, \cdot), \nu\big)$ are 
  a measurable selection of optimal transport plans for the problem
  $\ot_{c_Y}\big(\gamma(x, \cdot), \nu\big)$ as in the proof of
  Proposition~\ref{prop:upper-bound-ot}, we observe
  \begin{equation}
    \pi c_Y
    =
    \int_S c_Y\,\dif\pi
    =
    \int c_Y\,\dif(f_\#\pi)
    =
    \int (\pi_x c_Y)\,\mu(\dif x)
    \ge
    \int (\pi_x^* c_Y)\,\mu(\dif x)
    =
    \td^Y_{c_Y}(\gamma).
  \end{equation}
  Taking the infimum on the left-hand side and applying \eqref{eq:td-c-infinity}
  implies $\td_{c_\infty}(\gamma) \ge \td^Y_{c_Y}(\gamma)$. Since we already know
  $\td_{c_\infty}(\gamma) \le \td^Y_{c_Y}(\gamma)$ by
  Proposition~\ref{prop:upper-bound-ot}, this proves equality.
  
  It is left to show that $\td_{c_\alpha}(\gamma) \to \td_{c_\infty}(\gamma)$ as
  $\alpha\to\infty$. Since $c_\infty \ge c_\alpha$ for all $\alpha > 0$, we
  know that $\td_{c_\infty}(\gamma) \ge \td_{c_\alpha}(\gamma)$ and it is
  accordingly sufficient to show $\liminf_{\alpha\to\infty}
  \td_{c_\alpha}(\gamma) = \td_{c_\infty}(\gamma)$.
  We thus fix $\bar\tau = \liminf_{\alpha\to\infty}
  \td_{c_\alpha}(\gamma)$ and consider a diverging sequence
  $(\alpha_n)_{n\in\NN}\subset (0, \infty)$ such that
  $\smash{\td_{c_{\alpha_n}}(\gamma) \to \bar\tau}$. Due to Prokhorov's theorem
  and the
  existence of optimal transport plans, we can assume that there are couplings
  $\pi_n^*$ and $\pi_\infty$ in $\CC(\gamma, \mu\otimes\nu)$ that satisfy
  $\td_{c_{\alpha_n}}(\gamma) = \pi_n^* c_{\alpha_n}$ and $\pi_n^* \weak
  \pi_\infty$.
  We use these properties to lead the assumption $\infty \ge
  \td_{c_\infty}(\gamma) > \bar\tau$ to a contradiction. If this assumption were
  true, observe that for any $k > 0$,
  \begin{equation}
    \infty > \bar\tau
    =
    \lim_{n\to\infty} \pi_n^* c_{\alpha_n}
    \ge
    \liminf_{n\to\infty} \alpha_n \pi_n^* c_X
    \ge
    k\,\liminf_{n\to\infty} \pi_n^* c_X
    \ge
    k\,\pi_\infty c_X,
  \end{equation}
  where the final inequality follows from the lower semi-continuity of the
  mapping $\pi \mapsto \pi h$ under weak convergence for lower semi-continuous
  integrands $h$ that are bounded from below \parencite[Proposition
  7.1]{santambrogio2015}.
  Since $k$ is arbitrary, we conclude $\pi_\infty c_X = 0$, which implies
  $\pi_\infty(S) = 1$ due to the positivity of $c_X$. Employing representation
  \eqref{eq:td-c-infinity} of $\td_{c_\infty}(\gamma)$, we find the
  contradiction
  \begin{equation}
    \td_{c_\infty}(\gamma) > \bar\tau
    =
    \lim_{n\to\infty}\pi_n^* c_{\alpha_n}
    \ge
    \liminf_{n\to\infty} \pi_n^* c_Y
    \ge
    \pi_\infty c_Y
    \ge
    \td_{c_\infty}(\gamma).
  \end{equation}
  This establishes $\bar\tau = \td_{c_\infty}(\gamma)$ and finishes the proof.
\end{proof}

\begin{proof}[Proof of Theorem~\ref{thm:marginal-transport-dependence-measurable}]
  We begin the proof by showing the following auxiliary statement.
  \begin{lemma}{}{marginal-transport-dependence}
    Let $X$ be a Polish space and $c$ a positive continuous cost function on
    $X$. For $\mu, \nu \in \PP(X)$ with $\supp\,\mu \subset\supp\,\nu$ and
    $\ot_c(\mu, \nu) < \infty$, it holds that
    \begin{equation}
      \ot_c(\mu, \nu) = (\mu \otimes \nu)c
      \qquad
      \Longleftrightarrow
      \qquad
      \mu = \delta_{x}~\text{for some $x\in X$}.
    \end{equation}
  \end{lemma}
  \begin{proof}
    The implication from right to left is trivial, since $\mu\otimes\nu$ is the
    only feasible coupling if $\mu$ is a point mass. To show the reverse
    direction, we assume that $\mu \neq \delta_x$ for any $x\in X$ and show that
    $\ot_c(\mu, \nu) = (\mu\otimes\nu)c$ is impossible. First, we pick two
    distinct points $x_1 \neq x_2$ from the support of $\mu$. By assumption,
    these points also lie in the support of $\nu$. To show that $\mu \otimes
    \nu$ cannot be an optimal transport plan, it is sufficient to show that
    $\supp\,(\mu\otimes\nu) = \supp\,\mu \times \supp\,\nu$ is not $c$-cyclically
    monotone \parencite{villani2008}. This is easy to see, since for $(x_1,
    x_2), (x_2, x_1) \in \supp\,(\mu\otimes\nu)$, we find
    \begin{equation}
      c(x_1, x_2) + c(x_2, x_1) > c(x_1, x_1) + c(x_2, x_2) = 0
    \end{equation}
    due to the positivity of $c$.
  \end{proof}

  Returning to the proof of
  Theorem~\ref{thm:marginal-transport-dependence-measurable}, we note that
  $\td^Y(\gamma)$ and $(\nu\otimes\nu)\,c_Y$ are equal iff
  \begin{equation}\label{eq:proof-marginal-transport-dependence}
    \int \ot_{c_Y}\big(\gamma(x, \cdot), \nu\big)\,\mu(\dif x)
    = 
    \int \big(\gamma(x, \cdot)\otimes\nu\big)\,c_Y\,\mu(\dif x).
  \end{equation}
  Evidently, these two values are the same if $\gamma = (\id,\varphi)_\#\mu$,
  since this implies $\gamma(x, \cdot) = \delta_{\varphi(x)}$ for $\mu$-almost all
  $x\in X$. So it only remains to show that equality in
  \eqref{eq:proof-marginal-transport-dependence} implies $\gamma
  = (\id,\varphi)_\#\mu$ for some $\mu$-almost surely defined measurable function
  $\varphi\colon X \to Y$.

  Since the right-hand side in \eqref{eq:proof-marginal-transport-dependence} is
  by assumption finite and $0 \le \ot_{c_Y}\big(\gamma(x, \cdot), \nu\big)
  \le \big(\gamma(x, \cdot)\otimes\nu\big)c_Y$ holds for each $x\in X$, we find
  that equality in \eqref{eq:proof-marginal-transport-dependence} can only hold
  if
  \begin{equation}
    \ot_{c_Y}\big(\gamma(x, \cdot), \nu\big)
    =
    \big(\gamma(x, \cdot)\otimes\nu\big)c_Y
  \end{equation}
  for $\mu$-almost all $x\in X$. As $\supp\,\gamma(x, \cdot) \subset \supp\,\nu$
  also holds for $\mu$-almost all $x\in X$ (see below), we can apply
  Lemma~\ref{lem:marginal-transport-dependence} and find that $\gamma(x, \cdot)
  = \delta_{\varphi(x)}$ for a $\mu$-almost surely defined function $\varphi$.
  The measurability of $\varphi$ follows from the measurability of the maps
  $x\mapsto\gamma(x, A)$ for all Borel sets $A\subset Y$, since
  $x\in\varphi^{-1}(A)$ is equivalent to $\gamma(x, A) = 1$
  for all $x\in X$ for which $\varphi$ is defined by the construction above.

  A brief argument to see that $\supp\,\gamma(x, \cdot) \subset \supp\,\nu$ for
  $\mu$-almost all $x\in X$ goes as follows: note that $\supp\,\gamma
  \subset \supp\,(\mu\otimes\nu) = \supp\,\mu \times \supp\,\nu$ and write
  \begin{equation}
    1
    =
    \int \dif\,\gamma
    =
    \int \ind_{\supp \nu}(y)\,\gamma(\dif x, \dif y)
    =
    \int \ind_{\supp \nu}(y)\,\gamma(x, \dif y)\,\mu(\dif x),
  \end{equation}
  which proves that $\gamma(x, \supp\,\nu) = 1$ for $\mu$-almost all $x\in X$.
\end{proof}

\begin{proof}[Proof of Proposition~\ref{prop:marginal-consistency}]
  We first note that $\td_{c_\alpha}(\gamma) < \infty$ for all $0 < \alpha \le
  \infty$ due to the finite $p$-th moment of $\gamma$. Next, we observe
  \begin{align*}
    &\EE\left[\big|\td_{c_{\alpha_n}}(\hat\gamma_n)^{1/p}
      - \td^Y_{c_Y}(\gamma)^{1/p}\big|\right] \leq \\
    &\hspace{3cm}
     \EE\left[\big|\td_{c_{\alpha_n}}(\hat\gamma_n)^{1/p}
      - \td_{c_{\alpha_n}}(\gamma)^{1/p}\big|\right]
      + \big|\td_{c_{\alpha_n}}(\gamma)^{1/p}- \td^Y_{c_Y}(\gamma)^{1/p}\big|.
  \end{align*}
  By Theorem~\ref{thm:marginal-transport-dependence}, the second summand
  converges to zero as $n\to \infty$. The first summand can be controlled by
  Theorem~\ref{thm:continuity-metric-power} and bound \eqref{eq:td-lipschitz},
  yielding
  \begin{equation}
    \EE\left[\big|\td_{c_{\alpha_n}}(\hat\gamma_n)^{1/p}
      - \td_{c_{\alpha_n}}(\gamma)^{1/p}\big|\right] 
    \leq
    3\,\EE\left[\ot_{c_{\alpha_n}}(\gamma, \hat\gamma_n)^{1/p}\right]
    \leq
    3\,\alpha_n\,\EE\left[\ot_{c_{1}}(\gamma, \hat\gamma_n)^{1/p}\right].
  \end{equation}
\end{proof}

\paragraph{Contracting couplings and maps.}
We next provide proofs and statements that were omitted in our work on
contractions in Section~\ref{sec:contractions}. This includes the proof of our
main results, Theorem~\ref{thm:contractions} and
\ref{thm:contractions-deterministic}, as well as the formulation of two
auxiliary statements (Lemma~\ref{lem:continuous-contraction} and
\ref{lem:uniform-extension}), which simplify Theorem~\ref{thm:contractions} if
sufficient regularity is imposed on the involved costs.

\begin{proof}[Proof of Theorem~\ref{thm:contractions}.]
  We first show the second part. Assuming that $\td(\gamma)
  = (\nu\otimes\nu)\,c_Y  < \infty$, we use the upper bound in
  Proposition~\ref{prop:upper-bound} to conclude $(\gamma\otimes\gamma)\,c_{XY}
  = (\nu\otimes\nu)\,c_Y$, which can be stated as
  \begin{equation}
    \int\!\min\!\big(c_X(x_1, x_2), c_Y(y_1, y_2)\big)\, \dif(\gamma\otimes\gamma)(x_1, y_1, x_2, y_2)
    =
    \int\! c_Y(y_1, y_2)\,\dif(\nu\otimes\nu)(y_1, y_2).
  \end{equation}
  This implies $(\gamma\otimes\gamma)(c_Y \le c_X) = 1$, since the left-hand
  side in the equality above would otherwise be strictly smaller than the
  right-hand side.
  Recalling that $c_X = h \circ k_X$ and $c_Y = h \circ d_Y$ for a strictly
  increasing $h$, we find $(\gamma\otimes\gamma)(d_Y \le k_X) = 1$, meaning that
  $\gamma$ is almost surely contracting.

  To prove the first statement, we show that $\td(\gamma) \ge
  (\nu\otimes\nu)\,c_Y$, which is sufficient to assert equality due to
  Proposition~\ref{prop:upper-bound}.
  Let $\pi\in\CC(\mu, \nu)$ be arbitrary. We define $\bar\lambda\in\PP\big((X\times
  Y)^2\times Y\big)$ via the relation $\dif\bar\lambda(x_1,y_1,x_2,y,y_2)
  = \gamma(x_2, \dif y)\,\dif\pi(x_1, y_1, x_2, y_2)$. One can readily establish
  that integrating out $y_2\in Y$ yields a measure $\lambda \in \CC(\gamma,
  \gamma) \subset\PP\big((X\times Y)^2\big)$. In particular, $\supp\,\lambda \subset
  \supp\,\gamma\times\supp\,\gamma$, which implies $\lambda(d_Y \le k_X) = 1$ by
  the assumption that $\gamma$ is contracting on its support. Combining this
  insight with the strict monotonicity of $h$ and the triangle inequality for
  $d_Y$, we find
  \begin{align}
    \pi c
    &=
    \int h\big(k_X(x_1,x_2) + d_Y(y_1,y_2)\big)\,
    \dif\bar\lambda(x_1, y_2, x_1, y, y_2) \\
    &\ge
    \int h\big(d_Y(y_1,y) + d_Y(y_1,y_2)\big)\,
    \dif\bar\lambda(x_1, y_2, x_1, y, y_2) \\
    &\ge
    \int h\big(d_Y(y,y_2)\big)\,\dif\bar\lambda(x_1, y_2, x_1, y, y_2) \\
    &= 
    \int c_Y(y,y_2)\,\dif\nu(y)\,\dif\nu(y_2) = (\nu\otimes\nu)\,c_Y.
  \end{align}
  One can check that the independence in the final line follows from the
  definition of $\bar\lambda$. Taking the infimum over
  $\pi\in\CC(\gamma,\mu\otimes\nu)$ now yields the desired result.
\end{proof}

\begin{lemma}{}{continuous-contraction}
  Let $X$ and $Y$ be Polish spaces and let $k_X$ and $d_Y$ be continuous cost
  functions on $X$ and $Y$. Then $\gamma$ is contracting on its support iff it
  is almost surely contracting (with respect to $k_X$ and $d_Y$).
\end{lemma}
\begin{proof}[Proof of Lemma~\ref{lem:continuous-contraction}]
  It is clear that a contracting coupling $\gamma$ is almost surely contracting
  since the set $(X\times Y)^2 \setminus (\supp\,\gamma)^2$ is a null set for
  $\gamma\otimes\gamma$. Let $\gamma$ therefore be almost surely contracting.
  If $k_X$ and $d_Y$ are continuous and there are $(x_1, y_1), (x_2, y_2) \in
  \supp\,\gamma$ with $d_Y(y_1, y_2) > k_X(x_1, x_2)$, then we also find an open
  neighbourhood $U\subset (X\times Y)^2$ of $(x_1, y_1, x_2,
  y_2)\in(\supp\,\gamma)^2 = \supp\,(\gamma\otimes\gamma)$ where this inequality
  is true. By the definition of the support, we conclude
  $(\gamma\otimes\gamma)(U) > 0$ and $\gamma$ thus fails to be almost surely
  contracting.
\end{proof}

\begin{lemma}{uniform extension}{uniform-extension}
  Let $X$ be Polish and $(Y, d_Y)$ a Polish metric space, and let $k_X\colon X\times
  X \to [0, \infty)$ be a continuous cost function on $X$. Any function
  $\tilde \varphi\colon D \to Y$ defined on a subset $D\subset X$ that satisfies
  \begin{equation}\label{eq:uniform-continuity}
    d_Y\big(\tilde\varphi(x_1), \tilde\varphi(x_2)\big)
    \le
    k_X(x_1, x_2)
  \end{equation}
  for all $x_1, x_2\in D$ can uniquely be extended to a function $\varphi\colon
  \bar D \to Y$ on the closure $\bar D$ of $D$ that satisfies
  \eqref{eq:uniform-continuity} for all $x_1, x_2\in\bar D$. In particular,
  $\varphi$ is continuous.
\end{lemma}
\begin{proof}[Proof of Lemma~\ref{lem:uniform-extension}]
  Let $x\in \bar D \setminus D$ and $(x_n)_{n\in\NN}$ be a sequence in $D$
  converging to $x$. Set $y_n = \tilde\varphi(x_n)$ for $n\in\NN$ and observe that
  $d_Y(y_n, y_m) \le k_X(x_n, x_m)$ for all $n,m\in\NN$. In particular, the sequence
  $(y_n)_n$ is Cauchy: if it were not Cauchy, there would exist an $\epsilon
  > 0$ and values $n_r, m_r \ge r$ for each $r \in\NN$ such that $d_Y(y_{n_r},
  y_{m_r}) \ge \epsilon$. However, observing $\lim_{r\to\infty} k_X(x_{n_r},
  x_{m_r}) \to 0$ due to continuity of $k_X$ leads this to a contradiction.
  Consequently, the sequence $(y_n)_{n\in\NN}$ converges to a unique limit $y
  \in Y$ due to the completeness of $(Y, d_Y)$.
  
  The limit point $y$ does not depend on the chosen sequence: if $(x'_n)_n$ is
  another sequence converging to $x$, and $y'$ is the corresponding limit in
  $Y$, continuity of $d_Y$ and $k_X$ make sure that
  \begin{equation}
    0
    \le
    d_Y(y, y')
    =
    \lim_{n\to\infty} d_Y(y_n, y'_n)
    \le
    \lim_{n\to\infty} k_X(x_n, x'_n)
    = 
    0.
  \end{equation}
  Therefore, a well-defined extension of $\tilde\varphi$ to $\bar D$ exists.
  For any $x, x' \in \bar D$, this extension $\varphi$ satisfies
  \begin{equation}
    d_Y\big(\varphi(x), \varphi(x')\big)
    =
    \lim_{n\to\infty} d_Y\big(\tilde\varphi(x_n), \tilde\varphi(x'_n)\big)
    \le
    \lim_{n\to\infty} k_X(x_n, x'_n)
    =
    k_X(x, x')
  \end{equation}
  as $n\to\infty$ for any sequence $(x_n, x'_n)_{n\in\NN} \subset D\times D$
  that converges to $(x, x')\in \bar D \times \bar D$.
\end{proof}

\begin{proof}[Proof of Theorem~\ref{thm:contractions-deterministic}]
  Let $\varphi\colon X\to Y$ be $\mu$-almost surely contracting such that
  $\gamma = (\id, \varphi)_\#\mu$. By a change of variables,
  \begin{equation}\label{eq:change-variables-contraction}
    (\gamma\otimes\gamma)(d_Y \le k_X)
    =
    (\mu\otimes\mu)(d_\varphi \le k_X)
    =
    1,
  \end{equation}
  where $d_\varphi(x_1, x_2) = d_Y\big(\varphi(x_1), \varphi(x_2)\big)$ for any
  $x_1, x_2\in X$.
  Therefore, $\gamma$ is almost surely contracting. Since $k_X$ and
  $d_Y$ are continuous, we can apply Lemma~\ref{lem:continuous-contraction} to
  conclude that $\gamma$ is contracting on its support. By
  Theorem~\ref{thm:contractions}, $\td(\gamma) = (\nu\otimes\nu)\,c_Y$ follows.

  For the reverse direction, we consult the upper bound in
  Proposition~\ref{prop:upper-bound-ot} to find that $\td(\gamma)
  = (\nu\otimes\nu)\,c_Y$ implies $\td^Y(\gamma) = (\nu\otimes\nu)\,c_Y$.
  Consequently, Theorem~\ref{thm:marginal-transport-dependence-measurable}
  establishes the existence of a measurable function $\varphi\colon X\to Y$ with
  $\gamma = (\id, \varphi)_\#\mu$. It is left to show that $\varphi$ is
  $\mu$-almost surely contracting. According to
  Lemma~\ref{lem:continuous-contraction} and Theorem~\ref{thm:contractions},
  $\gamma$ is contracting on its support. We can therefore consider the set $A
  = (\id, \varphi)^{-1}(\supp\,\gamma)$ and conclude that both $\mu(A)
  = \gamma(\supp\,\gamma) = 1$ and
  \begin{equation}
    d_Y\big(\varphi(x_1), \varphi(x_2)\big)
    \le
    d_X(x_1, x_2)
    \qquad\text{for all}\qquad
    x_1, x_2\in A,
  \end{equation}
  since all tuples $\big(x, \varphi(x)\big)\in X\times Y$ for $x\in A$ are
  elements of the support of $\gamma$.
\end{proof}

\paragraph{Properties of the transport correlation.}
We next present proofs of Proposition~\ref{prop:transport-correlation} to
\ref{prop:isometric-transport-correlation} regarding basic properties of the
transport correlation. Most of the claimed properties are direct consequences of
previously established results. Additional arguments are mainly required for
properties 5 and 6 in Proposition~\ref{prop:transport-correlation} and property
2 in Proposition~\ref{prop:isometric-transport-correlation}. Recall from
Section~\ref{sec:coefficients} that we focus on additive costs of the form
\begin{equation}
  c = (\alpha d_X + d_Y)^p
\end{equation}
on Polish metric spaces $(X, d_X)$ and $(Y, d_Y)$ for $\alpha, p > 0$.

\begin{proof}[Proof of Proposition~\ref{prop:transport-correlation}]
  Properties 1 and 2, which characterize when $\tc_\alpha(\gamma)$ equals $0$
  and $1$, follow from Theorem~\ref{thm:independence} in
  Section~\ref{sec:definitions} and from
  Theorem~\ref{thm:contractions-deterministic} in
  Section~\ref{sec:contractions}. The invariance in property 3 is
  a consequence of Proposition~\ref{prop:invariance}.
  If $\gamma$ is restricted to a set with fixed marginal $p^Y_\#\gamma = \nu$,
  then $(\nu\otimes\nu)\,d_Y^p$ is constant and convexity of $\gamma \mapsto
  \tc_\alpha(\gamma)^p = \td(\gamma) / \smash{(\nu\otimes\nu)\,d_Y^p}$ is guaranteed
  by convexity of $\td$ as stated in Proposition~\ref{prop:convexity}. This
  shows property 4.

  We now turn to the continuity in property 5. Let $\nu_n$ and $\nu$ denote the
  respective second marginals of $\gamma_n$ and $\gamma$.
  For convenience, we write $\delta_n =
  = (\nu_n\otimes\nu_n)\,d_Y^p$. Our first goal is to show that $p$-weak
  convergence of $\gamma_n$ to $\gamma$ implies $\delta_n \to \delta
  = (\nu\otimes\nu)\,d^p_Y$ as $n\to\infty$.
  To do so, we fix some $y_0\in Y$ and define the function
  \begin{equation}
    f(y_1, y_2)
    =
    2^p\big(d_Y(y_1, y_0)^p + d_Y(y_0, y_2)^p\big)
  \end{equation}
  for $y_1, y_2\in Y$.  Noting that $(a + b)^p \le 2^p(a^p + b^p)$ for $a, b \ge
  0$, we find $d_Y^p \le f$ by application of the triangle inequality.
  Consequently,
  \begin{equation}
    \delta_n
    =
    (\nu_n\otimes\nu_n)\,d_Y^p
    \le
    (\nu_n\otimes\nu_n)\,f
    =
    2^{p+1}\nu_n d_Y(\cdot, y_0)^p
    \to
    2^{p+1}\nu d_Y(\cdot, y_0)^p
  \end{equation}
  as $n\to\infty$. The convergence in this inequality follows from the fact that
  $\nu_n$ converges $p$-weakly if $\gamma_n$ converges $p$-weakly, which was
  shown in the proof of Proposition~\ref{prop:continuity}.
  Reaching back to Lemma~\ref{lem:dominated-uniform-integrability} (setting $g
  = d_Y^p$, $f = f$, and $\mu_n = \nu_n\otimes\nu_n$), we conclude that
  $\delta_n$ indeed converges to $\delta$ as $n\to\infty$.
  Next, since $c = (\alpha d_X + d_Y)^p \le \max(1, \alpha)^p (d_X + d_Y)^p$, we
  can apply Proposition~\ref{prop:continuity} and find $\td_c(\gamma_n) \to
  \td_c(\gamma)$. Setting $c_n = (\alpha_n d_X + d_Y)^p$, it is straightforward
  to see that $\|c / c_n - 1\|_\infty \to 0$ due to $\alpha_n \to \alpha > 0$
  for $n\to\infty$, which lets us use Proposition~\ref{prop:varying-costs} to
  obtain
  \begin{equation}
    \lim_{n\to\infty}\tc_{\alpha_n}(\gamma_n)^p
    =
    \lim_{n\to\infty}\frac{\td_{c_n}(\gamma_n)}{\delta_n}
    =
    \frac{\td_{c}(\gamma)}{\delta}
    =
    \tc_\alpha(\gamma)^p.
  \end{equation}
  Finally, the monotonicity of $\alpha \mapsto \tc_\alpha(\gamma)$ in property
  6 is trivial, since $c$ increases with $\alpha$. The concavity of
  $\alpha\mapsto\tc_\alpha(\gamma)^p$ for $p \le 1$ and fixed $\gamma\in\CC(\mu,
  \nu)$ with marginals $\mu\in\PP(X)$ and $\nu\in\PP(Y)$ follows from the fact
  that pointwise infima over concave functions are again concave. Indeed, we
  have that
  \begin{equation}
    \tc_\alpha(\gamma)^p
    =
    \frac{1}{(\nu\otimes\nu)d_Y^p}\,
    \inf_{\pi\in\CC(\gamma, \mu\otimes\nu)} \pi (\alpha d_X + d_Y)^p,
  \end{equation}
  where the mapping $\alpha \mapsto \pi (\alpha d_X + d_Y)^p$ is concave for
  any fixed $\pi$ if $p \le 1$.
\end{proof}

\begin{proof}[Proof of Proposition~\ref{prop:marginal-transport-correlation}]
  Due to Theorem~\ref{thm:marginal-transport-dependence}, results established
  for $\td$ under generic lower semi-continuous costs $c$ also hold for the
  marginal transport dependency with $\alpha = \infty$. Therefore, properties 1,
  3, and 4 follow from the general results stated in
  Theorem~\ref{thm:independence}, Proposition~\ref{prop:invariance}, and
  Proposition~\ref{prop:convexity}. Note that we can allow dilatations $f_Y$ in
  property 3 (instead of just isometries), since we normalize by
  $(\nu\otimes\nu)\,d_Y^p$ in
  Definition~\ref{def:marginal-transport-correlation} of $\tc_\infty$, which
  neutralizes the dilatation factor $\beta > 0$.  Property 2 relies on the
  specific cost structure for $\alpha = \infty$ and was derived separately in
  Theorem~\ref{thm:marginal-transport-dependence-measurable}.
\end{proof}

\begin{proof}[Proof of Proposition~\ref{prop:isometric-transport-correlation}]
  Properties 1, 4, and 5 follow in the same way as in the proof of
  Proposition~\ref{prop:transport-correlation}. The symmetry property 6 is
  trivial and directly visible from the definition of $\tc_*$.
  Property 3 is a consequence of
  the general invariance of the transport dependence under isometries
  (Proposition~\ref{prop:invariance}). We can extend this to dilatations $f_X$
  and $f_Y$ since we divide by the respective diameters in
  Definition~\ref{def:isometric-transport-correlation} of $\tc_*$, which
  nullifies any scaling factors.

  Regarding property 2, let $\gamma\in\CC(\mu, \nu)$ for $\mu\in\PP(X)$ and
  $\nu\in\PP(Y)$.
  We equip the spaces $X$ and $Y$ with the scaled metrics
  \begin{equation}
    \tilde d_X = d_X / \big((\mu\otimes\mu)\,d_X^p\big)^{1/p}
    \qquad\text{and}\qquad
    \tilde d_Y = d_Y / \big((\nu\otimes\nu)\,d_Y^p\big)^{1/p}.
  \end{equation}
  According to Theorem~\ref{thm:contractions-deterministic} (combined with
  Lemma~\ref{lem:uniform-extension}), a value of $\tc_*(\gamma) = 1$ is
  equivalent to there being a contraction $\varphi$ from $(\supp\,\mu, \tilde
  d_X)$ to $(Y, \tilde d_Y)$ that satisfies $(\id, \varphi)_\#\mu = \gamma$. This
  in particular means $\varphi_\#\mu = \nu$.  Defining $\tilde d_\varphi(x_1,
  x_2) = \tilde d_Y\big(\varphi(x_1), \varphi(x_2)\big)$ for any $x_1, x_2\in
  \supp\,\mu$, we know that $\tilde d_\varphi \le d_X$ since $\varphi$ is
  a contraction.
  By a change of variables, we calculate
  \begin{equation}
    (\mu\otimes\mu)\,\tilde d_\varphi^p
    = 
    (\nu\otimes\nu)\,\tilde d_Y^p
    = 1 =
    (\mu\otimes\mu)\,\tilde d_X^p
  \end{equation}
  and assert $\tilde d_\varphi = \tilde d_X$ to hold
  $(\mu\otimes\mu)$-almost surely. Since $\tilde d_X$ and $\tilde d_\varphi$ are
  continuous, we can conclude that $\smash{\varphi\colon (\supp\,\mu, \tilde
  d_X) \to (Y, \tilde d_Y)}$ is an isometry. Thus, $\varphi$ is a dilatation
  under the metrics $d_X$ and $d_Y$. Note that its dilatation factor $\beta$ is
  uniquely given by $\alpha_*$ defined in equation \eqref{eq:isometric-alpha}.
  Due to the symmetry of the setting, the same arguments also hold for $\psi$
  instead of $\varphi$ by exchanging the roles of $X$ and $Y$.
\end{proof}

\section{Lower complexity adaptation}
\label{app:lca}

In the following, we derive the property of lower complexity adaptation (LCA)
for the estimators proposed in Section~\ref{sec:estimation}.
We restrict our analysis to bounded and continuous cost functions $c\colon (X
\times Y)^2 \to [0, 1]$.
The proof strategy is inspired by \textcite{hundrieser2022empirical}, but we use
adapted arguments to exploit the additional randomness introduced by the
sampling procedure. As preparation, we revisit the dual formulation for the
optimal transport cost and introduce some tools from empirical process theory. 

\paragraph{Preliminaries.}

The duality theory of optimal transport crucially depends on the notion of the
$c$-conjugacy of functions. We present a selection of definitions and
results on this topic, adapted to the above setting.
For a given $g \colon X\times Y \to \RR$ bounded from above, we define its
\emph{$c$-transform} to be
\begin{equation}\label{eq:ctransform}
  g^c(x, y) = \inf_{(x', y')\in X\times Y} c(x, y, x', y') - g(x', y')
\end{equation}
for any $(x,y)\in X\times Y$.
Each function that can be written as a $c$-transform is called
\emph{$c$-concave}. Since $c$ is continuous, all $c$-concave functions are
measurable. The set of (standardized) $c$-concave functions is denoted by
\begin{equation}
  \FF_c
  =
  \big\{g^c\, | \, g \colon X\times Y \to \RR, \,\sup\nolimits_{x\in X} g(x) = 0\big\}.
\end{equation}
Based on the boundedness of the cost function, it is possible to show that both
$f$ and $f^c$ are absolutely bounded by $1$ for each $f\in\FF_c$.
Since we need to restrict $\FF_c$ to the support of the coupling $\gamma\in\PP(X
\times Y)$ in order to exploit the LCA property, we also define the
domain-restricted function classes
\begin{equation}
  \FF_c(\gamma) = \{ f|_{\supp\,\gamma}\,|\,f\in\FF_c \}
\end{equation}
for $\gamma\in\PP(X\times Y)$. A function $f\in\FF_c(\gamma)$ can be
$c$-transformed via definition \eqref{eq:ctransform} by taking the infimum over
$(x', y')\in\supp\,\gamma$.
For any $\eta_1,\eta_2\in\PP(X\times Y)$, it now follows from standard optimal
transport theory \parencite[Theorem~5.10]{villani2008} that strong duality in
the form
\begin{align}
  \ot_c(\eta_1, \eta_2)
  =
  \sup_{f\in\FF_c} \eta_1 f + \eta_2 f^c
\end{align}
holds. If $\supp\,\eta_1 \subset \supp\,\gamma$, then the consistent behavior of
optimal transport under restriction of the base spaces (see, e.g.,
\cite[Theorem~5.19]{villani2008}, or \cite[Lemma~3]{staudt2022uniqueness}) even
allows us to conclude
\begin{align}
  \ot_c(\eta_1, \eta_2)
  =
  \sup_{f\in\FF_c(\gamma)} \eta_1 f + \eta_2 f^c.
\end{align}
The fact that we can optimize over $\FF_c(\gamma)$ instead of $\FF_c$ is what
gives rise to the lower complexity adaptation of statistical optimal transport,
since it enables us to reason about the optimal transport problem in terms of
the support of $\gamma$.

Another tool we require is the \emph{uniform metric entropy} of a class $\FF$ of
real-valued functions. It is defined as the logarithm of the uniform covering
number $\covering(\epsilon, \FF)$, which denotes the minimal number of
sets with diameter $2\epsilon$ required to cover $\FF$ in the uniform norm for $\epsilon > 0$.
Our results are formulated under the assumption that the uniform metric entropy
of $\FF_c(\gamma)$ is upper bounded by
\begin{equation}\label{eq:lca-entropy-bound}
  \log\covering\big(\epsilon, \FF_c(\gamma)\big)
  \lesssim
  \epsilon^{-k}
\end{equation}
for some $k > 0$.
Then, we will deduce upper bounds on the convergence rates of the form
\begin{equation}\label{eq:lca-rate}
  r_k(n)
  :=
  \begin{cases}
    n^{-1/2} & \text{if}~k < 2 \\
    n^{-1/2}\log(n) & \text{if}~k = 2 \\
    n^{-1/k} & \text{if}~k > 2.
  \end{cases}
\end{equation}
Clearly, the value of $k$ depends on the properties of $c$ as well as the
support of $\gamma$.
For example, if $X$ and $Y$ are smooth manifolds and $c$ is twice continuously
differentiable while the support of $\gamma$ is a compact subset of a smooth
submanifold of $X\times Y$ of dimension $s\in\NN$, one can derive $k = 2/s$. As
a general rule of thumb, we find $k = s / \alpha$ if $\alpha \in (0, 2]$ denotes
the \emph{Hölder-smoothness} of $c$ and $s$ the \emph{intrinsic dimension} of
$\gamma$ (unfortunately, higher degrees of smoothness of $c$ than $\alpha = 2$
cannot be exploited).
For precise definitions and results along those lines, see Section~3 in
\textcite{hundrieser2022empirical}.

\paragraph{LCA for various estimators.}

We start with proving the LCA property of the plug-in estimator
$\tau(\hat\gamma_n)$, where $\hat\gamma_n$ is the empirical measure
$\frac{1}{n}\sum_{i=1}^n \delta_{(\xi_i, \zeta_i)}$ for $(\xi_i,
\zeta_i)_{i=1}^n\sim\gamma^{\otimes n}$. According to inequality
\eqref{eq:lca-metric}, this estimator exhibits lower complexity adaptation under
additive metric costs. The result below confirms this property also for general
bounded costs.

\begin{theorem}{LCA, product estimator}{lcaproduct}
  Let $X$ and $Y$ be Polish spaces, $c\colon X\times Y \to [0, 1]$
  continuous, and $\gamma \in \PP(X\times Y)$. If bound
  \eqref{eq:lca-entropy-bound} holds for $k > 0$, then
  \begin{equation}
    \Exp\,|\tau(\hat\gamma_n) - \tau(\gamma)|
    \lesssim
    r_k(n).
  \end{equation}
\end{theorem}

\begin{proof}
  Let $\mu\in\PP(X)$ and $\nu\in\PP(Y)$ be the marginal distributions of
  $\gamma$. We introduce the random variables $(\xi, \zeta) \sim
  \mu\otimes\nu$, assumed to be independent of the samples $(\xi_1, \zeta_1),
  (\xi_2, \zeta_2), \ldots, (\xi_n, \zeta_n)$. Since $\supp\,\hat\gamma_n
  \subset \supp\,\gamma$ with probability $1$, we can use the same argument as
  in the proof of \textcite[Theorem~2.2]{hundrieser2022empirical}, to derive
  \begin{equation}
    \big| T_c(\hat\gamma_n, \hat\mu_n\otimes\hat\nu_n) - T_c(\gamma, \mu\otimes\nu) \big|
    \le
    \sup_{f\in\FF_c(\gamma)} | (\hat\gamma_n - \gamma) f |
    +
    \sup_{f\in\FF_c(\gamma)} | (\hat\mu_n\otimes\hat\nu_n - \mu\otimes\nu) f^c|.
  \end{equation}
  Note that the class $\FF_c(\gamma)$ equipped with the sup norm is separable
  since it has finite covering numbers, which means that the right hand side is
  measurable.
  The expectation over the first term on the right hand side can be treated like
  in \textcite{hundrieser2022empirical} and we find
  \begin{equation}
    \Exp \sup_{f\in\FF_c(\gamma)} |(\hat\gamma_n - \gamma) f|
    \lesssim
    r_k(n)
  \end{equation}
  with a universal constant. To address the second term, we write
  $\FF_c^c(\gamma) = \{f^c\,|\,f\in\FF_c(\gamma)\}$ and apply Lemma~2.1 in
  \textcite{hundrieser2022empirical} to establish that $\covering\big(\epsilon,
  \FF_c(\gamma)\big) = \covering\big(\epsilon, \FF_c^c(\gamma)\big)$. Noting
  that $(\xi_i, \zeta_j)$ is equal to $(\xi, \zeta_j)$ in distribution for $i
  \neq j$, we conclude
  \begin{align}
    \Exp \sup_{f\in\FF_c(\gamma)} |(\hat\mu_n\otimes\hat\nu_n - \mu\otimes\nu) f^c|
    &=
   \Exp \sup_{g\in\FF_c^c(\gamma)}
      \left| \frac{1}{n^2} \sum_{i,j=1}^n g(\xi_i, \zeta_j) - \Exp\,g(\xi, \zeta) \right| \\
    &\le
   \frac{1}{n}\sum_{i=1}^n ~\Exp \sup_{g\in\FF_c^c(\gamma)}
      \left| \frac{1}{n} \sum_{j=1}^n g(\xi_i, \zeta_j) - \Exp\,g(\xi, \zeta) \right| \\
    &\le
   \Exp \sup_{g\in\FF_c^c(\gamma)}
      \left| \frac{1}{n} \sum_{j=1}^n g(\xi, \zeta_j) - \Exp\,g(\xi, \zeta) \right| + \frac{2}{n},
  \end{align}
  where we have exploited that $\|g\|_\infty \le 1$ in the last inequality in order
  to replace the diagonal terms corresponding to $j = i$ in the sum.
  Consequentially, the random variables $g(\xi, \zeta_j)$ for $j\in\{1, \ldots,
  n\}$ are i.i.d.\ distributed and bounded. Denoting $Z_g
  = \frac{1}{\sqrt{n}} \sum_{i=1}^n g(\xi, \zeta_j)$, we can apply
  Hoeffding's inequality to derive sub-Gaussianity
  \begin{equation}
    \Prob\big(|Z_g - \Exp\,Z_g| \ge t\big) \le 2\,e^{-t^2/2\|g\|_\infty^2}
  \end{equation}
  with respect to the uniform norm $\|\cdot\|_\infty$. Therefore, we can resort to
  chaining arguments like in Theorem~5.22 in \textcite{wainwright2019high} to deduce
  \begin{align}
    &\Exp \sup_{f\in\FF_c(\gamma)} |(\hat\mu_n\otimes\hat\nu_n - \mu\otimes\nu) f^c| \\
    &\hspace{8em}\le
     O\big(n^{-1}\big) + \frac{1}{\sqrt{n}}~\Exp \sup_{g\in\FF_c^c(\gamma)} |Z_g - \Exp\,Z_g| \\
    &\hspace{8em}\lesssim
     O\big(n^{-1}\big) + \frac{1}{\sqrt{n}}\int_{\delta/4}^1 \sqrt{\log\,\covering\big(\epsilon, \FF_c(\gamma)\big)}\,\dif\epsilon + \frac{1}{\sqrt{n}}~\Exp\sup_{\substack{g,g'\in\FF_c^c(\gamma)\\ \|g-g'\|_\infty \le \delta}} |Z_{g} - Z_{g'}| \\
    &\hspace{8em}\lesssim
     O\big(n^{-1}\big) + \frac{1}{\sqrt{n}}\int_{\delta/4}^1 \sqrt{\log\,\covering\big(\epsilon, \FF_c(\gamma)\big)}\,\dif\epsilon + \delta
  \end{align}
  for any $\delta > 0$. We remark that the statement in
  \textcite{wainwright2019high} is not explicitly formulated for the absolute
  value of $Z_g - \Exp\,Z_g$. However, this can easily be adapted by enlarging
  the function class the supremum is taken over by the functions
  $\{-g\,|\,g\in\FF_c^c(\gamma)\}$, which at most increases the covering number
  by a factor of two.
  Theorem~2.2 in \textcite{hundrieser2022empirical} now shows that the
  right hand side in the display above is bounded by a multiple of $r_k(n)$ for
  a suitable choice of $\delta$, which establishes the claim of the theorem.
\end{proof}

We now turn to the sampling based estimators proposed in
Section~\ref{sec:estimation}. Our proof strategy allows for both sampling with
and without replacement. The proof works by exploiting sub-Gaussianity of sums
with randomly sampled indices.

\begin{lemma}{}{sampling-bound}
  Let $a = (a_{ij})_{i,j=1}^n \in \RR^{n\times n}$ for $n \in \NN$ and consider
  the random variable $Z = \frac{1}{\sqrt{n}}\sum_{i=1}^n a_{I_i J_i}$ with
  (random) indices $I = (I_i)_{i=1}^n$ and $J = (J_i)_{i=1}^n$ that satisfy
  either
  \begin{enumerate}[topsep=1ex, parsep=0.5ex, font={\normalfont\color{green!30!black}},label=\theenumi)]
    \item
      $I_i = i$ for all $i\in\{1, \ldots, n\}$ and
      $J$ is sampled from $\{1, \ldots, n\}$, or
    \item
      $(I, J)$ is sampled from $\{1, \ldots, n\}^2$,
  \end{enumerate}
  where the sampling is uniform with or without replacement and the roles of $I$
  and $J$ in case 1) can also be exchanged. Then there is a constant $b > 0$
  only depending on the sampling scheme such that
  \begin{equation}\label{eq:sampling-bound}
    \Prob\big(|Z - \Exp\,Z| \ge t\big) \le 2\,e^{-t^2/2b\|a\|_\infty^2},
  \end{equation}
  where $\|a\|_\infty = \max_{i,j} |a_{ij}|$ denotes the absolute value of the
  largest entry.
\end{lemma}
\begin{proof}
  In case of sampling with replacement, $Z$ is the (normalized) sum of $n$
  independent random variables with values in $[-\|a\|_\infty, \|a\|_\infty]$.
  Thus, the classical Hoeffding inequality suffices to derive inequality
  \eqref{eq:sampling-bound} for $b = 1$. For sampling without replacement in both components
  (scenario 2), the Hoeffding inequality holds as well (see, e.g.,
  \cite{bardenet2015concentration}). Finally, sampling without replacement in
  only one component is, for example, treated in \textcite{chatterjee2007stein},
  where Proposition~1.1 implies the desired result after shifting and scaling
  the data (rough estimates show $b \le 12$).
\end{proof}

\begin{theorem}{LCA, sampling estimators}{lcasampling}
  Let $X$ and $Y$ be Polish spaces and $c\colon X\times Y \to [0, 1]$ be
  continuous. For $\gamma \in \CC(\mu, \nu)$ with $\mu\in\PP(X)$ and
  $\nu\in\PP(Y)$, assume that bound \eqref{eq:lca-entropy-bound} holds for $k
  > 0$. Define the estimator
  \begin{equation}
    (\mu\otimes\nu)_n
    =
    \frac{1}{n}\sum_{i=1}^n \delta_{(\xi_{I_i}, \zeta_{J_i})}
  \end{equation}
  with (random) indices $I = (I_i)_{i=1}^n$ and $J = (J_i)_{i=1}^n$ as in
  Lemma~\ref{lem:sampling-bound}, independent of the data. Then
  \begin{equation}
    \Exp\,\big|\ot_c\big(\hat\gamma_n, (\mu\otimes\nu)_n\big) - \tau(\gamma)\big|
    \lesssim
    r_k(n).
  \end{equation}
\end{theorem}

\begin{proof}
  Let $\hat\Exp$ denote the expectation with respect to the random sampling
  $(I,J)$ only. We will show that
  \begin{equation}\label{eq:conditional-rate}
    \hat\Exp\, \big|T_c\big(\hat\gamma_n, (\mu \otimes\nu)_n\big) - T_c(\hat\gamma_n, \hat\mu_n \otimes \hat\nu_n)\big| 
    \lesssim
    r_k(n)
  \end{equation}
  holds for all realizations of the data (with a common constant), which, in
  combination with Theorem~\ref{thm:lcaproduct}, suffices to show the claim.
  We first note
  \begin{equation}
    \hat{\Exp}\,(\mu\otimes\nu)_n = \hat\mu_n\otimes\hat\nu_n,
  \end{equation}
  which follows from the observation that $\hat\Exp\,(\mu\otimes\nu)_n$ places
  the same amount of mass on any point $(\xi_i, \zeta_j)_{i,j=1}^n$ due to the
  symmetry of the sampling procedure. Denoting $Z_g = \smash{\frac{1}{\sqrt{n}}}
  \sum_{i=1}^n g(\xi_{I_i}, \zeta_{J_i})$ for $g \in \FF_c^c(\gamma)$, we apply
  analogous arguments as in the proof of Theorem~\ref{thm:lcaproduct} to derive
  \begin{equation}
    \hat\Exp\, \big|T_c\big(\hat\gamma_n, (\mu \otimes\nu)_n\big) - T_c(\hat\gamma_n, \hat\mu_n \otimes \hat\nu_n)\big| 
    \le
    \frac{1}{\sqrt{n}}~\hat\Exp
      \sup_{g\in\FF_c^c(\gamma)} \left| Z_g - \hat\Exp\,Z_g \right|.
  \end{equation}
  According to Lemma~\ref{lem:sampling-bound}, the process
  $(Z_g)_{g\in\FF_c^c(\gamma)}$ is sub-Gaussian with respect to the uniform norm
  (scaled by a constant $b > 0$) when conditioned on the observations. Following
  the remaining argumentation in the proof of Theorem~\ref{thm:lcaproduct} (with
  $\Exp$ replaced by $\hat\Exp$), bound \eqref{eq:conditional-rate} follows,
  where we emphasize that the right hand side does not depend on the data anymore.
  Note also that the value of $b$, which amounts to a scaling in the radius of the
  covering, only affects \eqref{eq:conditional-rate} in the implicit constant.
\end{proof}

Finally, we want to highlight the fact that repeated sampling from the
observations does not deteriorate the convergence rate of the proposed sampling
estimators. In practice, combining several (independent) samples leads to
a reduced variability of the estimate. This effect is especially pronounced for
small to moderate values of $n$ (see Section~\ref{sec:applications}).

\begin{corollary}{LCA, repeated sampling}{lcarepeated}
  Let $X$ and $Y$ be Polish spaces and $c\colon X\times Y \to [0, 1]$ be
  continuous. For $\gamma \in \CC(\mu, \nu)$ with $\mu\in\PP(X)$ and
  $\nu\in\PP(Y)$, assume that bound \eqref{eq:lca-entropy-bound} holds for $k
  > 0$. For $j \in \{1, \ldots, r\}$ with $r \in\NN$, let
  $\smash{(\mu\otimes\nu)_{n, j}}$ be one of the sampling estimators defined in
  Theorem~\ref{thm:lcasampling} and set
  \begin{equation}
    (\mu\otimes\nu)_n^r
    =
    \frac{1}{r}\sum_{j=1}^r (\mu\otimes\nu)_{n, j}.
  \end{equation}
  Then
  \begin{equation}
    \Exp\,\big|\ot_c\big(\hat\gamma_n, (\mu\otimes\nu)_n^r\big) - \tau(\gamma)\big)\big|
    \lesssim
    r_k(n).
  \end{equation}
\end{corollary}

\begin{proof}
  With the same argument from the proof of
  \textcite[Theorem~2.2]{hundrieser2022empirical}, that has already been applied in
  Theorem~\ref{thm:lcaproduct} and \ref{thm:lcasampling} to bound the difference
  of optimal transport costs, we find
  \begin{align}
    \Exp\,\big|\ot_c\big(\hat\gamma_n, (\mu\otimes\nu)_n^r\big) - \ot_c\big(\hat\gamma_n, \hat\mu_n\otimes\hat\nu_n\big)\big|
    &\le
    \Exp\,\sup_{f\in\FF_c(\gamma)} \big|\big((\mu\otimes\nu)_n^r - \hat\mu_n\otimes\hat\nu_n\big)f^c\big| \\
    &\le
    \frac{1}{r}\sum_{j=1}^r
    \Exp\,\sup_{f\in\FF_c(\gamma)} \big|\big((\mu\otimes\nu)_{n,j} - \hat\mu_n\otimes\hat\nu_n\big)f^c\big| \\
    &\lesssim r_k(n),
  \end{align}
  where we made use of the triangle inequality of the absolute value and
  reused the proof of Theorem~\ref{thm:lcasampling} in the final step. By
  Theorem \ref{thm:lcaproduct}, which shows that $T_c(\hat\gamma_n,
  \hat\mu_n\otimes\hat\nu_n)$ approaches $\td(\gamma)$ with rate $r_k(n)$, we
  can conclude the claim.
\end{proof}

\section{Analytic computation of the transport dependency}
\label{app:explicit}

In this appendix, we present a simple setting where the transport dependency can
be calculated explicitly. Let $c_{r}$ for $r\in\NN$ denote the squared Euclidean
cost in $\RR^r$ and let $\delta\colon \RR^r \to \RR^{2r}$ be the mapping to the
diagonal, $\delta(x) = (x, x)$. Consider $\gamma = \delta_{\#}\mu \in \CC(\mu,
\nu)$ for $\mu = \nu = \smash{\bigotimes_{i=1}^r} \mu_i \in \PP(\RR^r)$, where
$\mu_i\in\PP(\RR)$ for each $1\le i \le r$. Expressed in random variables, this
joint distribution corresponds to the case $\xi = \zeta$. We assume that
$\tau(\gamma) = \tau_{c_{2r}}(\gamma) < \infty$. 
By construction, the support of $\gamma$ is
contained in the affine subspace 
\begin{equation*} 
  \delta\big(\RR^r\big) 
  = 
  \big\{\delta(x)\,|\,x\in\RR^{r}\big\}\subset \RR^{2r}, 
\end{equation*} 
and the orthogonal projection $p\,\colon \RR^{2r}\to \delta\big(\RR^r\big)$
onto this subspace is $p(y) = \delta(\tilde{y})$ with $\tilde{y}\in\RR^r$,
$\tilde{y}_i = (y_i + y_{r+i})/2$. 
According to \cite[Proposition~2.3]{hundrieser2022empirical}, we find that the
optimal transport cost in this scenario can be decomposed as
\begin{align}
  \tau(\gamma)
  &=
  \ot_{c_{2r}}(\gamma, \mu\otimes\mu) \\
  &=
  \ot_{c_{2r}}\big(\gamma, p_\#(\mu\otimes\mu)\big)
  + 
  \ot_{c_{2r}}\big(p_\#(\mu\otimes\mu), \mu\otimes\mu\big) \\
  &=
  \ot_{c_{2r}}\big(\gamma, p_\# (\mu\otimes\mu)\big)
  +
  \int \|y - p(y)\|^2 \,\dif(\mu\otimes\mu)(y).
\end{align}
Moreover, letting $\tilde\mu_i = t_\# (\mu_i \otimes\mu_i)$ with $t(a, b)
= (a+b)/2$ and setting $\tilde\mu = \smash{\bigotimes_{i=1}^r} \tilde\mu_i$,
it follows that $p_\#(\mu\otimes\mu) = \delta_\#\tilde\mu$. Observing
$c_{2r}\big(\delta(x), \delta(y)\big) = 2\,c_r(x, y)$ for any $x,y\in\RR^r$, we
now apply Lemma~\ref{lem:ot-invariance} to conclude
\begin{equation*} 
  \ot_{c_{2r}}\big(\gamma, p_\#(\mu\otimes\mu)\big)
  =
  \ot_{c_{2r}}(\delta_\# \mu, \delta_\# \tilde\mu)
  = 
  2\,\ot_{c_r}(\mu, \tilde\mu)
  =
  2\,\sum_{i=1}^{r}\ot_{c_1}(\mu_i, \tilde\mu_i).
\end{equation*} 
Combining the previous two equations yields
\begin{align*} 
  \tau(\gamma) 
  &= 
  \ot_{c_{2r}}\big(\gamma, p_\#(\mu\otimes\mu)\big) + \int \|y - p(y)\|^2 \,\dif(\mu\otimes\mu)(y) \\ 
  &= 
  2 \sum_{i=1}^r \ot_{c_1}(\mu_i, \tilde\mu_i) + \frac{1}{2} \sum_{i=1}^d \int (y_1 - y_2)^2 \,\dif(\mu_i\otimes\mu_i)(y_1, y_2), 
\end{align*} 
which reduces calculating $\tau(\gamma)$ to one-dimensional optimal transport
problems and integrals over $\mu_i\otimes\mu_i$. If $\mu_i
= \mathrm{Unif}[0, 1]$ for all $i$, for example, the involved quantities can
be calculated explicitly and we find 
\begin{equation*} 
  \tau(\gamma) 
  = 
  \frac{r}{60} + \frac{r}{12} = \frac{r}{10}. 
\end{equation*}

\clearpage

\section{Additional Simulations}
\label{app:simulations}

This appendix contains a range of figures that supplement
Section~\ref{sec:applications} and further illustrate the behaviour of the
transport correlation on noisy datasets.

\begin{figure}[t!]
  \centering
  \includegraphics[width=.8\textwidth]{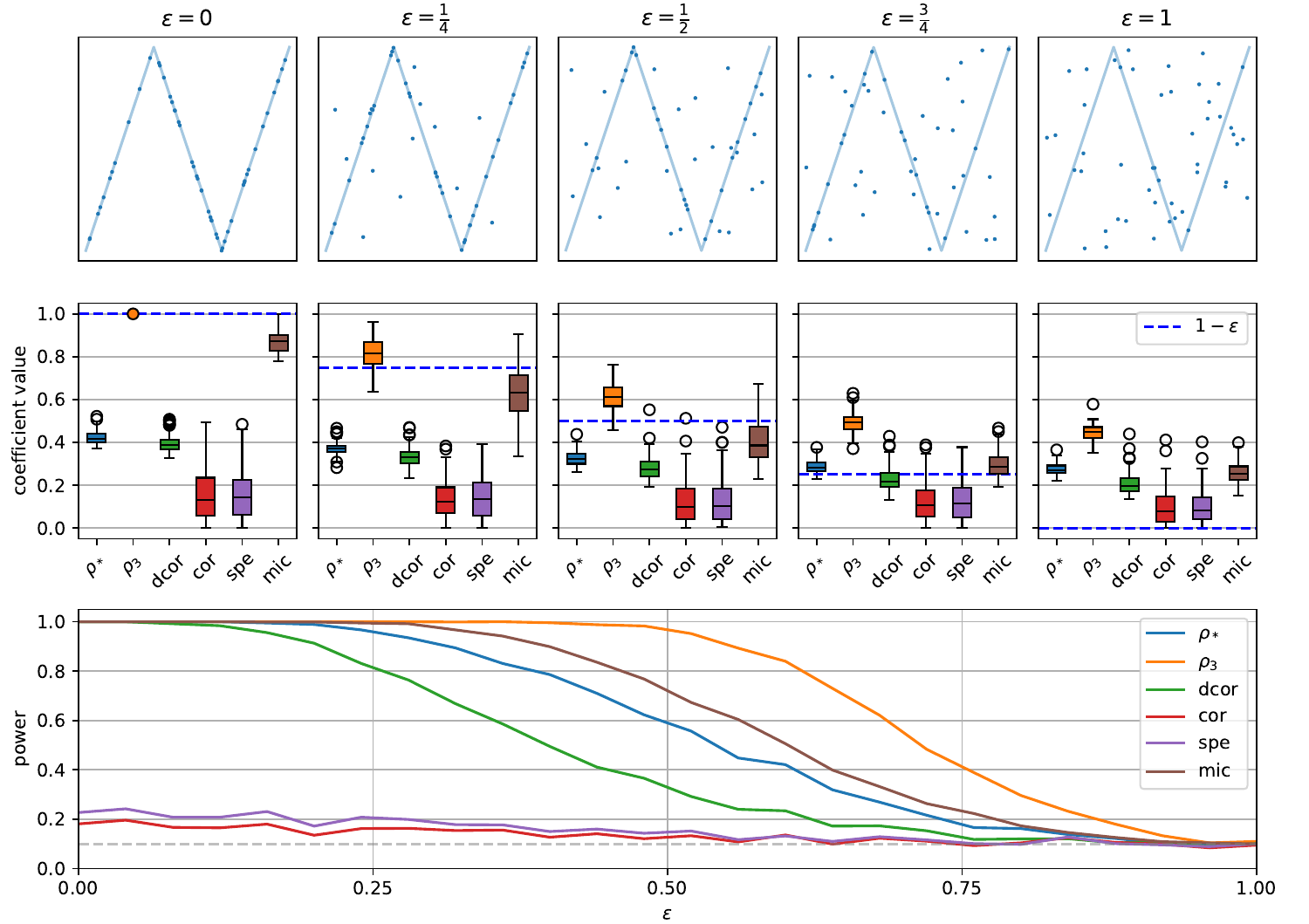}
  \caption{}
  \label{fig:appendix01}
  \end{figure}

\begin{figure}[t!]
  \centering
  \includegraphics[width=.8\textwidth]{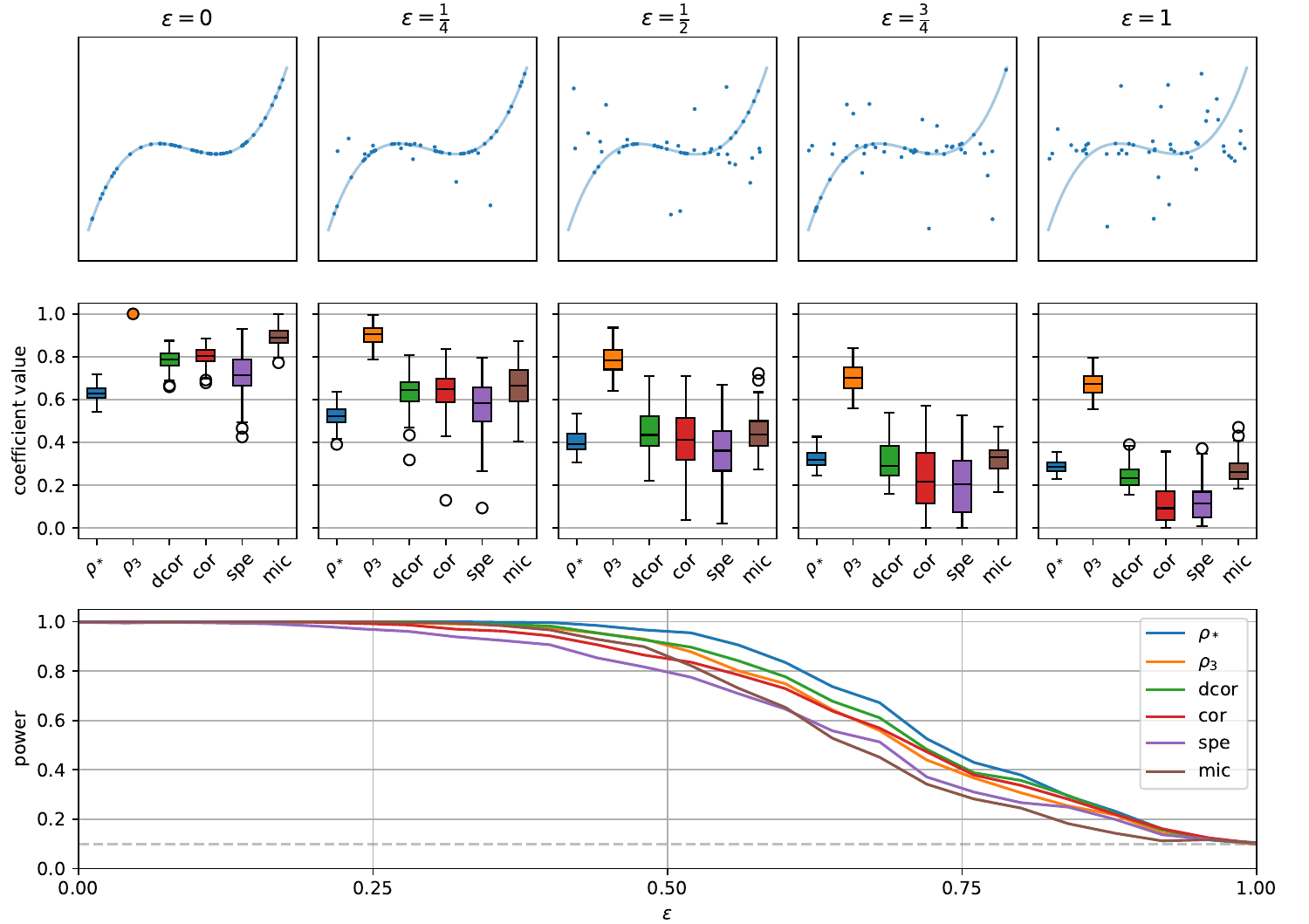}
  \caption{}
  \label{fig:appendix03}
\end{figure}

\begin{figure}
  \centering
  \includegraphics[width=.8\textwidth]{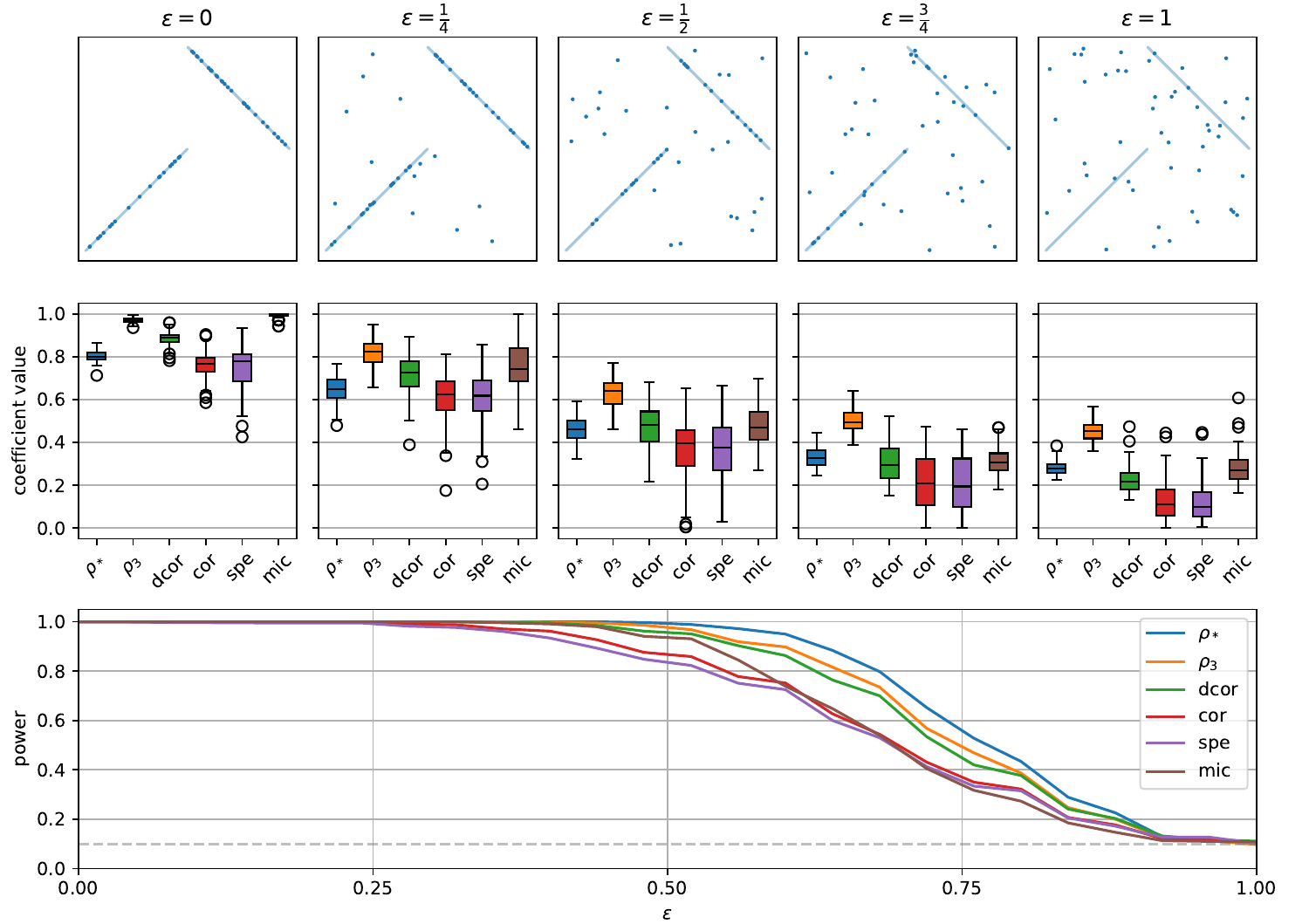}
  \label{fig:appendix04}
  \caption{}
\end{figure}

\begin{figure}
  \centering
 \includegraphics[width=.8\textwidth]{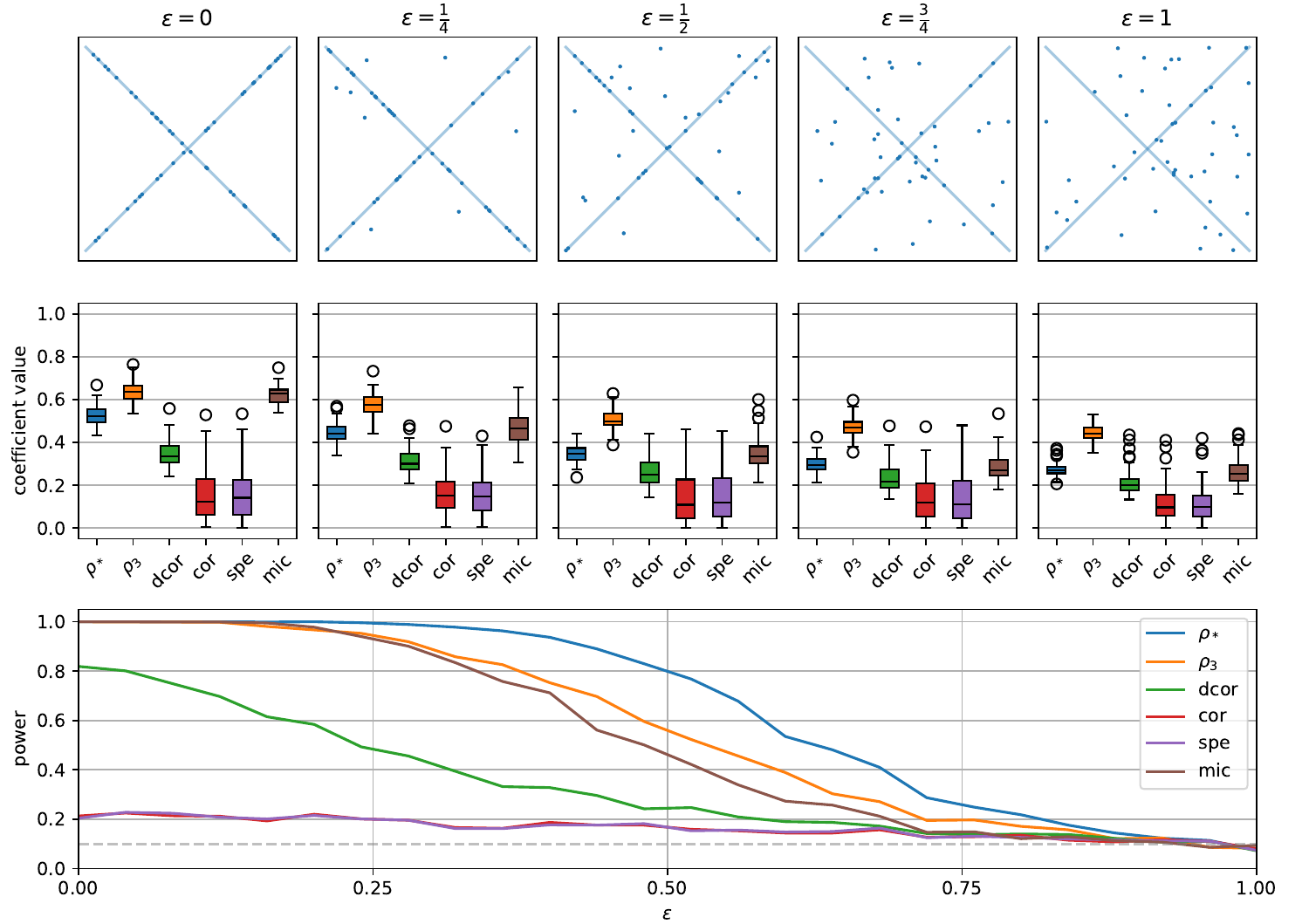}
 \label{fig:appendix05}
 \caption{}
\end{figure}

\begin{figure}
  \centering
  \includegraphics[width=.8\textwidth]{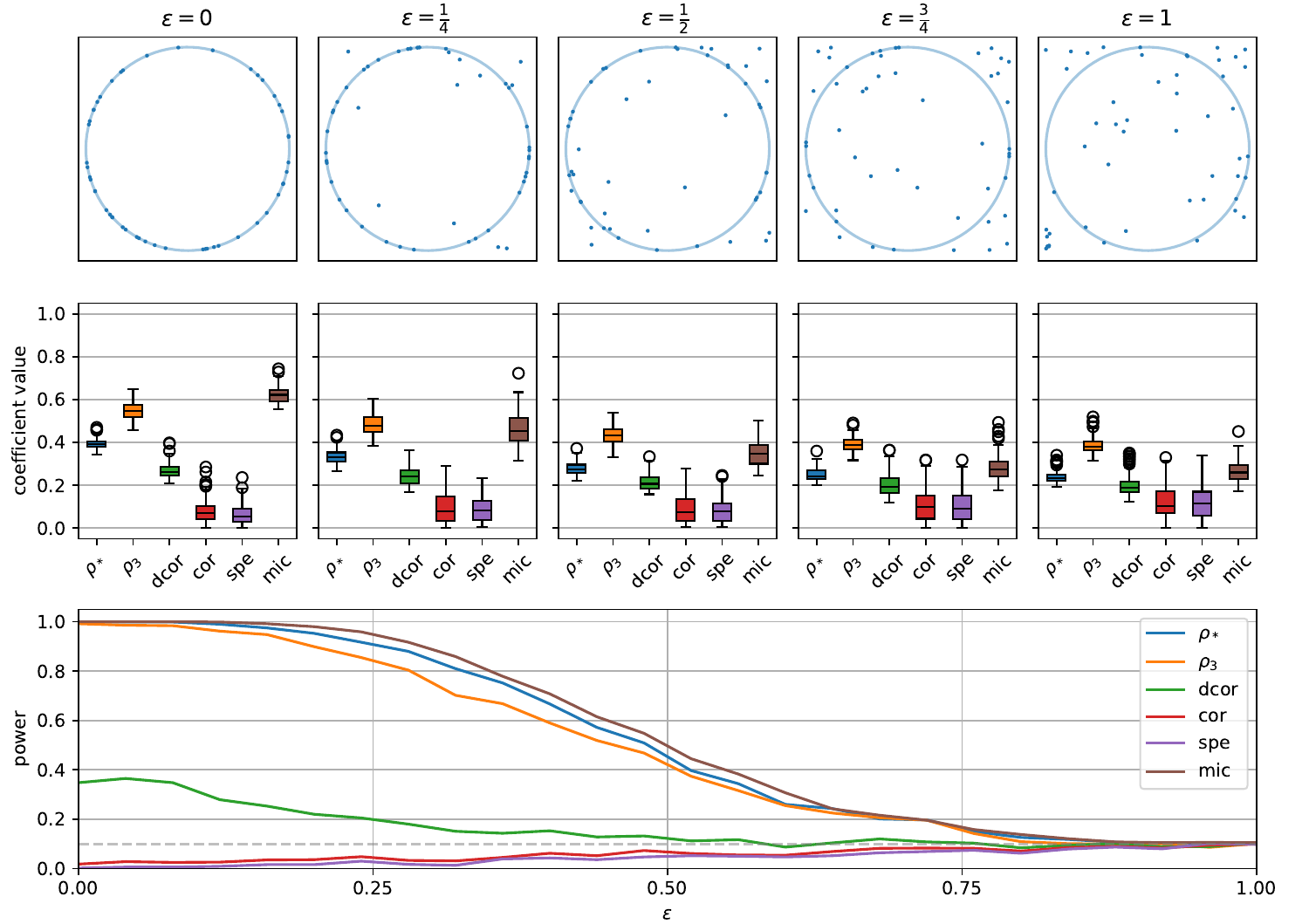}
  \label{fig:appendix06}
  \caption{}
\end{figure}

\begin{figure}
  \centering
  \includegraphics[width=.8\textwidth]{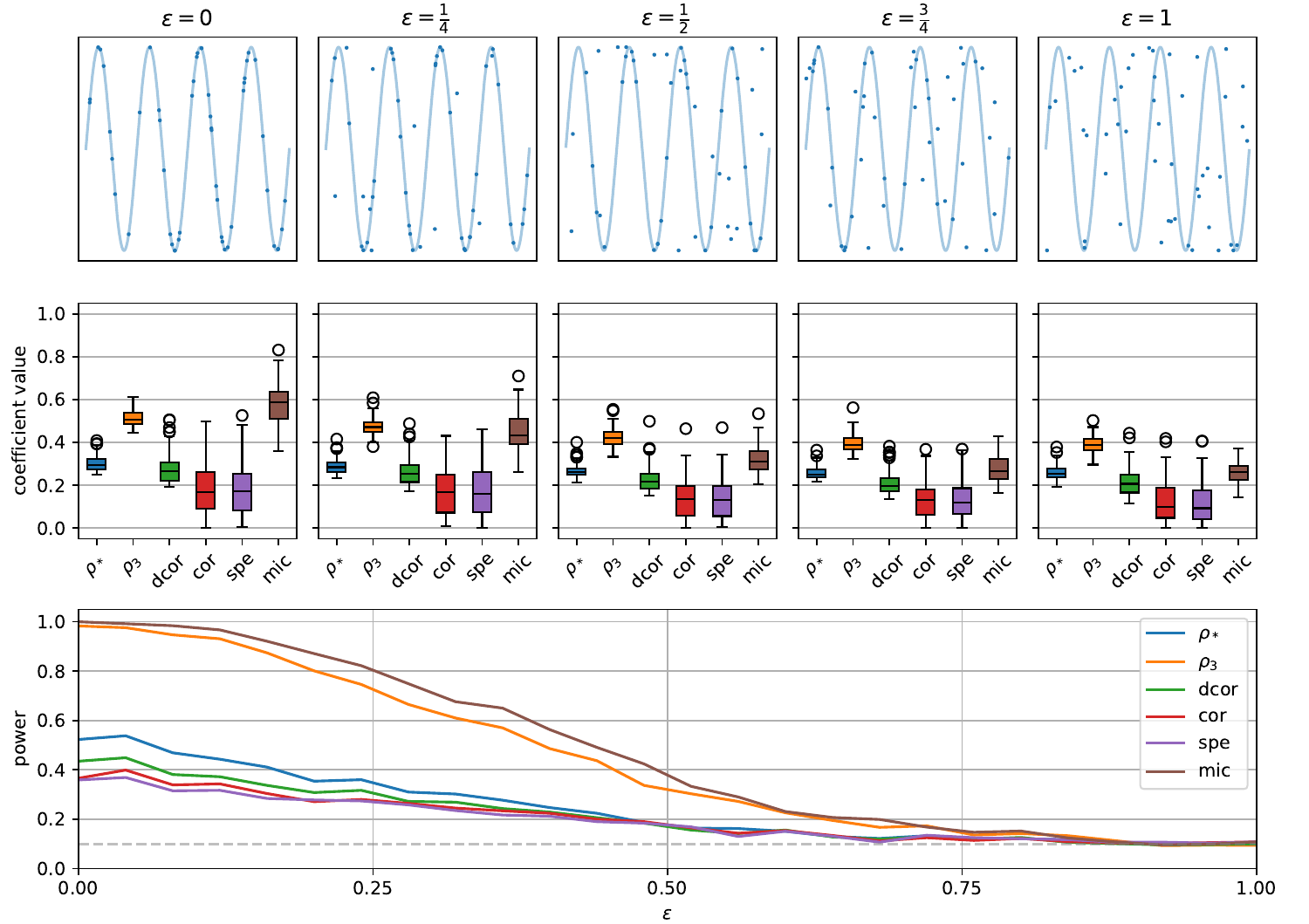}
  \label{fig:appendix07}
  \caption{}
\end{figure}

\begin{figure}
  \centering
  \includegraphics[width=.8\textwidth]{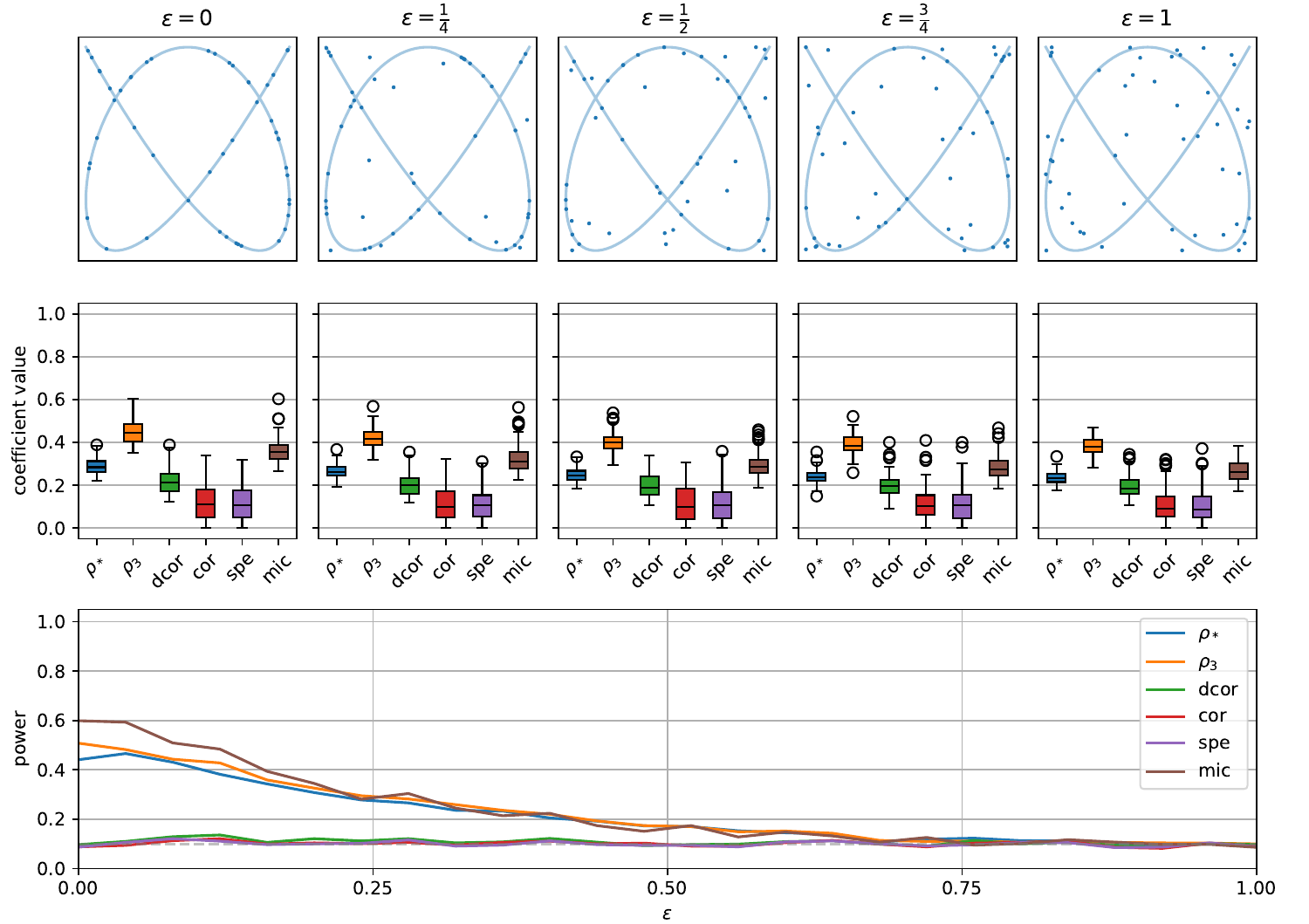}
  \caption{}
  \label{fig:appendix08}
\end{figure}

\paragraph{Convex noise.}
Like Figure~\ref{fig:noise1} and~\ref{fig:noise2} in
Section~\ref{sec:applications}, the upcoming Figures~\ref{fig:appendix01} to
\ref{fig:appendix08} compare $\tc_*$, $\tc_3$, and several other
dependency coefficients under the convex noise
model~\eqref{eq:convex-contamination} for different geometries $\gamma$.
One particularly noteworthy observation is that $\tc_*$ has a higher (or at
least the same) discriminative power compared to the distance correlation,
Pearson correlation, and Spearman correlation for all geometries considered.

\begin{figure}
  \centering
  \includegraphics[width=.8\textwidth]{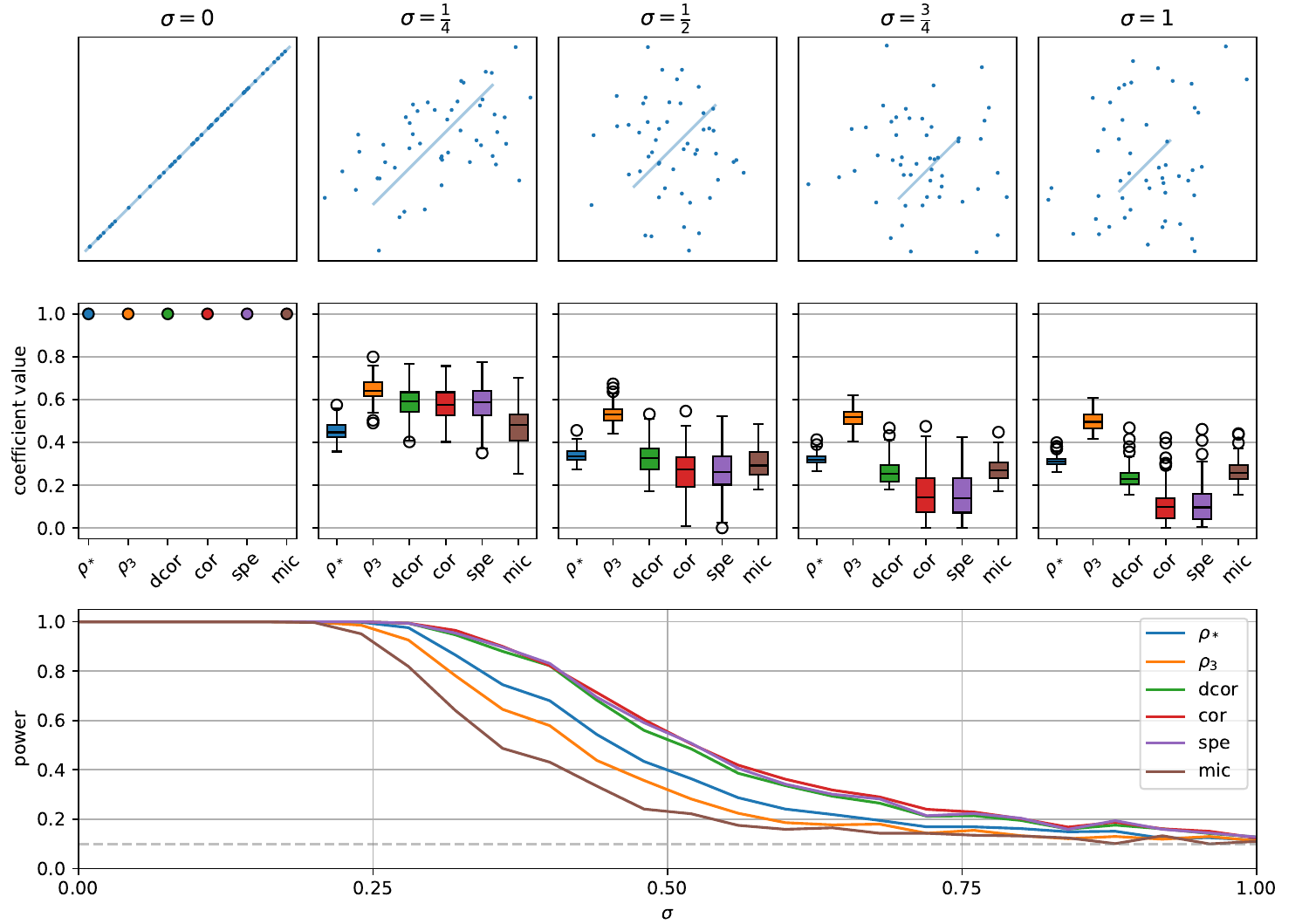}
  \caption{}
  \label{fig:appendix20}
\end{figure}

\begin{figure}
  \centering
  \includegraphics[width=.8\textwidth]{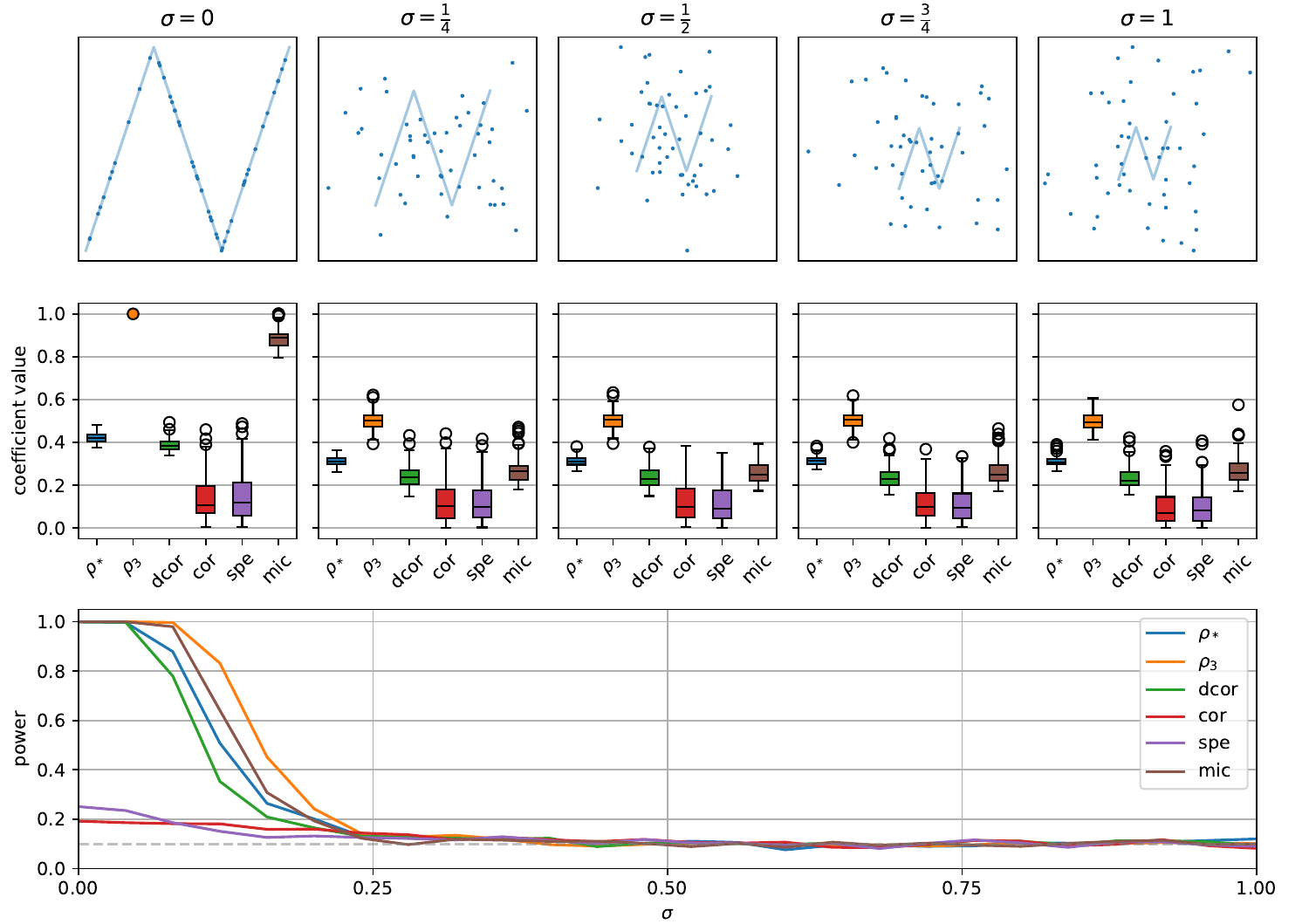}
  \label{fig:appendix21}
  \caption{}
\end{figure}

\begin{figure}
  \centering
  \includegraphics[width=.8\textwidth]{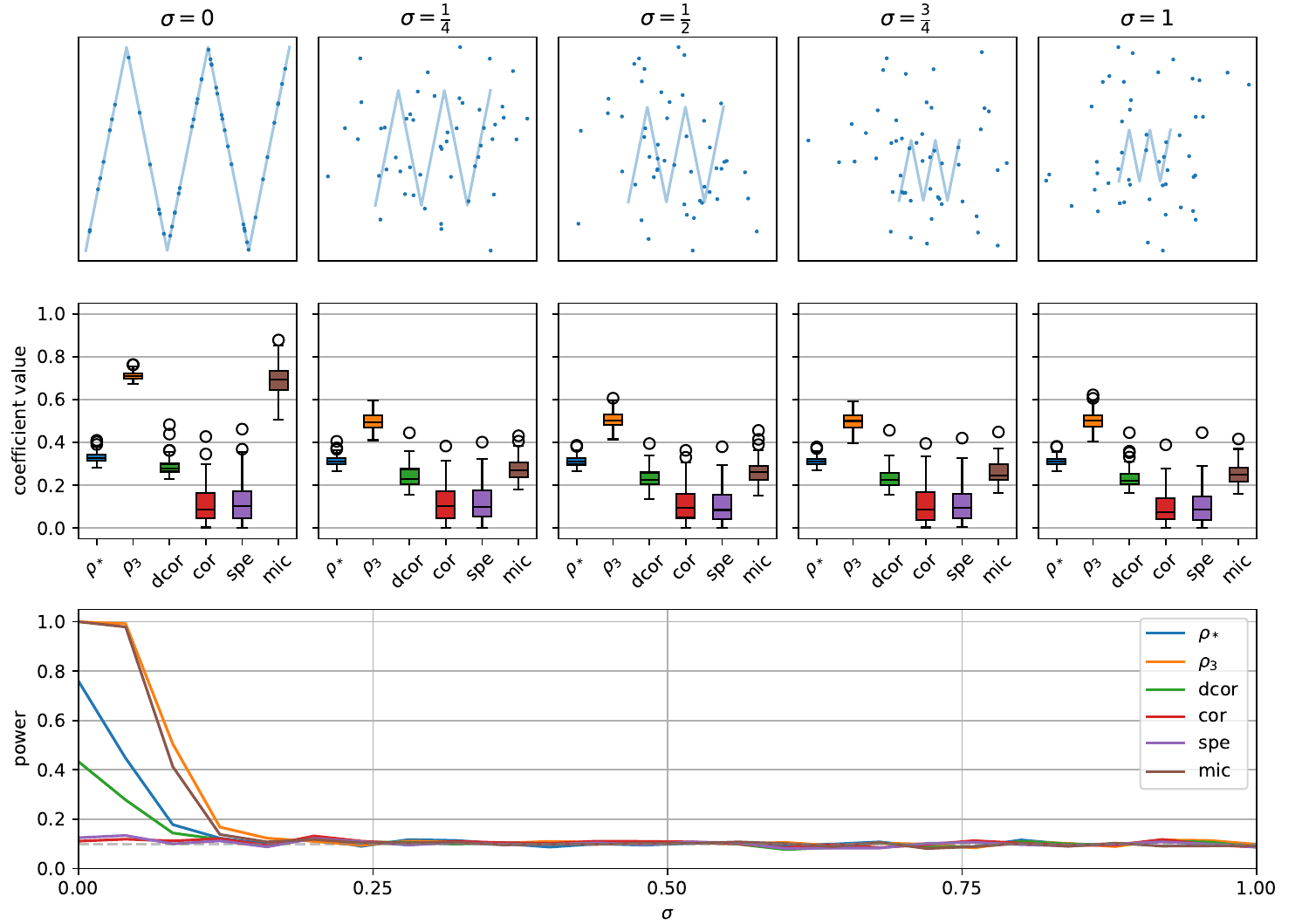}
  \label{fig:appendix22}
  \caption{}
\end{figure}

\begin{figure}
  \centering
  \includegraphics[width=.8\textwidth]{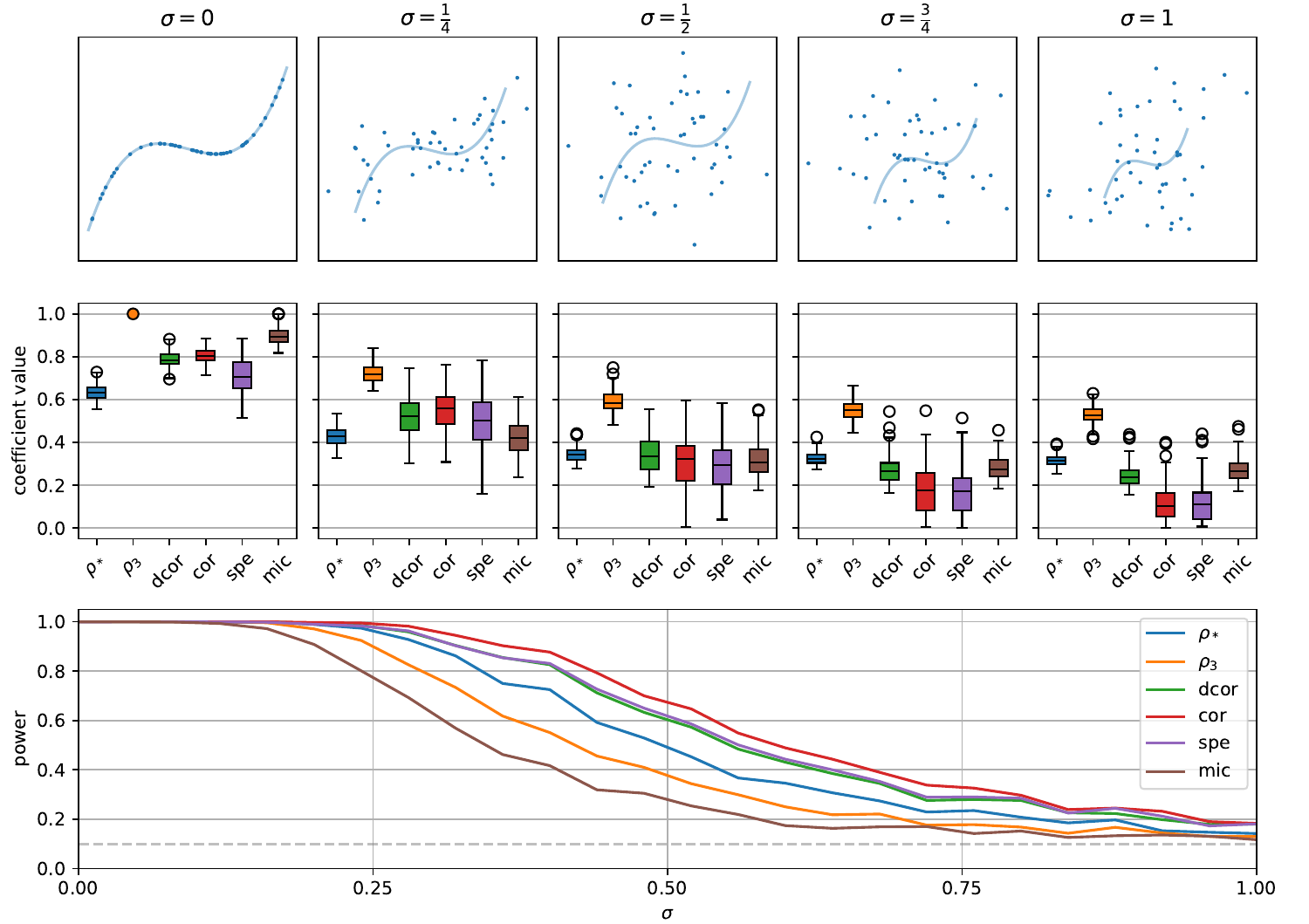}
  \label{fig:appendix23}
  \caption{}
\end{figure}

\begin{figure}
  \centering
  \includegraphics[width=.8\textwidth]{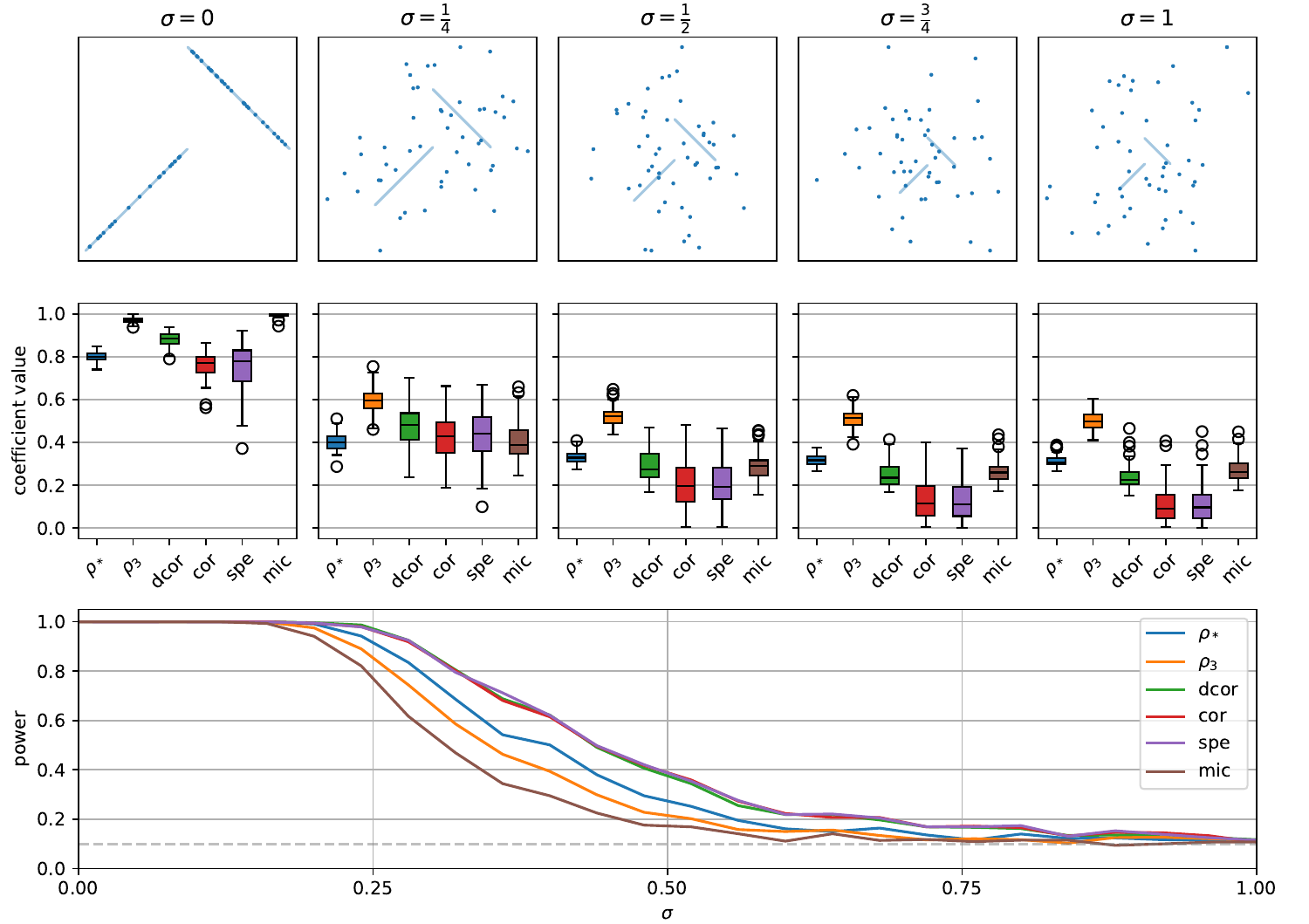}
  \label{fig:appendix24}
  \caption{}
\end{figure}

\begin{figure}
  \centering
  \includegraphics[width=.8\textwidth]{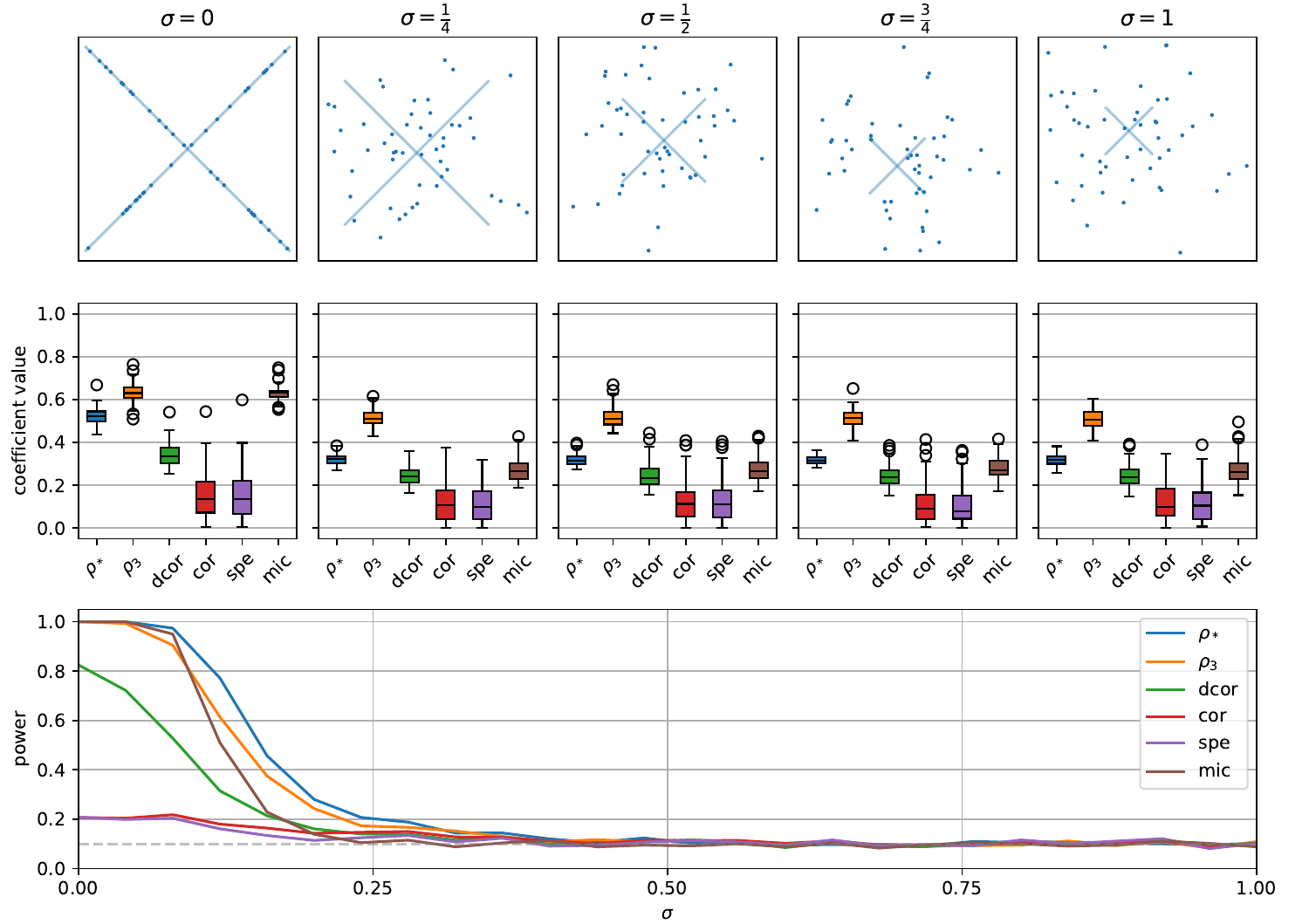}
  \label{fig:appendix25}
  \caption{}
\end{figure}

\begin{figure}
  \centering
  \includegraphics[width=.8\textwidth]{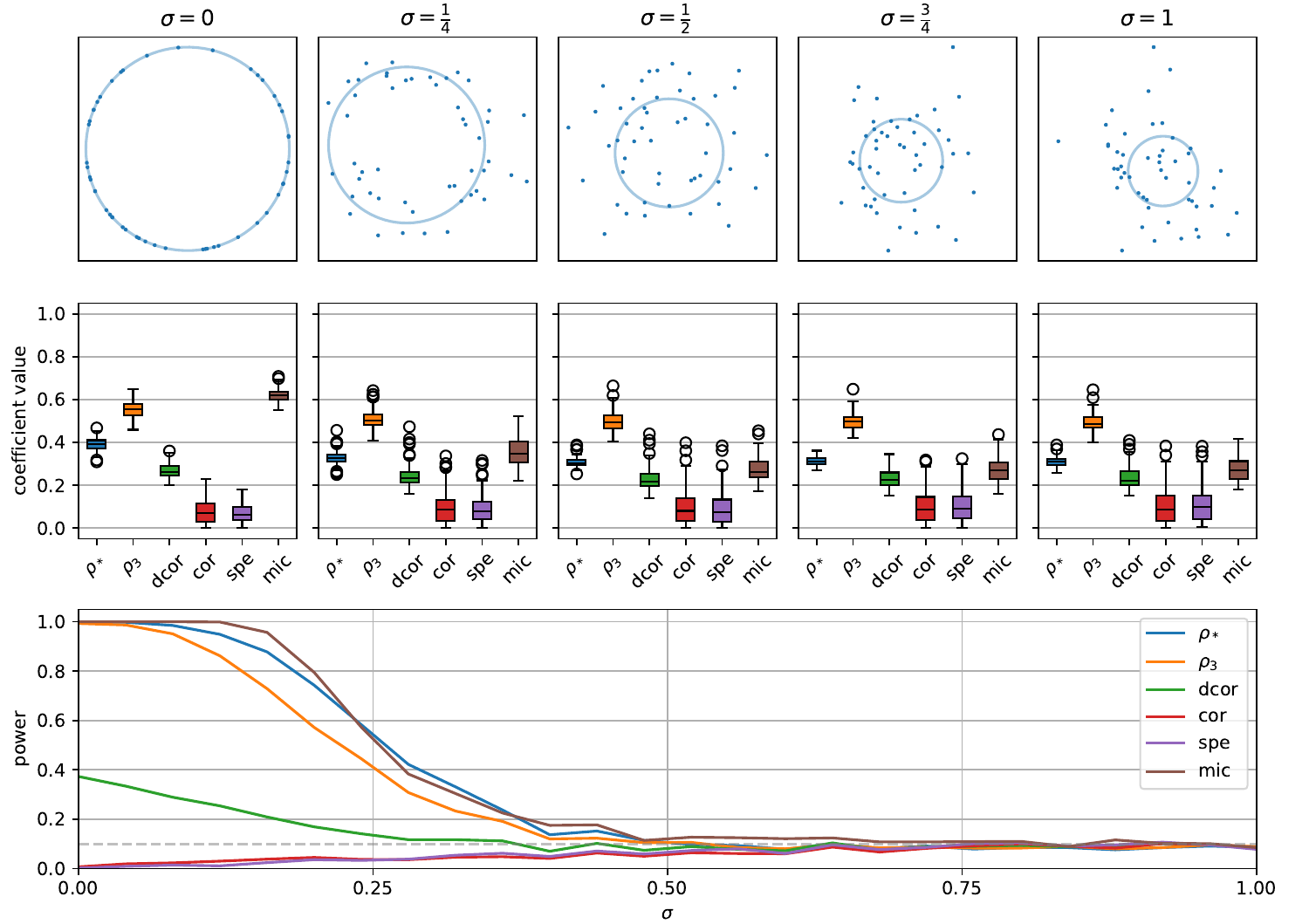}
  \label{fig:appendix26}
  \caption{}
\end{figure}

\begin{figure}
  \centering
  \includegraphics[width=.8\textwidth]{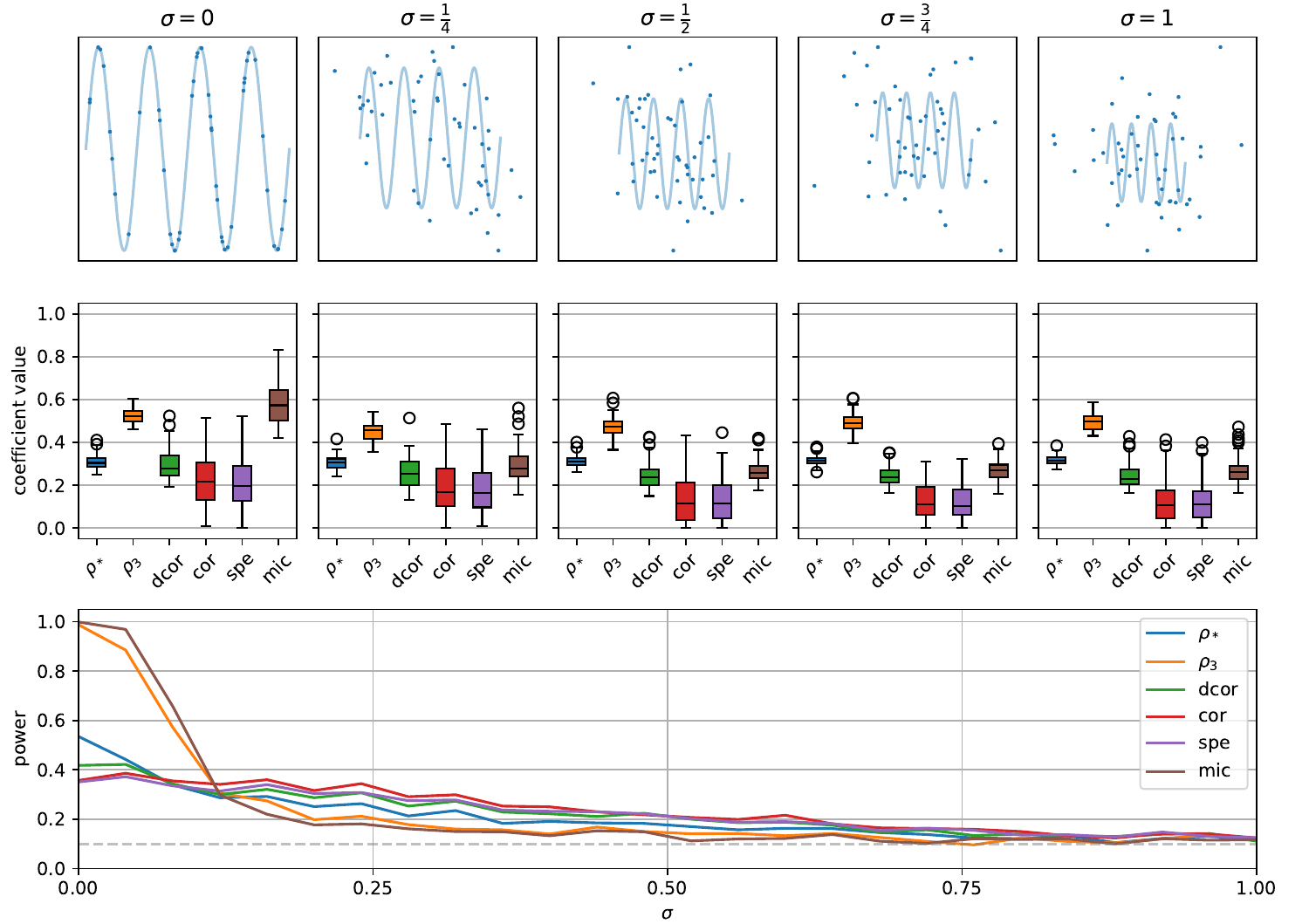}
  \label{fig:appendix27}
  \caption{}
\end{figure}

\begin{figure}
  \centering
  \includegraphics[width=.8\textwidth]{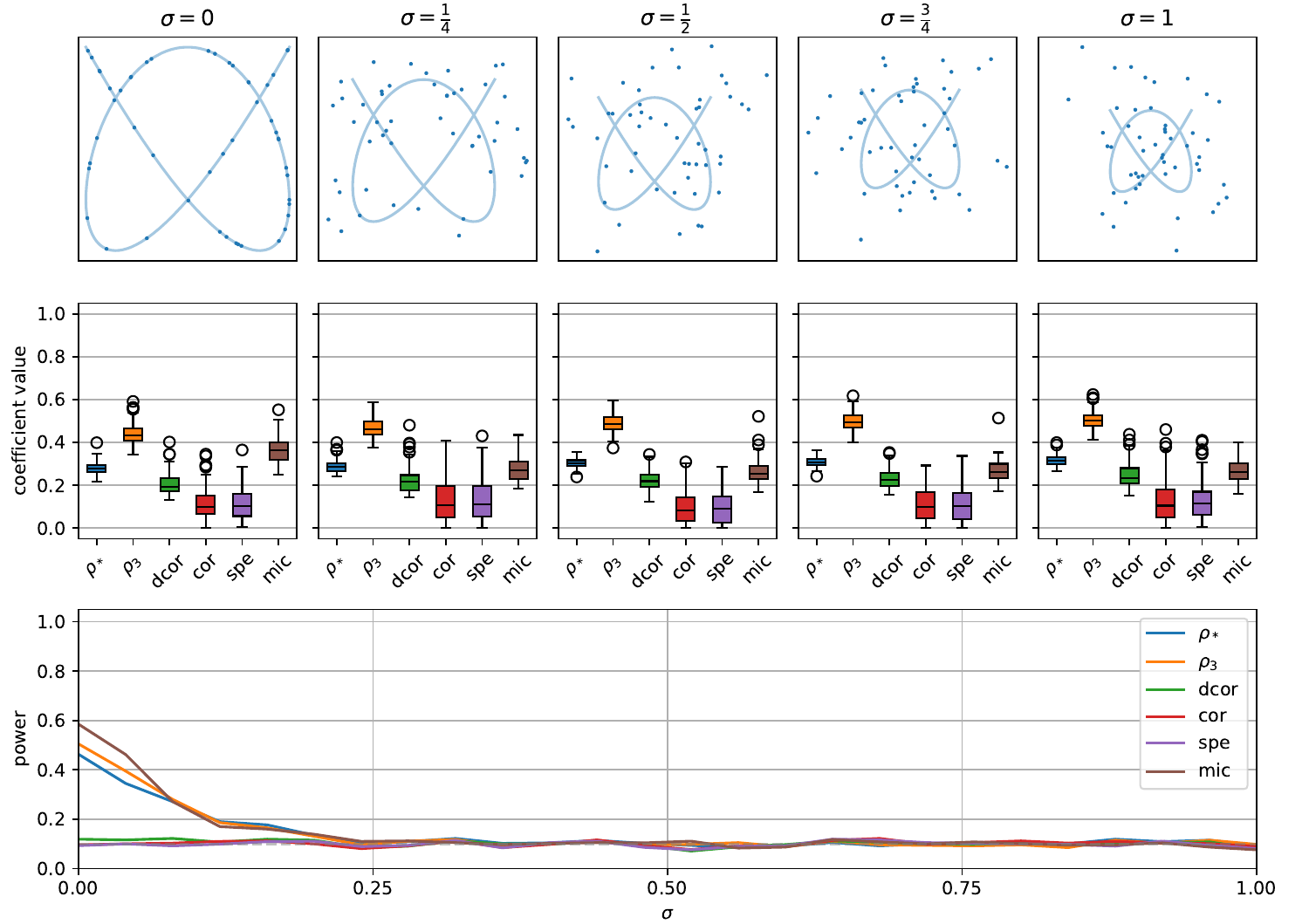}
  \caption{}
  \label{fig:appendix28}
\end{figure}

\paragraph{Gaussian additive noise.}
We repeated our previous simulations under a Gaussian additive noise model
instead of convex noise (Figure~\ref{fig:appendix20} to
Figure~\ref{fig:appendix28}). For a given level of noise $\sigma > 0$ and
a given base distribution $\gamma$, we consider the noisy relationships $(\xi,
\zeta)\sim \gamma^\sigma$, where
\begin{equation}
  \gamma^\sigma= \gamma * \kappa
  \qquad\text{with}\qquad
  \kappa
  =
  \mathcal{N}(0, \sigma\, \mathrm{Id}_2).
\end{equation}
In this setting, we notice that the Pearson and Spearman correlations as well as
the distance correlation have a higher power than $\tc_*$ if $\gamma$ exhibits
a clear monotonic tendency.
If, on the other hand, the linear correlation of the underlying distribution
$\gamma$ is low, then we observe similar trends as for the convex noise model.

\begin{figure}
  \centering
  \includegraphics[width=.8\textwidth]{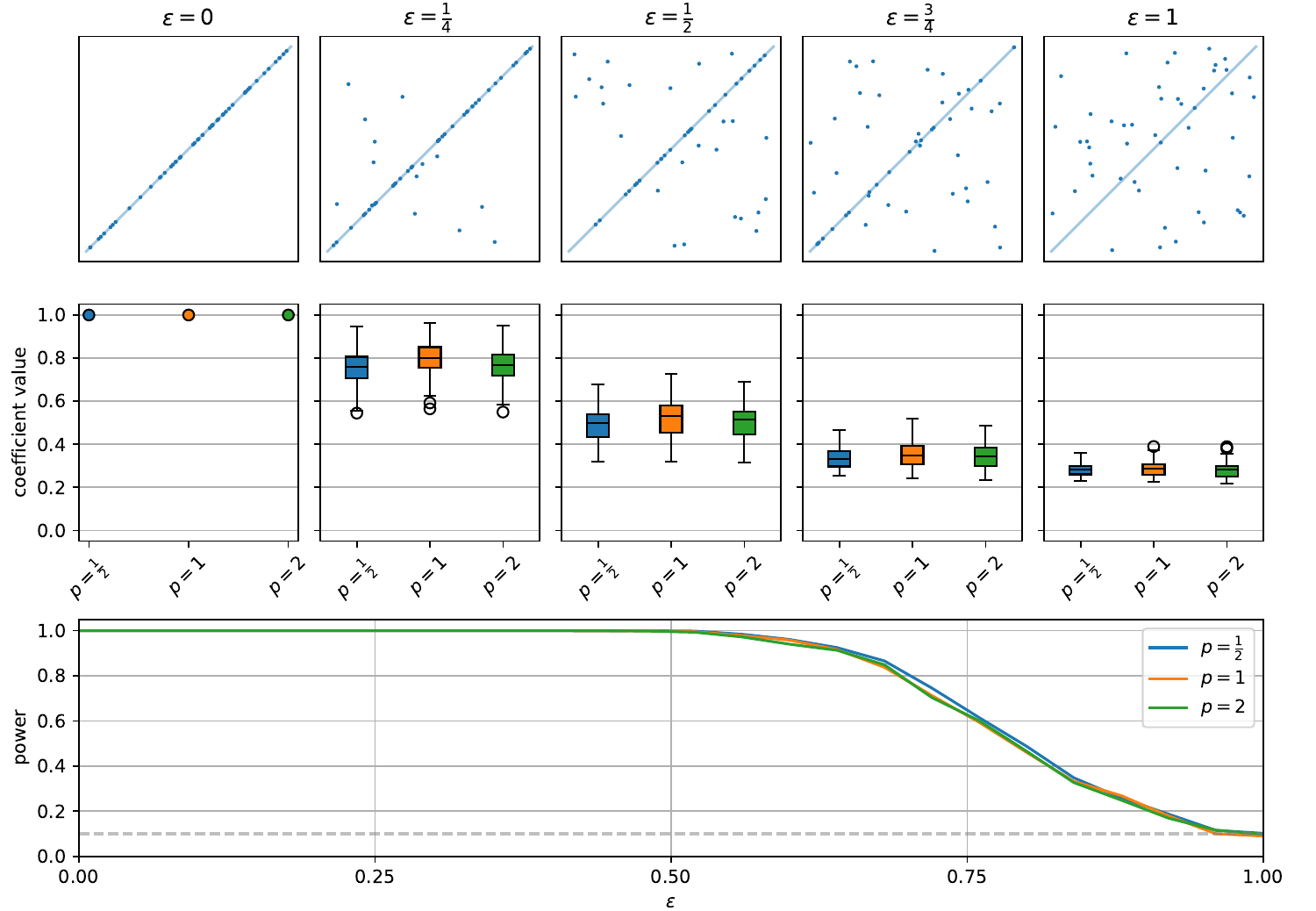}
  \caption{}
  \label{fig:different-p00}
\end{figure}

\begin{figure}
  \centering
  \includegraphics[width=.8\textwidth]{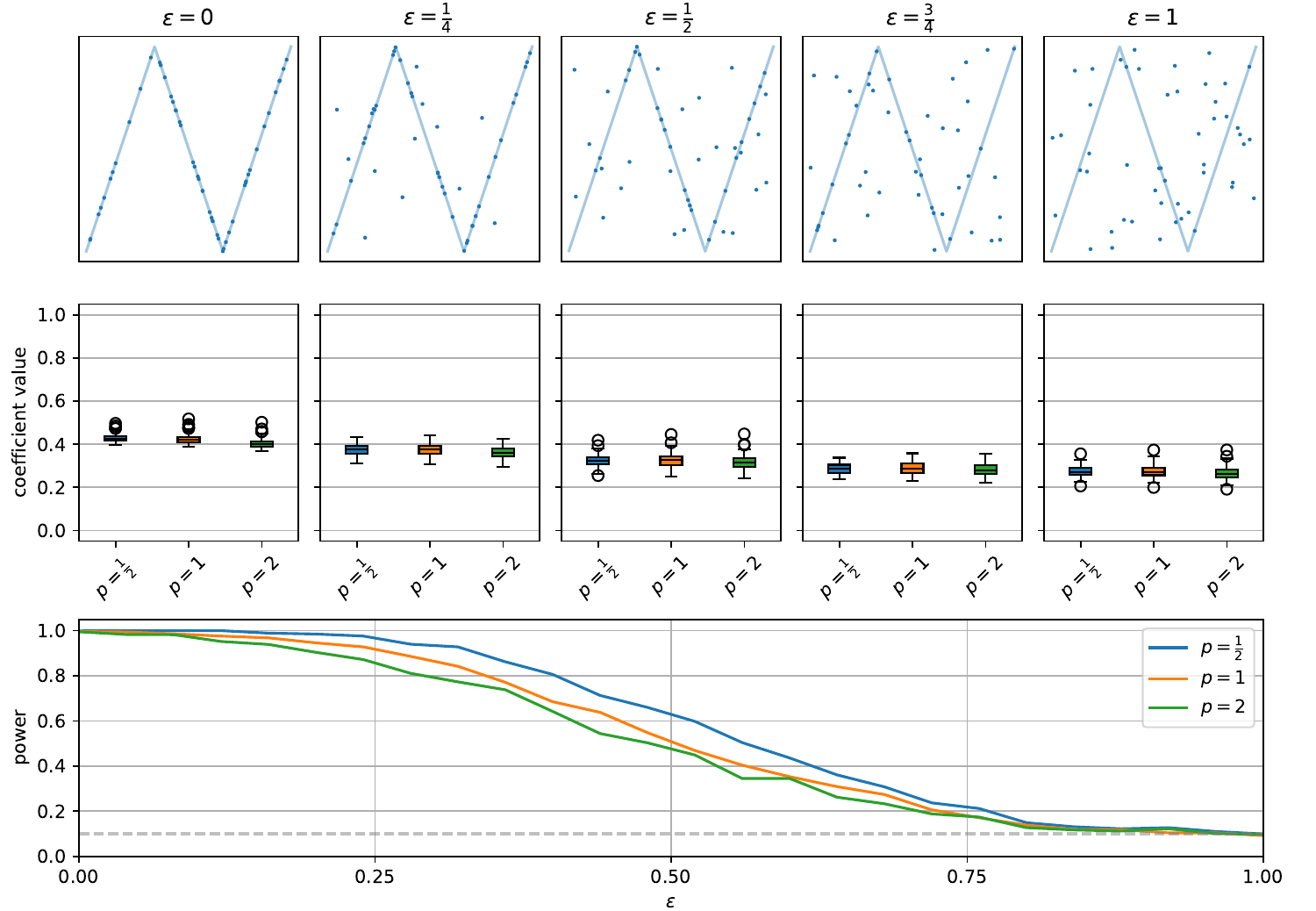}
  \label{fig:different-p01}
  \caption{}
\end{figure}

\begin{figure}
  \centering
  \includegraphics[width=.8\textwidth]{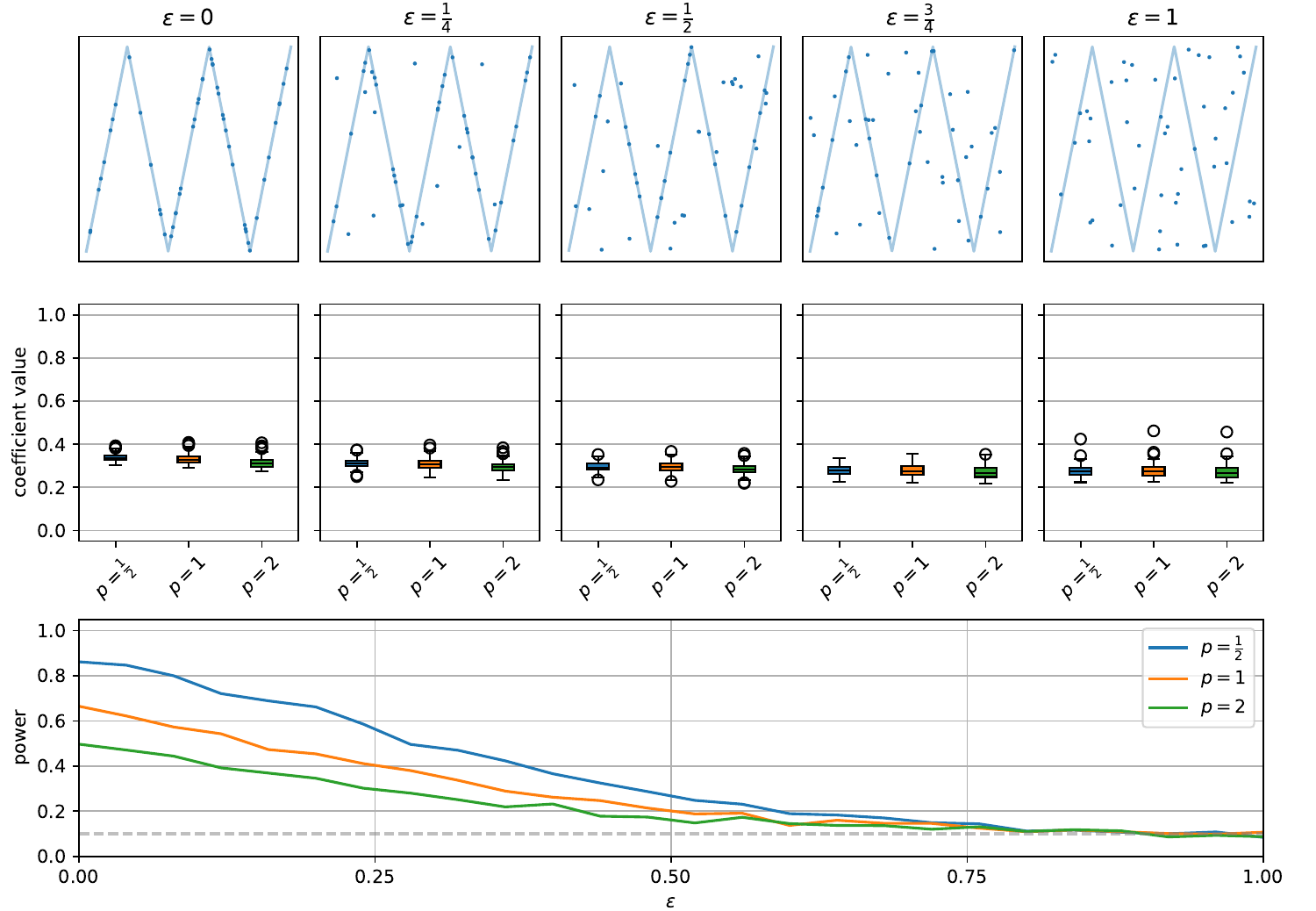}
  \label{fig:different-p02}
  \caption{}
\end{figure}

\begin{figure}
  \centering
  \includegraphics[width=.8\textwidth]{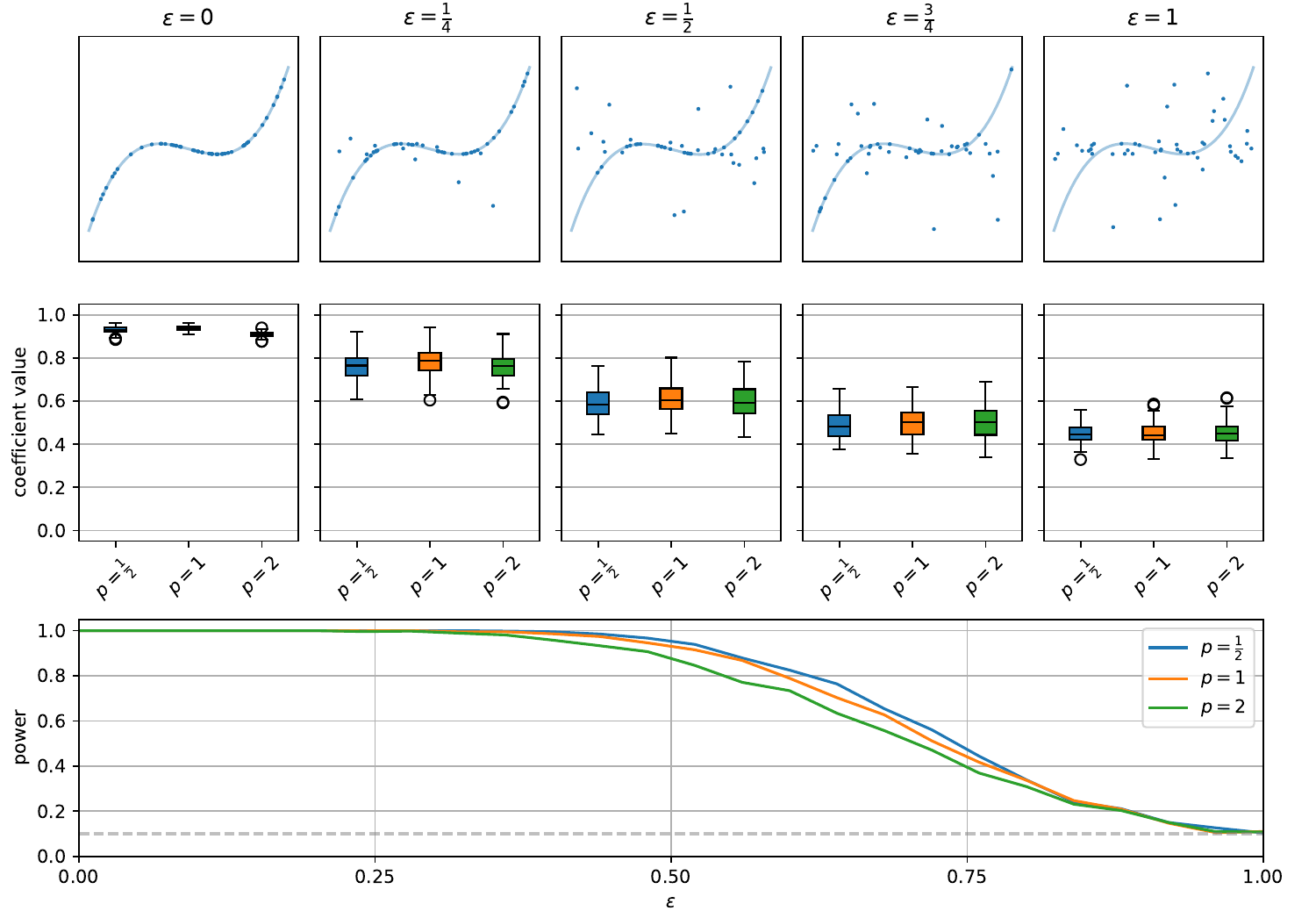}
  \label{fig:different-p03}
  \caption{}
\end{figure}

\begin{figure}
  \centering
  \includegraphics[width=.8\textwidth]{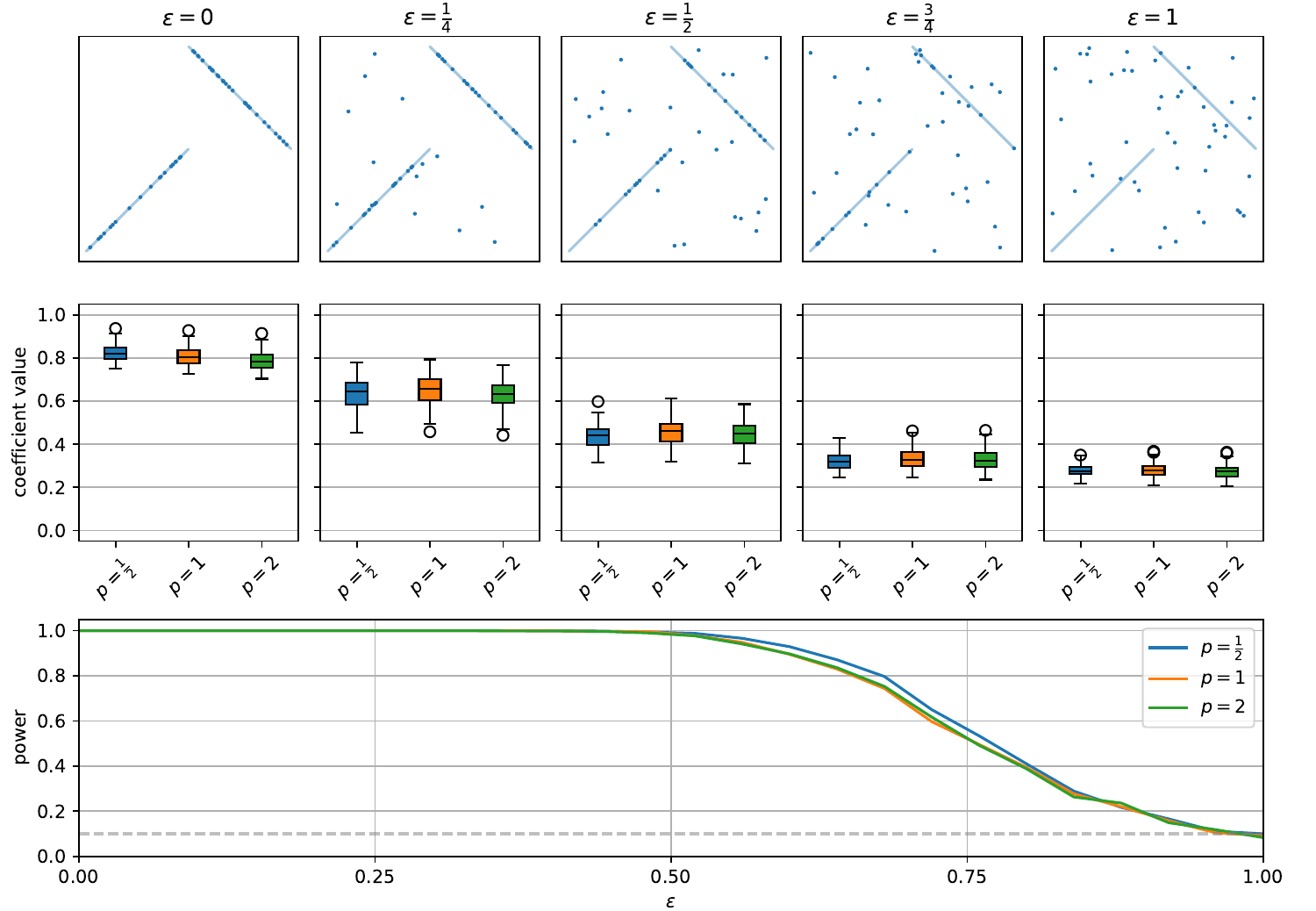}
  \label{fig:different-p04}
  \caption{}
\end{figure}

\begin{figure}
  \centering
  \includegraphics[width=.8\textwidth]{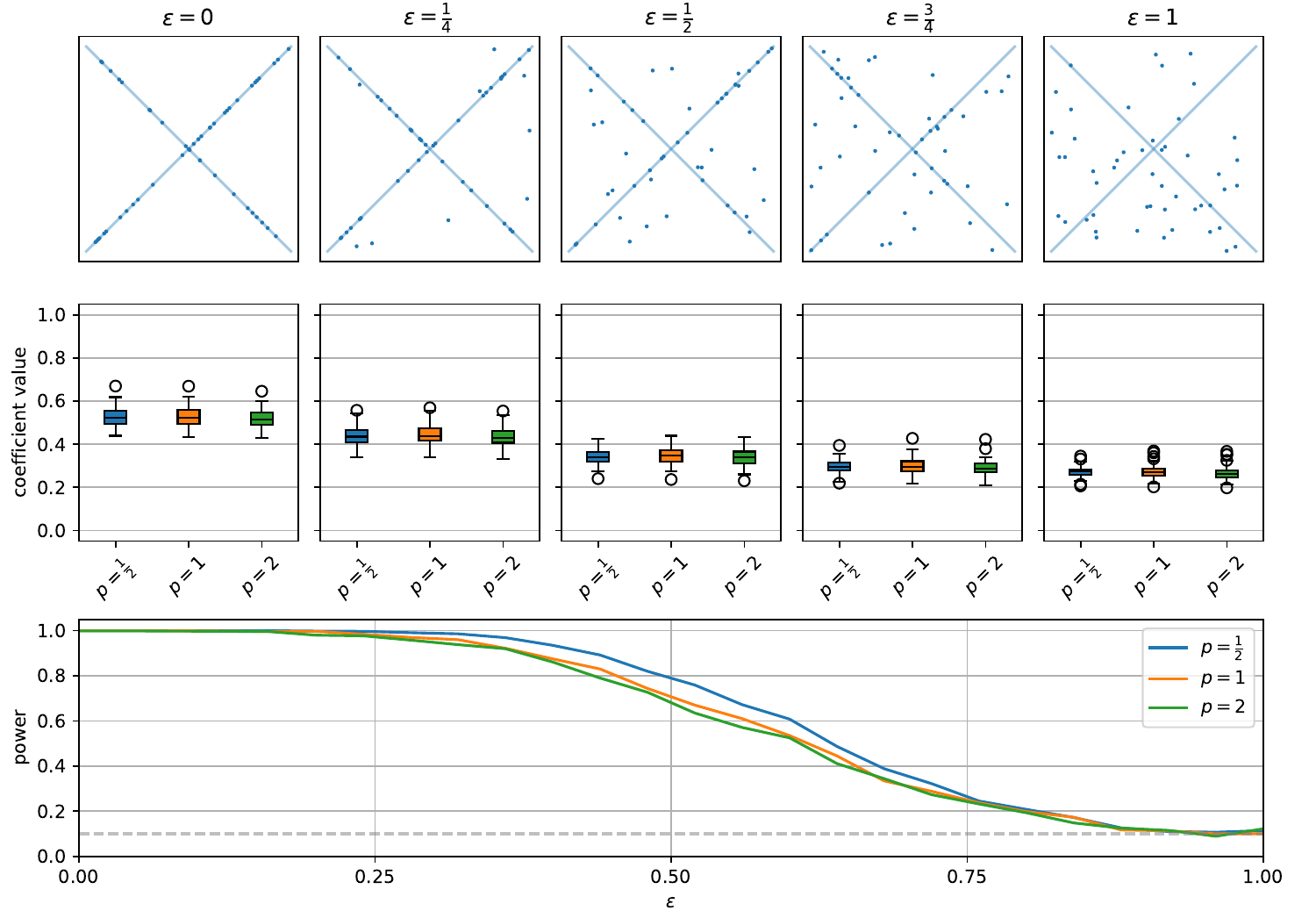}
  \label{fig:different-p05}
  \caption{}
\end{figure}

\begin{figure}
  \centering
  \includegraphics[width=.8\textwidth]{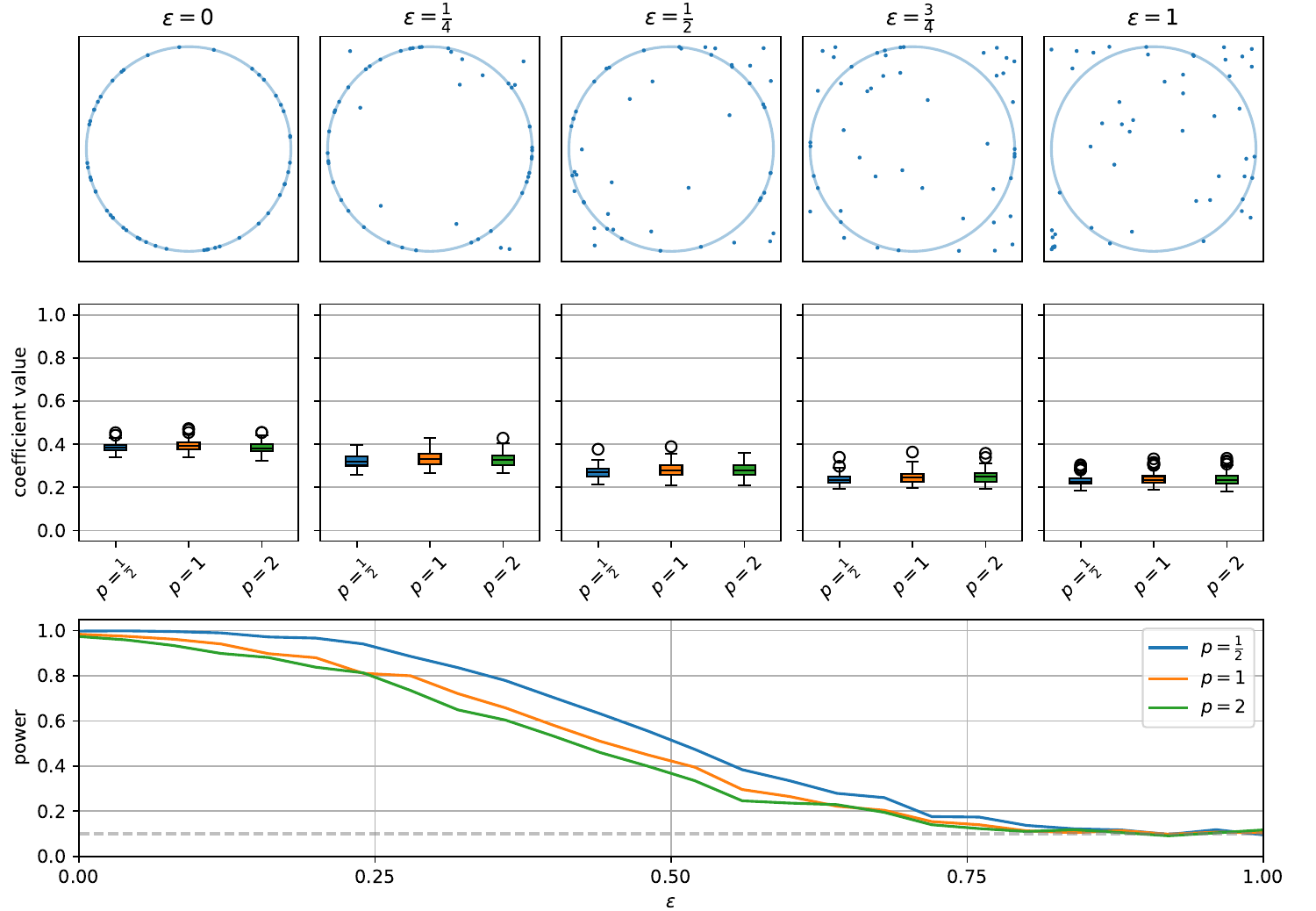}
  \label{fig:different-p06}
  \caption{}
\end{figure}

\begin{figure}
  \centering
  \includegraphics[width=.8\textwidth]{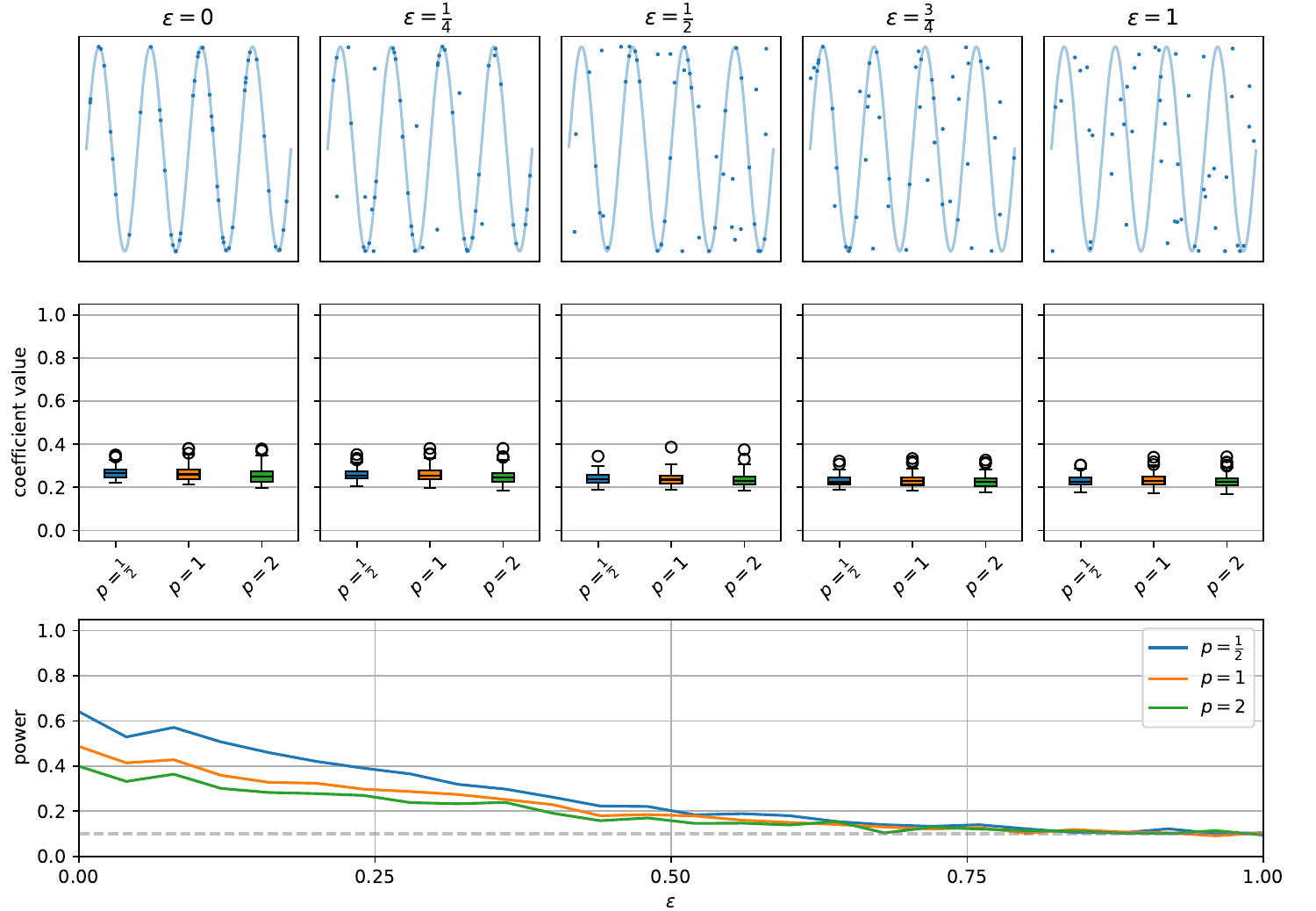}
  \label{fig:different-p07}
  \caption{}
\end{figure}

\begin{figure}
  \centering
  \includegraphics[width=.8\textwidth]{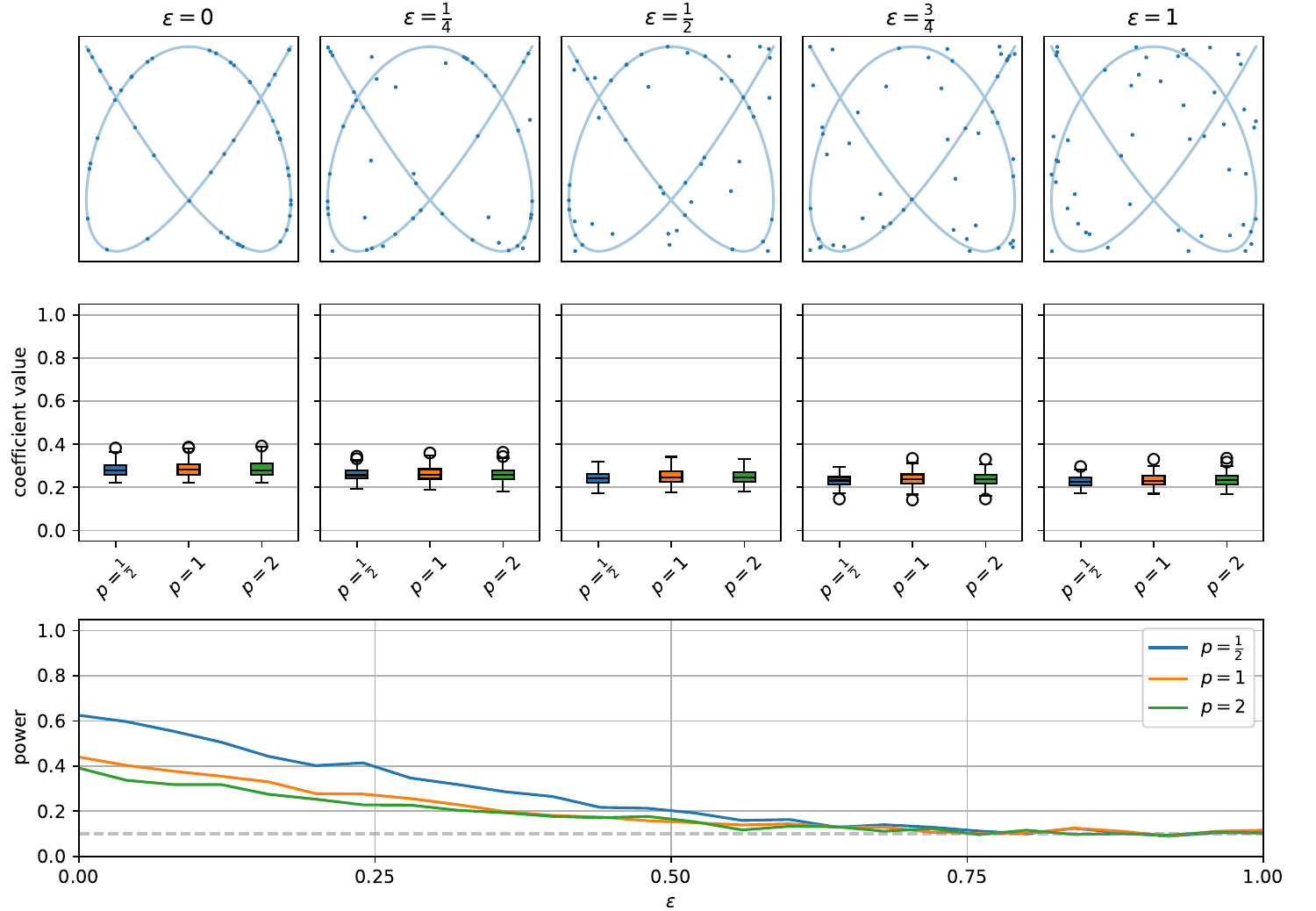}
  \caption{}
  \label{fig:different-p08}
\end{figure}

\begin{figure}
  \centering
  \includegraphics[width=.8\textwidth]{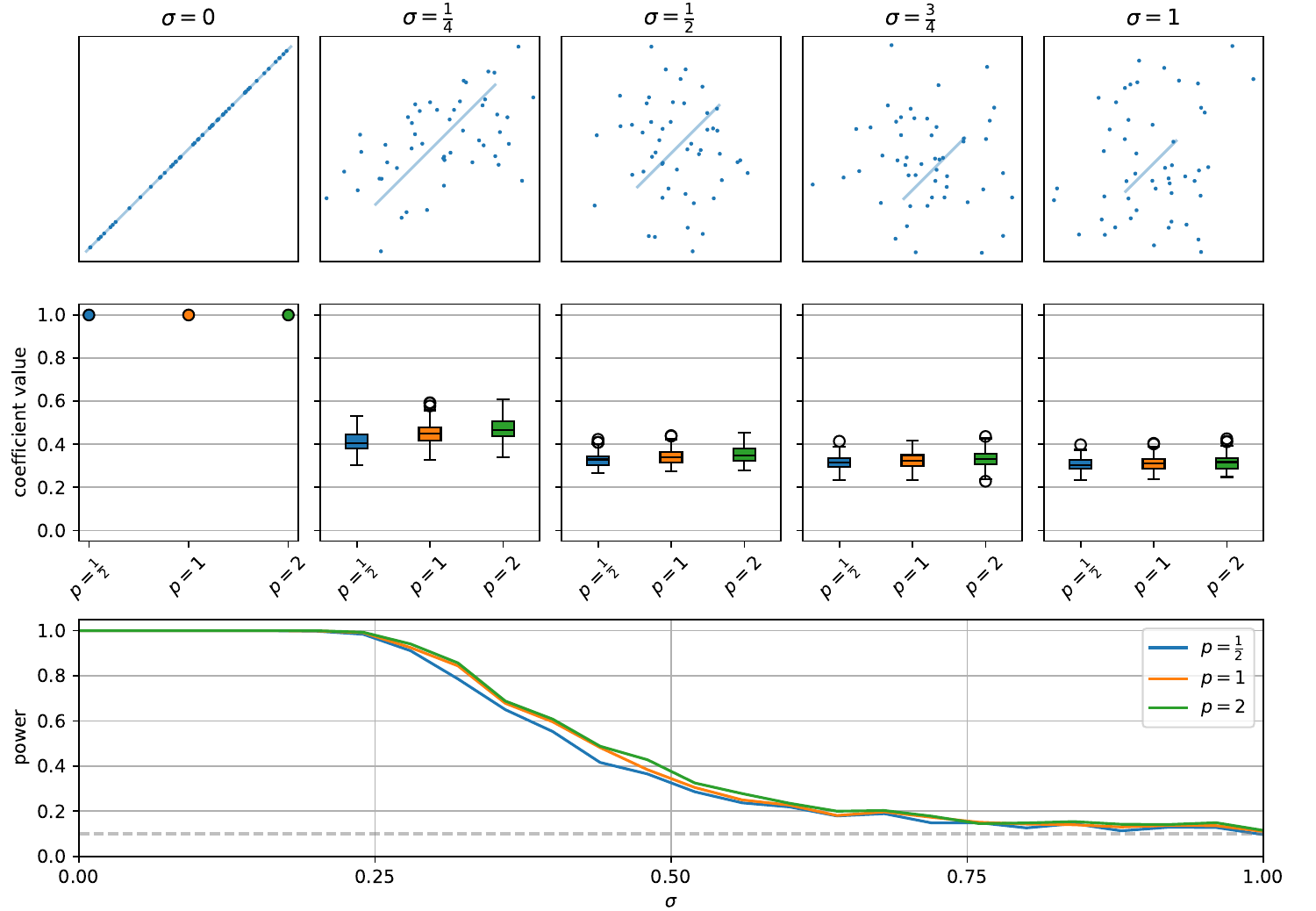}
  \caption{}
  \label{fig:different-p20}
\end{figure}

\begin{figure}
  \centering
  \includegraphics[width=.8\textwidth]{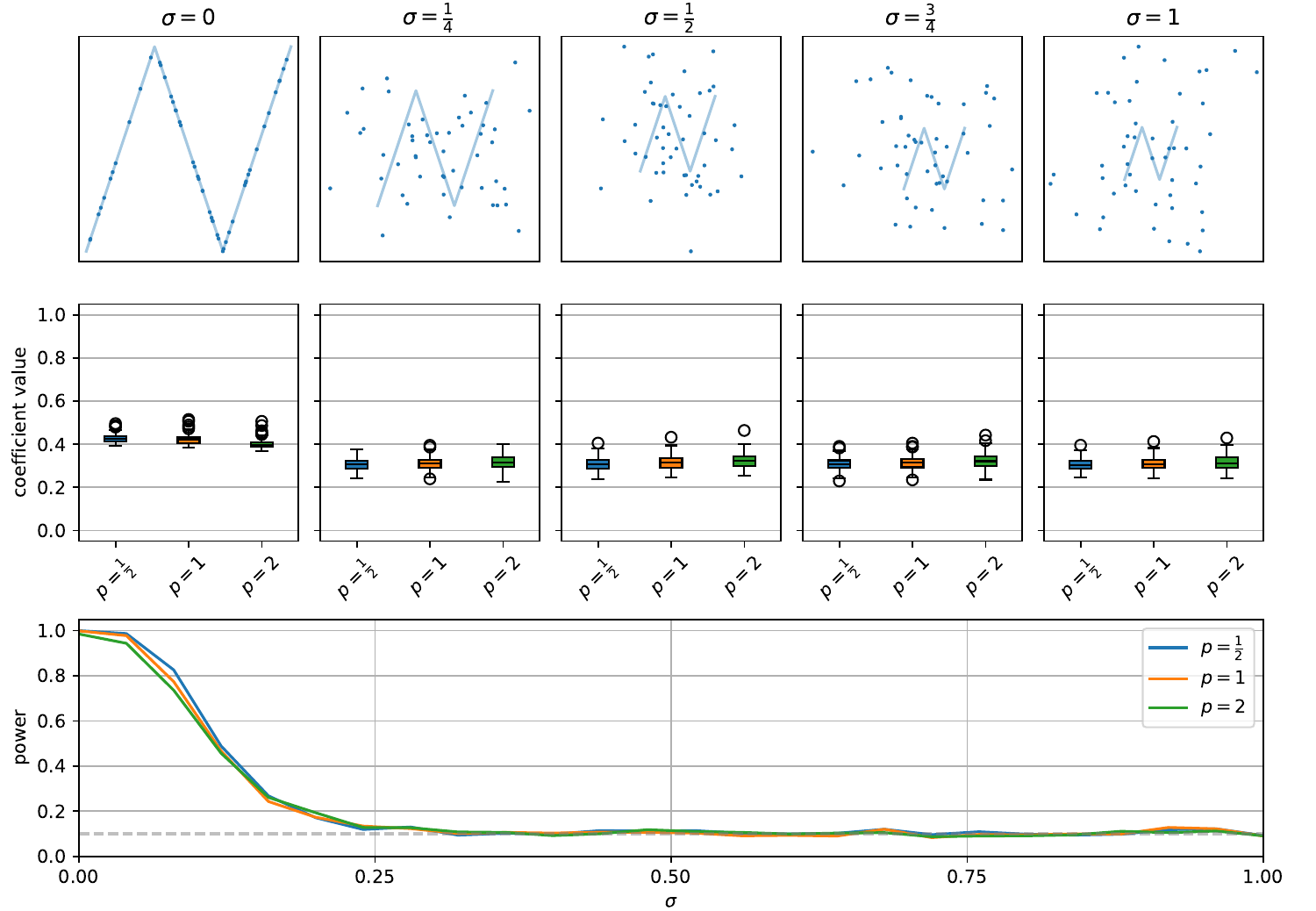}
  \label{fig:different-p21}
  \caption{}
\end{figure}

\begin{figure}
  \centering
  \includegraphics[width=.8\textwidth]{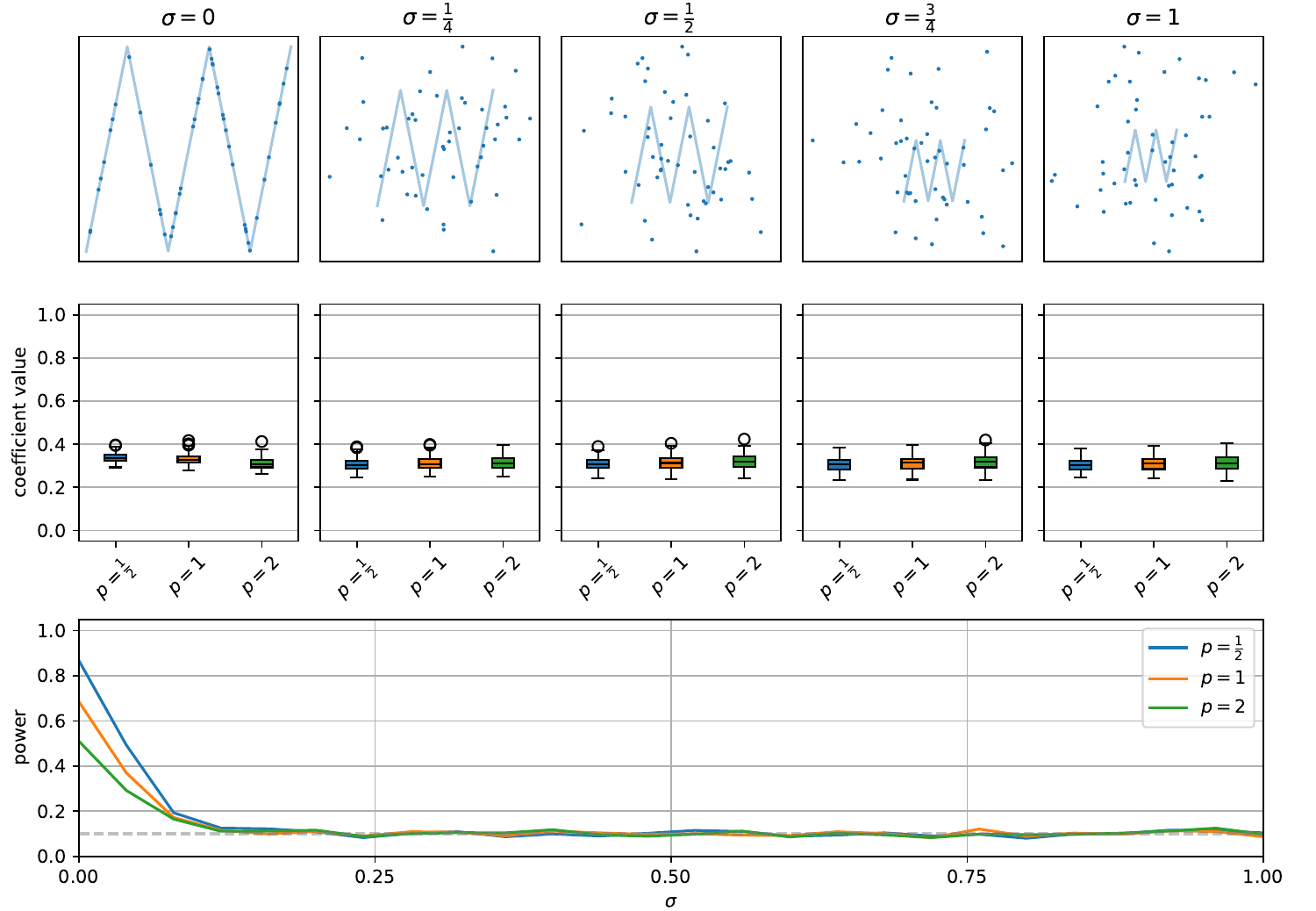}
  \label{fig:different-p22}
  \caption{}
\end{figure}

\begin{figure}
  \centering
  \includegraphics[width=.8\textwidth]{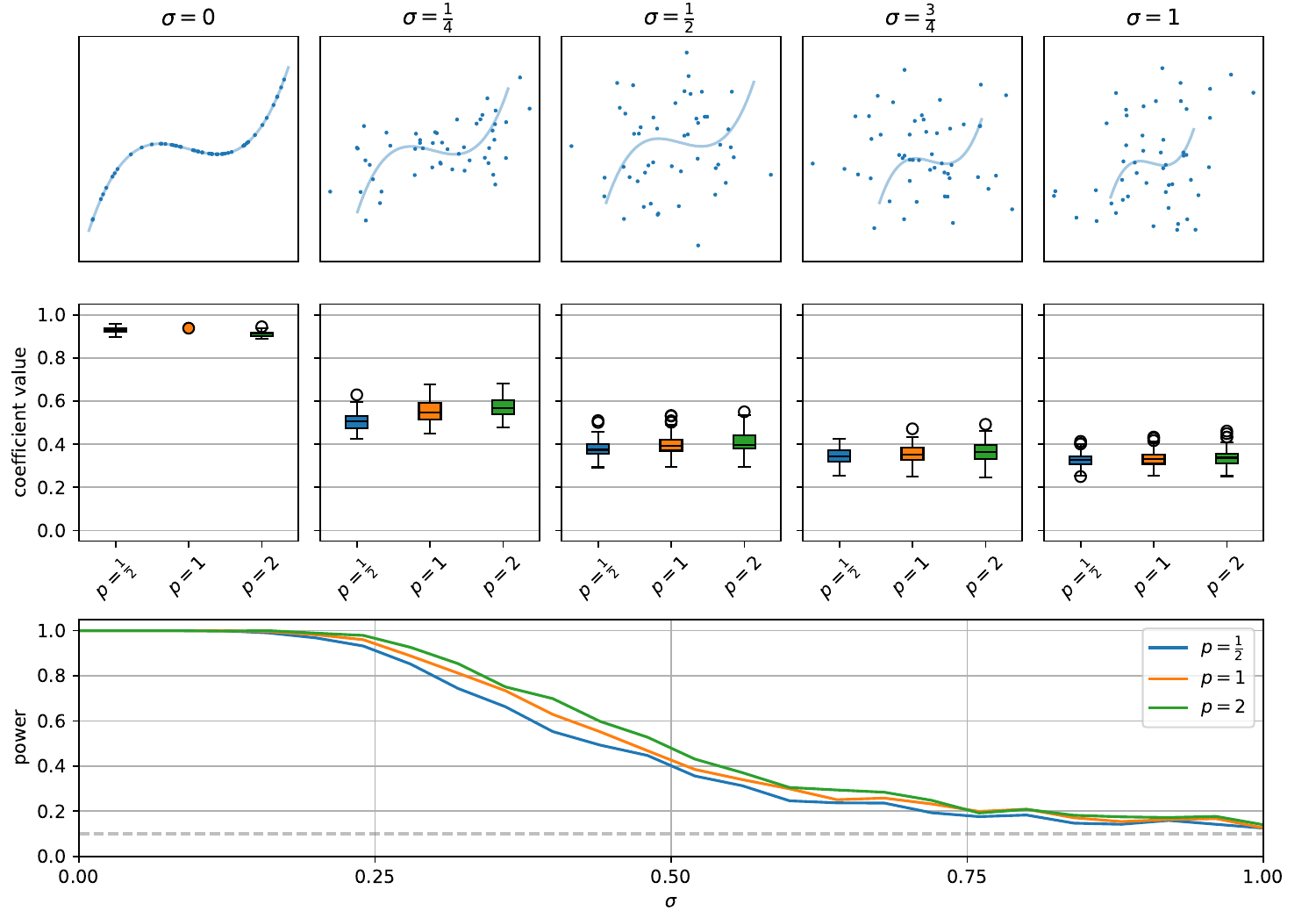}
  \label{fig:different-p23}
  \caption{}
\end{figure}

\begin{figure}
  \centering
  \includegraphics[width=.8\textwidth]{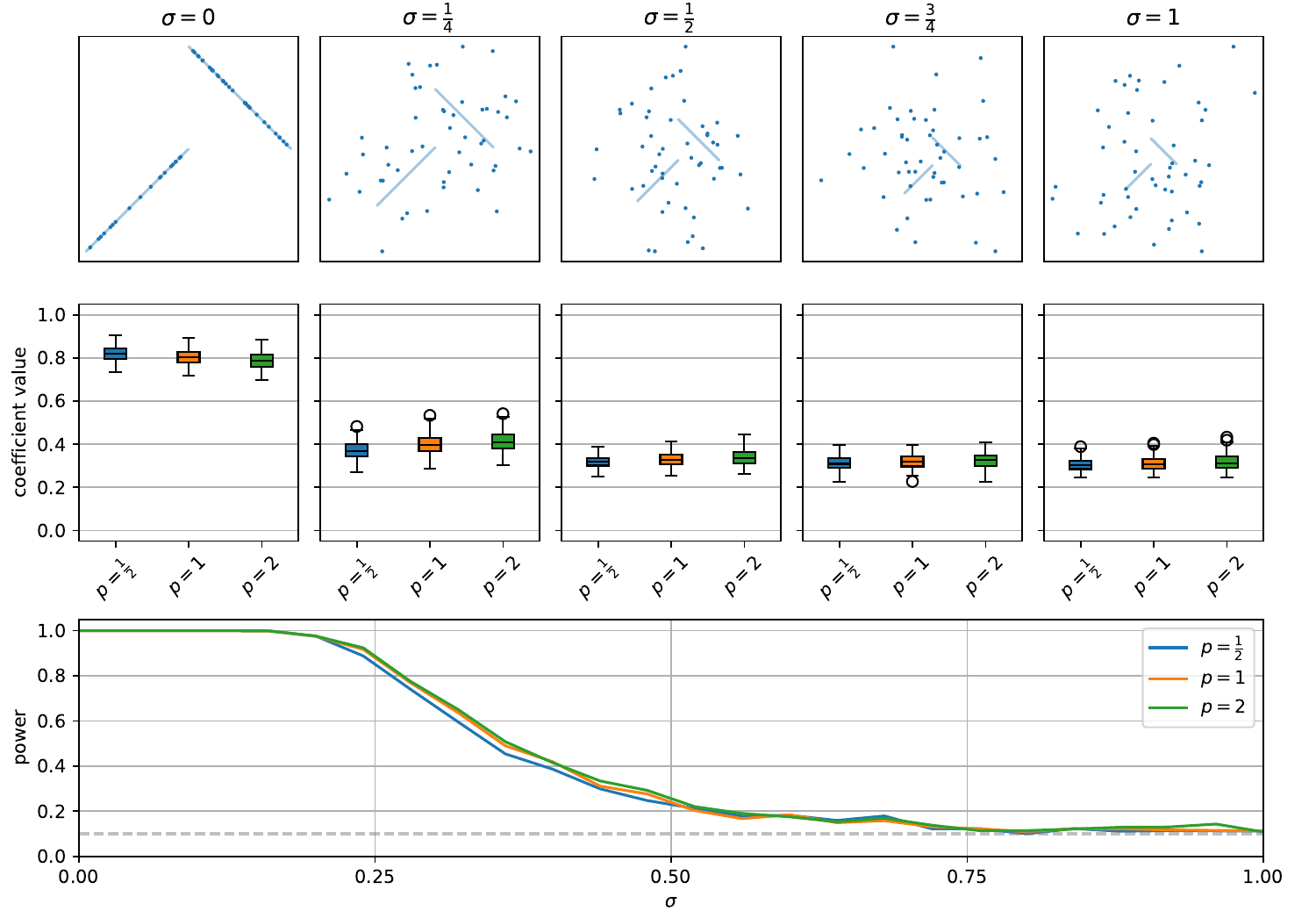}
  \label{fig:different-p24}
  \caption{}
\end{figure}

\begin{figure}
  \centering
  \includegraphics[width=.8\textwidth]{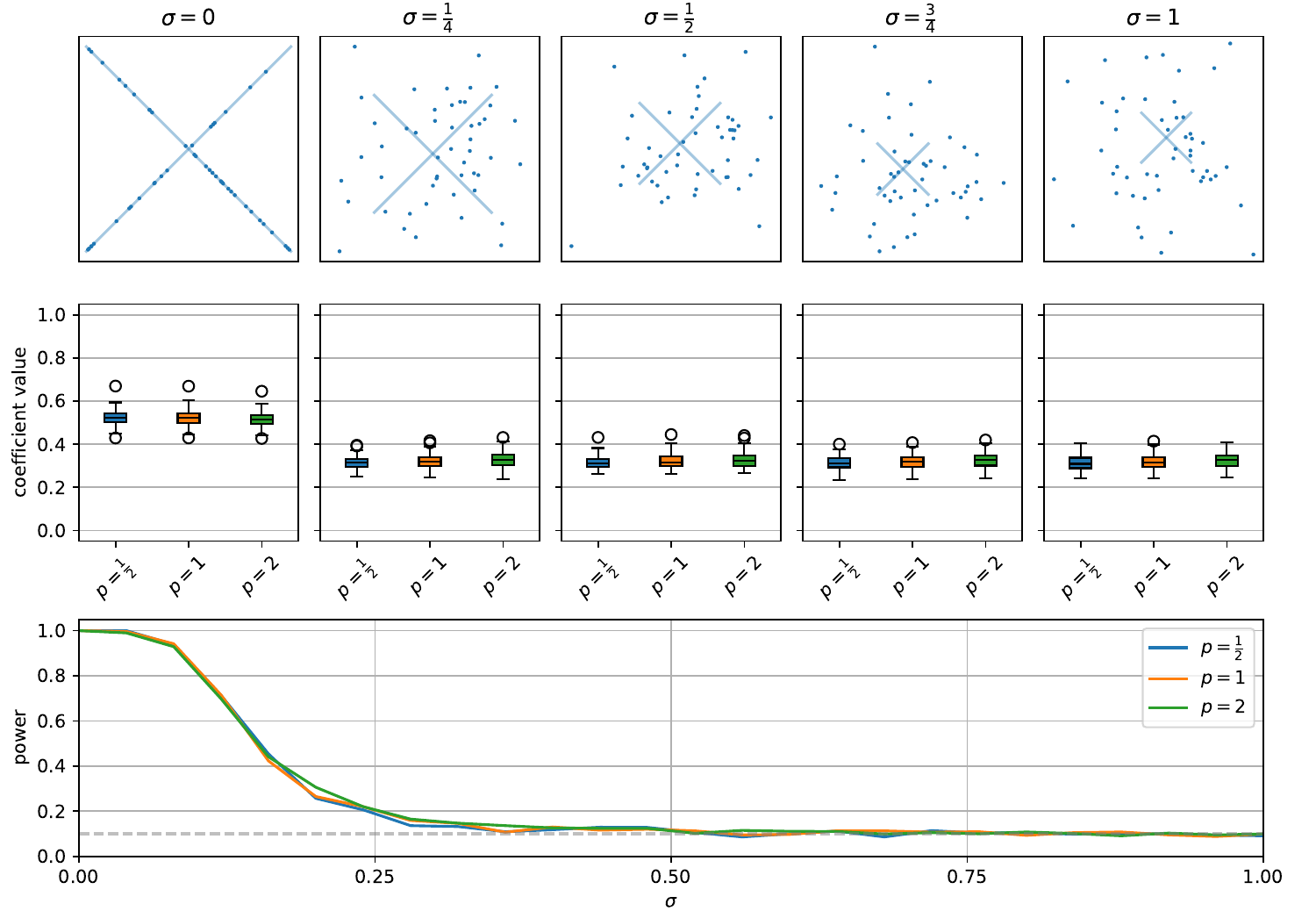}
  \label{fig:different-p25}
  \caption{}
\end{figure}

\begin{figure}
  \centering
  \includegraphics[width=.8\textwidth]{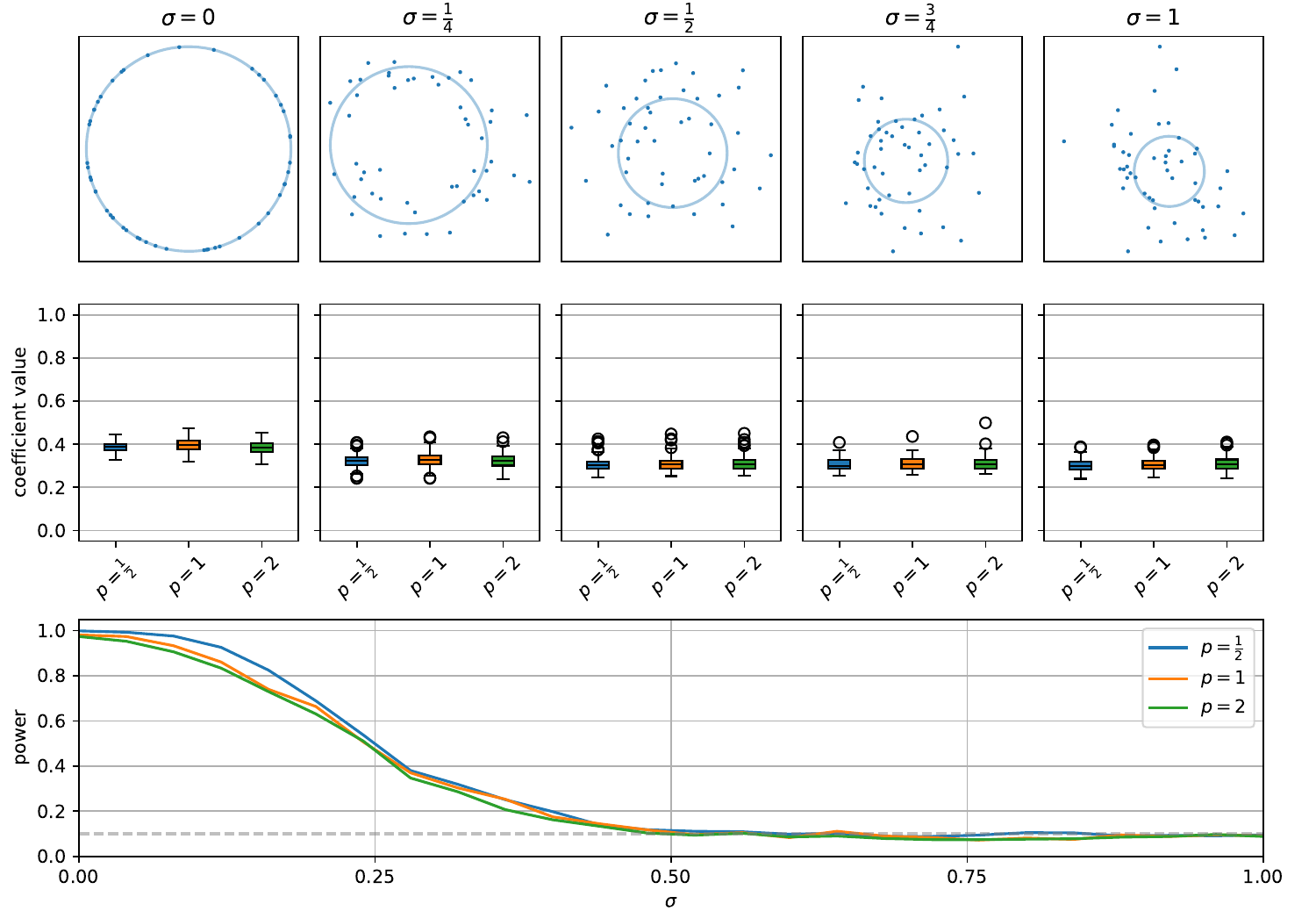}
  \label{fig:different-p26}
  \caption{}
\end{figure}

\begin{figure}
  \centering
  \includegraphics[width=.8\textwidth]{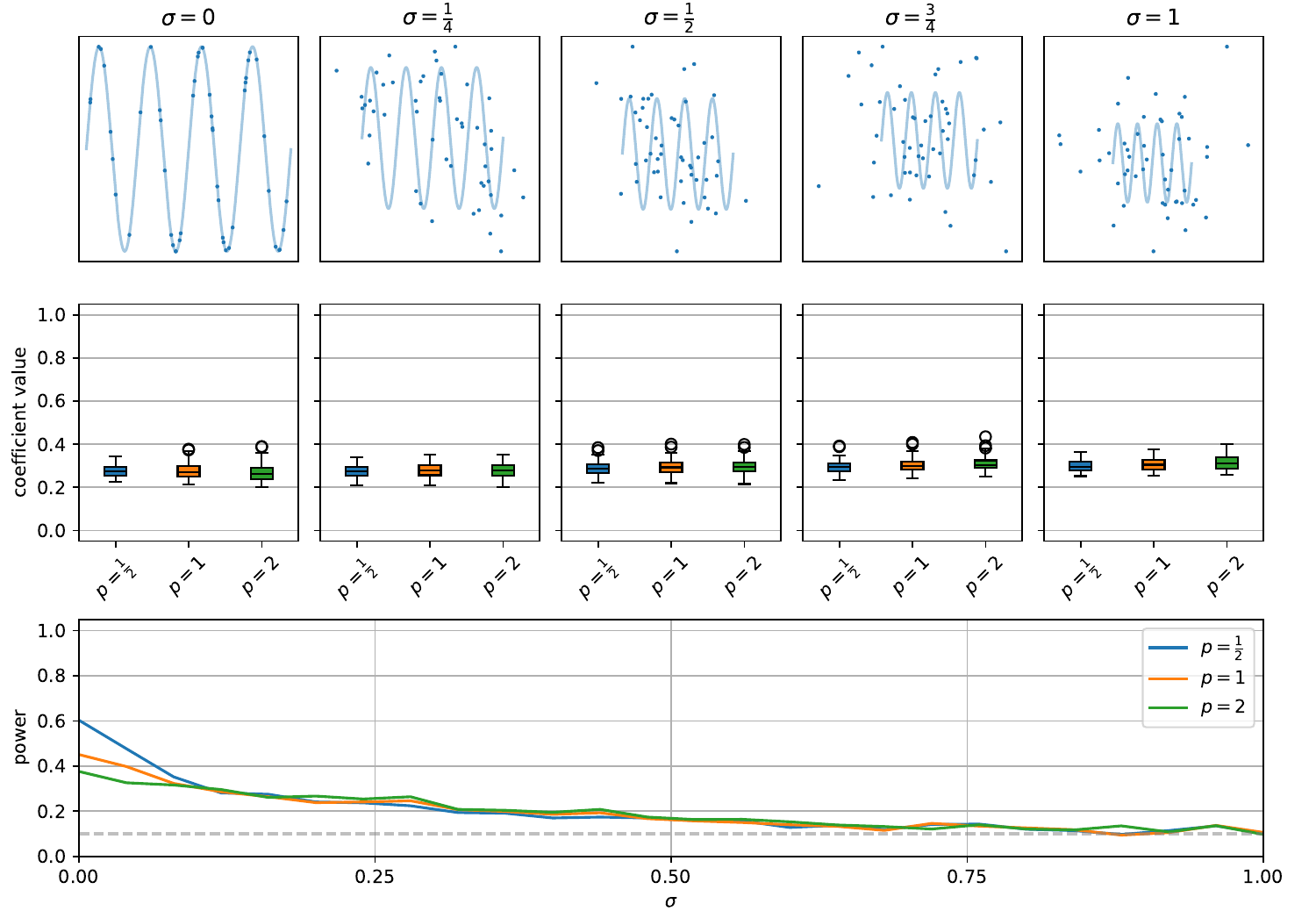}
  \label{fig:different-p27}
  \caption{}
\end{figure}

\begin{figure}
  \centering
  \includegraphics[width=.8\textwidth]{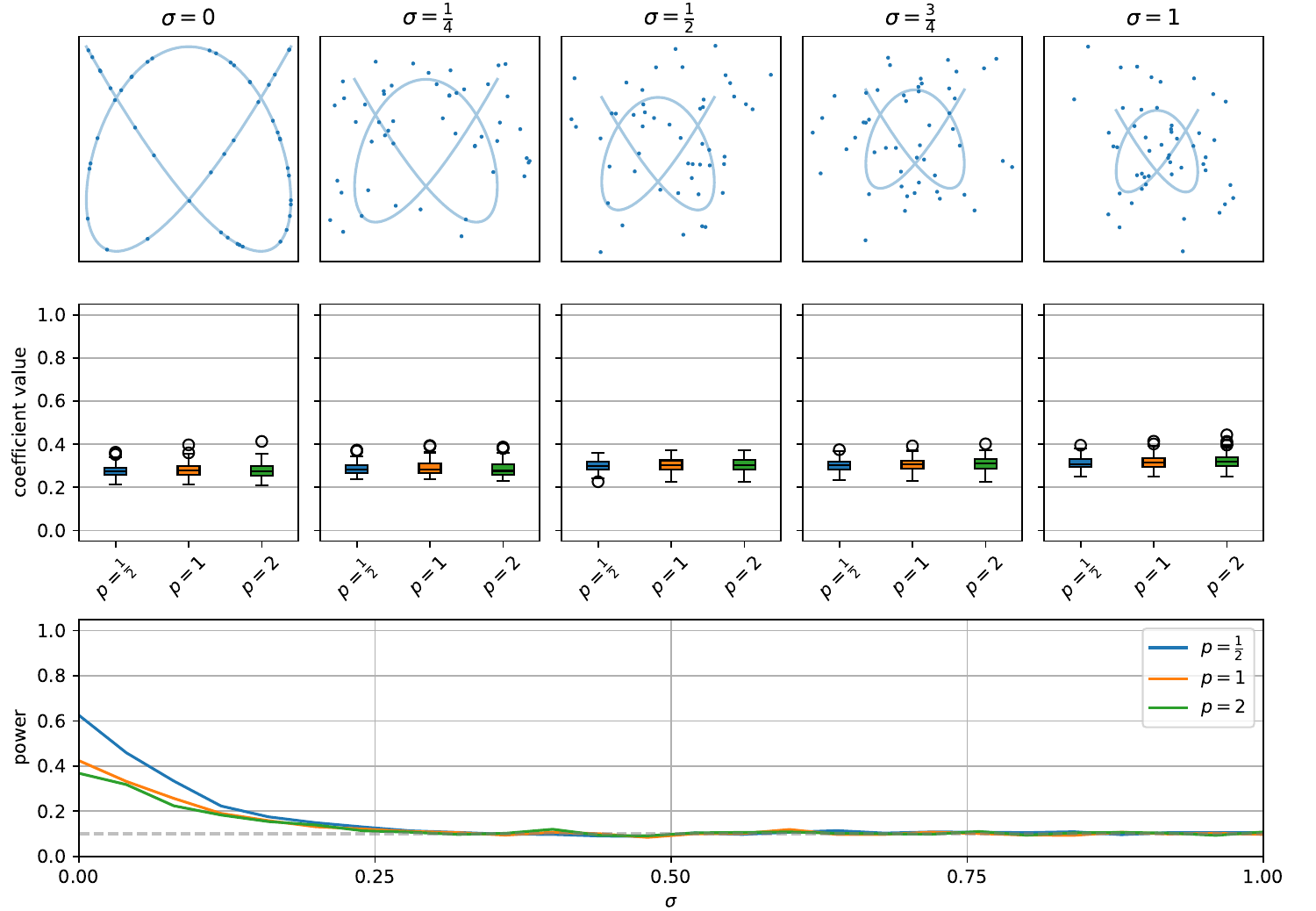}
  \caption{}
  \label{fig:different-p28}
\end{figure}

\paragraph{Influence of $\boldsymbol{p}$.}
Finally, we use the same distributions and noise models as above to study the
effect of the parameter $p$ on the transport correlation coefficient $\tc_*$,
see Definition~\ref{def:isometric-transport-correlation}.
From Figure~\ref{fig:different-p00} to Figure~\ref{fig:different-p28}, we 
compare the values of $\tc_*$ when the parameter $p$ is set equal to $0.5$, $1$,
and $2$, respectively.
The setting is the same as for the previous comparisons, the only
difference being that just $9$ (instead of $29$) random permutations are used
for the permutation tests. One observation that holds for all geometries
$\gamma$ is that the box plots of the estimates of $\td$ are generally very
similar for all $p$.
In contrast, we observe that the parameter $p=0.5$ often provides a (slightly)
higher discriminative power against independence than the choices $p = 1$ or $p
= 2$. 



\end{appendix}

\end{document}